\documentclass[manmat,nospthms]{svjour}

\usepackage{lscape}		

\usepackage{mathptmx}		

\usepackage{amsmath}		
\usepackage{amsthm}		
\usepackage{amsfonts}		
\usepackage{amssymb}		
\usepackage[mathscr]{eucal}	

\theoremstyle{plain}
\newtheorem{mytheorem}{Theorem}[section]
\newtheorem{myproposition}[mytheorem]{Proposition}
\newtheorem{mylemma}[mytheorem]{Lemma}
\newtheorem{mycorollary}[mytheorem]{Corollary}

\theoremstyle{definition}
\newtheorem{mydefinition}[mytheorem]{Definition}

\theoremstyle{remark}

\newtheorem*{myremark*}{Remark}

\renewcommand\ge\geqslant
\renewcommand\le\leqslant

\DeclareMathOperator{\Ord}{Ord}
\DeclareMathOperator{\ord}{ord}
\newcommand{\divides}{\mid}
\def\tfrac#1#2{{\textstyle\frac{#1}{#2}}}

\numberwithin{equation}{section}

\makeatletter
\newenvironment{sizemultline}[1]{%
  \skip@=\baselineskip
  #1%
  \baselineskip=\skip@
  \multline
}{\endmultline \ignorespacesafterend}  
\makeatother

\begin{document}

\author{Robert S. Maier}
\title{Nonlinear differential equations satisfied by\\ certain classical modular forms}
\titlerunning{Nonlinear differential equations for classical modular forms}
\combirunning{Robert S. Maier: Nonlinear differential equations for modular forms}
\institute{Depts.\ of Mathematics and Physics, University of Arizona, Tucson AZ 85721, USA\\\email{rsm@math.arizona.edu}}

\date{}

\maketitle

\begin{abstract}
A unified treatment is given of low-weight modular forms on $\Gamma_0(N),$
$N=2,3,4,$ that have Eisenstein series representations.  For each~$N,$
certain weight-$1$ forms are shown to satisfy a coupled system of nonlinear
differential equations, which yields a single nonlinear third-order
equation, called a generalized Chazy equation.  As byproducts, a table of
divisor function and theta identities is generated by~means of
$q$-expansions, and a transformation law under~$\Gamma_0(4)$ for the second
complete elliptic integral is derived.  More generally, it is shown how
Picard--Fuchs equations of triangle subgroups of~${\it PSL}(2,\mathbf{R}),$
which are hypergeometric equations, yield systems of nonlinear equations
for weight-$1$ forms, and generalized Chazy equations.  Each triangle group
commensurable with~$\Gamma(1)$ is treated.
\end{abstract}

\section{General introduction}
\label{sec:genintro}

In this article we systematically derive ordinary differential equations
(ODEs) that are satisfied by certain elliptic modular forms and their
roots.  The latter are respectively single-valued holomorphic functions,
and potentially multivalued ones, on the upper half plane
$\mathfrak{H}=\{\Im\tau>0\}$.  Among the classical modular groups that will
appear are the full modular group $\Gamma(1)={\it PSL}(2,\mathbf{Z}),$ the
Hecke congruence subgroups $\Gamma_0(N),$ $N=2,3,4,$ and their Fricke
extensions $\Gamma_0^+(N)<{\it PSL}(2,\mathbf{R})$.

There are two distinct sorts of ODE satisfied by forms on classical modular
groups, distinguished by their independent variables.
\begin{enumerate}
  \item If the independent variable is the period
    ratio~$\tau\in\mathfrak{H},$ the ODE will typically be nonlinear.
    Classical examples include Ramanujan's coupled ODEs for Eisenstein
    series on~$\Gamma(1),$ Rankin's fourth-order ODE for the modular
    discriminant~$\Delta$ (a~weight-$12$ form on~$\Gamma(1)$), and Jacobi's
    third-order ODE for the theta-null functions
    $\vartheta_2,\vartheta_3,\vartheta_4$ (which are weight-$\frac12$ forms
    on~$\Gamma_0(4)$).
  \item If the independent variable is a $\mathbf{P}^1(\mathbf{C})$-valued
    Hauptmodul (function field generator) for a modular group~$\Gamma\!,$
    the ODE will be linear~\cite{Stiller88}.  Classical examples include
    Jacobi's hypergeometric ODE for the first complete elliptic
    integral~$\mathsf{K}$ (a~weight-$1$ form on~$\Gamma_0(4)$), viewed as a
    function of the modulus~$k$ or its square~$k^2$ (a~Hauptmodul
    for~$\Gamma_0(4)$); and Picard--Fuchs equations satisfied by periods of
    elliptic families rationally parametrized by other Hauptmoduls.
\end{enumerate}
The stress in this article is on deriving ODEs of the first (nonlinear)
class, which are less closely tied than the second to the classical theory
of special functions.  The computations will make heavy use of
$q$-expansions and the theory of modular forms.  However, special function
methods such as hypergeometric transformations will prove useful in dealing
with forms on the Fricke extensions~$\Gamma_0^+(N)$.

The derivation of ODEs for classical modular forms has recently been
considered by Ohyama (e.g., in~\cite{Ohyama96}) and
Zudilin~\cite{Zudilin2003}.  Our approach differs from theirs by being
extensively `modular,' in~that it exploits dimension formulas (coming
ultimately from the Riemann--Roch theorem), explicit $q$-expansions and
number-theoretic interpretations of their coefficients, etc.  Also, Ohyama
focuses on deriving systems of nonlinear ODEs satisfied by weight-$2$
\emph{quasi-modular} forms, analogous to the Eisenstein series $E_2$
on~$\Gamma(1)$.  His coupled ODEs are of an interesting quadratic type,
called Darboux--Halphen systems~\cite{Ablowitz2006,Harnad2000}.  In the the
present article the fundamental dependent variables are weight-$1$ modular
forms, sometimes multivalued, which are analogous to $E_4^{1/4}\!,\,$
$E_6^{1/6}\!,\,$ and~$\Delta^{1/12}$; and we regard the resulting
differential systems as more fundamental than Darboux--Halphen ones, though
quasi-modular forms play a role.  Zudilin focuses on deriving systems, and
also linear ODEs of the second type distinguished above, which are
satisfied by forms on $\Gamma_0(N),$ $\Gamma_0^+(N),$ or on subgroups
of~$\Gamma_0^+(N)$.  We~are able to clarify the modular underpinnings of
his ODEs, and derive several more such equations.

The weight-$1$ forms studied below include triples of forms on
$\Gamma_0(N),$ $N=2,3,4,$ which we denote
$\mathscr{A}_r,\mathscr{B}_r,\mathscr{C}_r,$ $r=4,3,2$.  They were
introduced as functions on~$\mathfrak{H}$ by the Borweins~\cite{Borwein91}
in their study of alternative AGM (arithmetic-geometric mean) iterations.
Their work was inspired by Ramanujan's theory of elliptic functions to
alternative bases, the base being specified by the `signature'~$r$.
(Cf. Berndt et al.~\cite{Berndt95}.)  Our approach places these functions
firmly in a modular setting (see also~\cite{Maier12}).  As a byproduct of
the analysis of these modular forms and their powers, we derive many
divisor function and theta identities.  A~minor example is Jacobi's Six
Squares Theorem; our proof of~it may be the most explicitly modular one
to~date.  (See Thm.~\ref{thm:6squares}.)  We also give a modular
interpretation of the \emph{second} complete elliptic
integral~$\mathsf{E},$ probably for the first time, by identifying it as a
weight-$1$ form on~$\Gamma_0(4)$ with an explicit, quasi-modular
transformation law.  (See Prop.~\ref{prop:e}.)

The fundamental goal of this article, however, is the development of a
modular theory of `nonlinear' special functions, by determining which
integrable nonlinear ODEs (\emph{generalized Chazy equations}, in our
terminology) can arise in certain well-specified modular contexts.  (See
Thms.\ \ref{thm:main} and~\ref{thm:triangular}; and the discussion
in~\S\,\ref{subsec:final4}.)

\section{Motivation and the first theorem}
\label{sec:intro}

As initial motivation, consider forms on the full modular group
$\Gamma(1)={\it PSL}(2,\mathbf{Z})$ and their differential relations.  Van
der Pol~\cite[\S\,13]{VanderPol51} and Rankin~\cite{Rankin56} proved that
the modular discriminant function~$\Delta$ on $\mathfrak{H}=\{\Im\tau>0\},$
which when viewed as a function of $q:=\exp(2\pi{\rm i}\tau),$
$\left|q\right|<1,$ is defined by
\begin{equation}
\Delta(q) = q\prod_{n=1}^\infty(1-q^n)^{24},
\end{equation}
satisfies the nonlinear, fourth-order homogeneous differential equation
\begin{equation}
\label{eq:Deltaeq}
2\,\Delta^3\Delta'''' - 10\,\Delta^2\Delta'\Delta''' - 3\,\Delta^2 \Delta''^2
+24\,\Delta\Delta'^2 \Delta''- 13\,\Delta'^4=0,
\end{equation}
where ${}'$ signifies the derivation $q\,{\rm d}/{\rm d}q = (2\pi{\rm
i})^{-1}{\rm d}/{\rm d}\tau$.  (The derivation ${\rm d}/{\rm d}\tau$ will
be indicated by a dot; the primes in~(\ref{eq:Deltaeq}) can be optionally
replaced by dots.)  This ODE is fairly well known, as is Jacobi's
third-order one for his theta-null functions $\vartheta_i,$
$i=2,3,4$~\cite{Jacobi1848}.  (For remarks on the latter,
see~\cite{Ehrenpreis94}.)  In the following, the parallels between them
will be brought out.

One approach to understanding the rather complicated Eq.~(\ref{eq:Deltaeq})
is to treat it as a corollary of a much nicer nonlinear third-order
ODE~\cite{Rankin76,VanderPol51}, namely
\begin{subequations}
\begin{gather}
\label{eq:chazy}
\dddot u - 12\,u\ddot u + 18\,{\dot u}^2 = 0,\\
\intertext{i.e.,}
2\,E_2''' - 2\,E_2E_2'' + 3\,{E_2'}^2 = 0.
\label{eq:chazyb}
\end{gather}
\end{subequations}
Here, $u=(2\pi{\rm i}/12)E_2 = \pi{\rm i}\,E_2/6,$ and $E_2$~is the second
(normalized) Eisenstein series on the full modular group.  The Eisenstein
series $E_k=E_k^{{\bf1},{\bf1}}$ on~$\Gamma(1)$ are
\begin{gather*}
  E_{k}(q)=1+a_k\sum_{n=1}^\infty \sigma_{k-1}(n)\,q^n
  =1+a_k\sum_{n=1}^\infty \frac{n^{k-1}q^n}{1-q^n},\\
\sigma_k(n)=\sum_{d\divides n} d^k,\qquad\qquad a_k=\frac2{\zeta(1-k)}=\frac2{L(1-k,{\bf1})}=-\,\frac{2k}{B_k},
\end{gather*}
where $B_k$~is the $k$th Bernoulli number; so $a_2,a_4,a_6,\dots$ are
$-24,240,-504,\dots$.  The nonlinear ODE~(\ref{eq:chazy}) is a so-called
Chazy equation, with the interesting analytic property of having solutions
with a natural boundary (e.g., $\Im\tau=0$ or~$\left|q\right|=1$), beyond
which they cannot be continued; much as is the case with a lacunary series.
(See \cite[pp.~342--3]{Ablowitz91} and~\cite{Chazy11}.)  Substituting $E_2
= \Delta'/\Delta$ into~(\ref{eq:chazyb}) yields~(\ref{eq:Deltaeq}).

Equation~(\ref{eq:chazyb}), in~turn, follows from a result of Ramanujan.
He introduced functions $P,Q,R$ on the disk $\left|q\right|<1,$ defined by
convergent $q$-series, which are identical to $E_2,E_4,E_6$.  That is, they
are respectively a quasi-modular form of weight $2$ and depth~$\le1,$ and
modular forms of weights $4$ and~$6$.  He determined the differential
structure on the ring $\mathbf{C}[E_2,E_4,E_6]$ by showing that
\begin{subequations}
\label{eq:PQR}
\begin{align}
  ({E_4}^3)' &= E_2\cdot {E_4}^3 - {E_4}^2{E_6},\label{eq:PQRa}\\
  ({E_6}^2)' &= E_2\cdot {E_6}^2 - {E_4}^2{E_6},\label{eq:PQRb}\\
  \Delta'  &= E_2\cdot\Delta,\label{eq:PQRc}\\
  12\,{E_2}'   &= E_2\cdot E_2 - E_4,\label{eq:PQRd}
\end{align}
\end{subequations}
where Eqs.~(\ref{eq:PQR}abc) are linearly dependent, since ${E_4}^3 =
{E_6}^2 + 12^3\Delta,$ which is an equality between weight-$12$ modular
forms.  By rewriting the system (\ref{eq:PQR}abd) into a single third-order
equation for~$E_2,$ one obtains Eq.~(\ref{eq:chazyb}).  It is worth noting
for later use that the system (\ref{eq:PQR}abcd) can be rewritten as
\begin{subequations}
\label{eq:PQR2}
  \begin{align}
    ({\mathcal{A}}^{12})' &= \mathcal{E}\cdot {\mathcal{A}}^{12} -
    {\mathcal{A}}^{8} {\mathcal{B}}^{6},\\
    ({\mathcal{B}}^{12})' &= \mathcal{E}\cdot {\mathcal{B}}^{12} -
    {\mathcal{A}}^{8} {\mathcal{B}}^{6},\\
    ({\mathcal{C}}^{12})' &= \mathcal{E}\cdot {\mathcal{C}}^{12},\\
    12\,\mathcal{E}' &= \mathcal{E}\cdot\mathcal{E} - {\mathcal{A}}^4,
  \end{align}
\end{subequations}
where $\mathcal{A}\!,\mathcal{B},\mathcal{C};\mathcal{E}$ are respectively
$E_4^{1/4}\!,\,E_6^{1/6}\!,\,(12^3\Delta)^{1/12};\allowbreak E_2$.  Of
these, $\mathcal{A}\!,\mathcal{B},\mathcal{C}$ are formally weight-$1$
forms for~$\Gamma(1)$; but the first two are multivalued on~$\mathfrak{H}$.
(Their $q$-expansions, the integer coefficients of which lack an
arithmetical interpretation, do~not converge on all of~$\left|q\right|<1$.)

Some of Ramanujan's results along this line were subsequently extended by
Ramamani (\cite{Ramamani70}; see also~\cite{Ramamani89}).  She introduced
three $q$-series somewhat similar to his $P,Q,R,$ and derived a coupled
system of first-order ODEs that they satisfy.  Recently, Ablowitz,
Chakravarty and Hahn (\cite{Ablowitz2006}; see also~\cite{Hahn2008}) showed
that her $q$-series define modular forms on the Hecke subgroup
$\Gamma_0(2)<\Gamma(1),$ including a weight-$2$ quasi-modular form
analogous to~$E_2,$ and derived a single nonlinear third-order ODE that it
satisfies.  This turns~out to be a Chazy-like equation, of a general type
first studied by Bureau~\cite{Bureau87}.

One may wonder whether these results can be generalized, by extending them
to other modular subgroups.  The question is answered in the affirmative by
Theorem~\ref{thm:main} below, which provides a unified treatment of certain
Eisenstein series on the subgroups $\Gamma_0(2),\Gamma_0(3),\Gamma_0(4)$.
For the latter two as~well as for~$\Gamma_0(2),$ a nonlinear third-order
ODE is satisfied by a quasi-modular form of weight~$2$.  A unified
treatment is facilitated by the fact that up~to isomorphism, these are the
only genus-zero proper subgroups of~$\Gamma(1)$ that have exactly three
inequivalent fixed points on~$\mathfrak{H}^*$; with the exception of the
principal modular subgroup~$\Gamma(2),$ which is conjugated
to~$\Gamma_0(4)$ by the $2$-isogeny $\tau\mapsto2\tau$ in~${\it
PSL}(2,\mathbf{R})$; and also the index-$2$ subgroup~$\Gamma^2\!,\,$ which
is a bit anomalous.  The statement of the theorem requires
\begin{mydefinition}
\label{def:u468}
If $u$ is a holomorphic function on~$\mathfrak{H},$ define functions
$u_4,u_6,u_8,\dotsc$ by $u_4:=\dot u-u^2$ and $u_{k+2}:=\dot u_k-kuu_k$.
Thus,
\begin{align*}
  u_4&= \dot u-u^2,\\
  u_6&= \ddot u-6\,u\dot u + 4\,u^3,\\
  u_8&= \dddot u-12\,u\ddot u-6\,\dot u^2 + 48\,u^2\dot u - 24\,u^4.
\end{align*}
A \emph{generalized Chazy equation} ${C}_p$ for~$u$ is a differential
equation of the form $p=0,$ where $p\in\mathbf{C}[u_4,u_6,u_8]$ is a
nonzero polynomial, homogeneous in that the weights of its monomials are
equal.  Here, the weight of $u_4^au_6^bu_8^c$ is $4a+6b+8c$.
\end{mydefinition}
\begin{myremark*}
  The classical Chazy equation, Eq.~(\ref{eq:chazy}), has $p=u_8+24u_4^2$.
  The so-called Chazy--XII class~\cite{Chazy11} includes equations $C_p$
  with $p=u_8+\text{const}\cdot u_4^2$.  This further generalization is
  prefigured by the treatment of Clarkson and Olver~\cite{Clarkson96}.
\end{myremark*}
\begin{mydefinition}
  For any $\chi\colon\mathbf{Z}/N\mathbf{Z}\to\mathbf{C},$ define the
  $\chi$-weighted divisor and conjugate divisor functions
  \begin{equation*}
  \sigma_k(n;\chi) = \sum_{d|n}\chi(d\bmod N)d^k,\qquad
  \sigma^{\mathrm{c}}_k(n;\chi) = \sum_{d|n}\chi({(n/d)}\bmod{N})d^k. 
  \end{equation*}
Such weighted divisor functions, with $\chi$~not necessarily a Dirichlet
character, have been considered by Glaisher~\cite{Glaisher1885b},
Fine~\cite[\S\S\,32 and~33]{Fine88}, and others.  The argument~$\chi$ will
usually be written~out in~full, as $\chi(0),\dots,\chi(N-1)$.
\end{mydefinition}

Results attached to
$\Gamma_0(2),\allowbreak\Gamma_0(3),\allowbreak\Gamma_0(4)$ will be
referred~to as belonging to Ramanujan's theories of signature $4,3,2,$
respectively.  The fixed points on~$\mathfrak{H}^*$ of each group include
(the equivalence classes~of) two cusps, namely the infinite cusp $\tau={\rm
i}\infty$ (i.e.,~$q=0$) and the cusp $\tau=0$; and also a third fixed
point, which for~$\Gamma_0(2)$ is the quadratic elliptic point $\tau=\rm
i,$ for~$\Gamma_0(3)$ is the cubic elliptic point $\tau=\zeta_3 :=
\exp(2\pi{\rm i}/3),$ and for~$\Gamma_0(4)$ is an additional cusp, namely
$\tau=1/2$.  (For a review of these facts, and for triangular fundamental
domains the vertices of~which are these fixed points, see,
e.g.,~\cite{Schoeneberg74}.)

\begin{mytheorem}
\label{thm:main}
    On each modular subgroup\/ $\Gamma_0(N),$ $N=2,3,4,$ i.e., for each of
    the corresponding signatures\/ $r=4,3,2,$ the following are true.
    \begin{enumerate}
    \item
      There is a quasi-modular form\/~$\mathscr{E}_r$ of weight\/ $2$ and
      depth\/ $\le1,$ equaling unity at the infinite cusp, such that\/
      $u=(2\pi{\rm i}/r)\mathscr{E}_r$ satisfies a generalized Chazy
      equation\/ ${C}_{p_r},$ for some polynomial\/~$p_r$.  Namely,
      \begin{alignat*}{2}
	\mathscr{E}_4(q) &= \tfrac13\bigl[4\,E_2(q^2)-E_2(q)\bigr]\\
	                 &= 1 + 8\sum_{n=1}^\infty \sigma_1(n;-1,1)q^n& &= 1 + 8\sum_{n=1}^\infty \sigma^{\mathrm{c}}_1(n;-3,1)q^n, \\
	\mathscr{E}_3(q) &= \tfrac18\bigl[9\,E_2(q^3)-E_2(q)\bigr]\\
	                 &= 1 + 3\sum_{n=1}^\infty \sigma_1(n;-2,1,1)q^n& &=1 + 3\sum_{n=1}^\infty \sigma^{\mathrm{c}}_1(n;-8,1,1)q^n,\\
	\mathscr{E}_2(q) &= \tfrac13\bigl[4\,E_2(q^4)-E_2(q^2)\bigr]\\
	                 &= 1 + 4\sum_{n=1}^\infty \sigma_1(n;-1,0,1,0)q^n& &= 1 + 8\sum_{n=1}^\infty \sigma^{\mathrm{c}}_1(n;-3,0,1,0)q^n,
      \end{alignat*}
      so that\/ $\mathscr{E}_2(q)=\mathscr{E}_4(q^2)$.  The 
      polynomials\/~$p_r\in\mathbf{C}[u_4,u_6,u_8]$ are
      \begin{align}
	p_4 &= u_4u_8 - u_6^2 + 8\,u_4^3,\label{eq:oldchazy1}\\
	p_3 &= u_4u_8^2 - u_6^2u_8 + 24\,u_4^3u_8 -15\,u_4^2u_6^2 + 144\,u_4^5,\label{eq:oldchazy2}\\
	p_2 &= u_4u_8 - u_6^2 + 8\,u_4^3,\label{eq:oldchazy3}
      \end{align}
      so that\/ $p_2 = p_4$.
    \item
      There is a triple of weight\/-$1$ modular forms\/
      $\mathscr{A}_r,\mathscr{B}_r,\mathscr{C}_r$ (allowed to have
      nontrivial [i.e., non-Dirichlet] multiplier systems, and also allowed
      to be multivalued in the above sense of being roots of conventional
      [single-valued] modular forms), such that
      \begin{enumerate}
      \item 
	${\mathscr{A}_r}^r = {\mathscr{B}_r}^r +
	{\mathscr{C}_r}^r,$ each term being a \emph{single-valued} weight\/-$r$ form.
      \item 
	$\mathscr{A}_r,\mathscr{B}_r,\mathscr{C}_r$ vanish respectively at
	(the equivalence classes of) the abovementioned third fixed point,
	the cusp\/ $\tau=0,$ and the cusp\/ $\tau={\rm i}\infty$; and they
	vanish nowhere else.  In each case, the order of vanishing
	(computed with respect to a local parameter for\/ $\Gamma_0(N)$)
	is\/~$1/r$.
      \item 
	${\mathscr{A}_r}^r,{\mathscr{B}_r}^r,{\mathscr{C}_r}^r,$ together
	with\/~$\mathscr{E}_r,$ satisfy the coupled system of nonlinear
	first-order equations
	\begin{align*}
	  ({\mathscr{A}_r}^r)' &= \mathscr{E}_r\cdot {\mathscr{A}_r}^r - {\mathscr{A}_r}^2{\mathscr{B}_r}^r,\\
	  ({\mathscr{B}_r}^r)' &= \mathscr{E}_r\cdot {\mathscr{B}_r}^r - {\mathscr{A}_r}^2{\mathscr{B}_r}^r,\\
	  ({\mathscr{C}_r}^r)' &= \mathscr{E}_r\cdot {\mathscr{C}_r}^r,\\
	  r\,\mathscr{E}_r'&= \mathscr{E}_r\cdot \mathscr{E}_r - {\mathscr{A}_r}\!{\vphantom{\mathscr{A}_r}}^{4-r}{\mathscr{B}_r}^r,
	\end{align*}
	from which the generalized Chazy equation\/ ${C}_{p_r}$
	for\/~$\mathscr{E}_r$ can be derived by elimination.  (The third
	equation says that\/ $u=\dot{\mathscr{C}}_r/\mathscr{C}_r$.)
      \end{enumerate}
  \end{enumerate}
\end{mytheorem}

\begin{myremark*}
  The results of van~der Pol--Rankin and Ramanujan, attached to
  $\Gamma(1),$ cannot be subsumed into Thm.~\ref{thm:main}; but see the
  more general Theorem~\ref{thm:triangular} below.
\end{myremark*}

\begin{myremark*}
  For the subgroup $\Gamma_0(2),$ i.e., when $r=4,$ the coupled ODEs of
  Theorem~\ref{thm:main}(2) are equivalent to those of
  Ramamani~\cite{Ramamani70}, Ablowitz et~al.~\cite{Ablowitz2006}, and
  Hahn~\cite{Hahn2008}.  (Their $\mathcal{P},e\text{
  [or~$\widetilde{\mathcal{P}}$]},\mathcal{Q}$ are the
  $\mathscr{E}_4,{{\mathscr{A}}_4}^2,{{\mathscr{B}}_4}^4$ of the theorem.)
  The nonlinear third-order ODE of Jacobi~\cite{Jacobi1848}, which is
  satisfied by his theta-null functions
  $\vartheta_2,\vartheta_3,\vartheta_4$ on~$\mathfrak{H},$ turns~out to be
  a corollary of the $r=2$ case of the theorem, since
  $\mathscr{A}_2,\mathscr{B}_2,\mathscr{C}_2$ can be chosen to equal
  ${\vartheta_3}^2\!,\,{\vartheta_4}^2\!,\,{\vartheta_2}^2$.
\end{myremark*}

\medskip
The body of this article is laid out as follows.  In \S\,\ref{sec:prelims},
the modular forms $\mathscr{A}_r,\mathscr{B}_r,\allowbreak\mathscr{C}_r$
are defined as eta products and $q$-series.  (These functions
on~$\left|q\right|<1$ were introduced by the Borweins~\cite{Borwein91} as
the theta functions of certain quadratic forms; see the Appendix.  They
play a role in Ramanujan's alternative theories of elliptic
functions~\cite{Berndt95}.  In~\cite{Maier12}, we interpreted them as forms
on $\Gamma_0(2),\allowbreak\Gamma_0(3),\allowbreak\Gamma_0(4)$.)  In
passing, we generate a table of $q$-expansions and divisor-function
identities (Table~\ref{tab:1}), of independent interest, and give a modular
proof of Jacobi's Six Squares Theorem.  In~\S\,\ref{sec:mainthm2}, we prove
Theorem~\ref{thm:main}(2) by exploiting the dimensionality of spaces of
modular forms, i.e., by applying linear algebra to the graded ring
$\mathbf{C}[\mathscr{A}_r,\mathscr{B}_r,\mathscr{C}_r]$.

Section~\ref{sec:theta} is a digression.  From the $r=2$ system, we derive
an elliptic integral transformation law, and differential relations for
theta-nulls that imply Jacobi's nonlinear third-order ODE\null.  Deriving
interesting identities is facilitated by the quasi-modular
form~$\mathscr{E}_2(q)$ equaling (up~to a transcendental constant factor)
the even function $K(q)E(q),$ i.e., the product of the classical first and
second complete elliptic integrals, viewed as functions of the nome~$q$.
No~satisfactorily `modular' transformation law for $E=E(q)$ has previously
been derived.

In \S\,\ref{sec:mainthm1}, we give a direct proof of the generalized Chazy
equations of Theorem~\ref{thm:main}(1).  They too can be derived by linear
algebra.  (Indeed, for each~$r,$ the functions $u_4,u_6,u_8$ are modular
forms of the specified weight, with trivial multiplier systems;
cf.~\cite[Lemma~5]{Takhtajan92}.)  We give a second proof that is less
explicitly modular, based on results of~\cite{Maier12}.  Each of
$\mathscr{A}_r,\mathscr{B}_r,\mathscr{C}_r$ satisfies a `hypergeometric'
Picard--Fuchs equation, which is a linear second-order ODE with three
singular points, the independent variable of which is a Hauptmodul for the
corresponding group~$\Gamma_0(N)$.  Moreover, $\tau$~is a ratio of
solutions of this equation (cf.~\cite{Ford51}).  These facts make possible
the second proof.  Theorem~\ref{thm:general} is an extension of
Theorem~\ref{thm:main}(1), or equivalently, a general result on solutions
of Gauss hypergeometric equations.  It reveals which generalized Chazy
equations can arise from genus-zero subgroups of ${\it PSL}(2,\mathbf{R})$
with three inequivalent fixed points.

In~\S\,\ref{sec:final}, a comparable extension of Theorem~\ref{thm:main}(2)
is obtained.  Theorem~\ref{thm:triangular}, derived using ODE manipulations
like those of Ohyama \cite{Ohyama96}, presents the system of nonlinear
first-order ODEs, satisfied by a triple of weight-$1$ modular forms
$\mathcal{A}\!,\mathcal{B},\mathcal{C},$ that arises from any specified
triangle subgroup of ${\it PSL}(2,\mathbf{R})$; i.e., from its
Picard--Fuchs equation.  As examples, we treat the nine triangle groups
commensurable with $\Gamma(1)$.  (Generalized Darboux--Halphen systems on
these groups have been obtained by Harnad and McKay~\cite{Harnad2000}.)
The systems we derive in \S\S\,\ref{subsec:final2} and~\ref{subsec:final3}
include a `Type~II' one that subsumes Ramanujan's
system~(\ref{eq:PQR2}abcd) on~$\Gamma(1),$ and also applies to the Fricke
extensions $\Gamma_0^+(N),$ $N=2,3$.  A `Type III' system, associated to
index-$2$ subgroups of these three groups, is derived as~well.

\section{Modular forms and divisor function identities}
\label{sec:prelims}

The modular forms ${\mathscr{A}}_r,{\mathscr{B}}_r,{\mathscr{C}}_r,$
$r=4,3,2,$ of which only $\mathscr{A}_4$ is multivalued on~$\mathfrak{H},$
will be defined here in~terms of the Dedekind eta function, rather than
univariate or multivariate theta functions.  In the Appendix, several of
the original definitions of the Borweins~\cite{Borwein91} are reproduced,
as are AGM identities these forms satisfy.

Being a (single-valued) form has its usual meaning.  On
${\mathfrak{H}}^*=\mathfrak{H}\cup\mathbf{P}^1(\mathbf{Q}),$ i.e.,
$\mathfrak{H}\cup\mathbf{Q}\cup\{{\rm i}\infty\},$ a holomorphic
function~$f$ is modular of integral weight~$k$ on some ${\Gamma<\Gamma(1)}$
if $f(\frac{a\tau+b}{c\tau+d}) =\hat\chi(a,b,c,d)(c\tau+d)^kf(\tau)$ for
all
$\pm\left(\begin{smallmatrix}a&b\\c&d\end{smallmatrix}\right)\in\Gamma$.
Here, $\hat\chi$~is a $\mathbf{C}^\times$-valued multiplier system, with
$\hat\chi(-a,-b,-c,-d)$ equaling $(-1)^k\chi(a,b,c,d)$.  The simplest case,
occurring if $\Gamma<\Gamma_0(N)$ for some~$N,$ is when
$\hat\chi(a,b,c,d)$ equals $\chi(d),$ the extension to~$\mathbf{Z}$ of some
Dirichlet character $\chi\colon(\mathbf{Z}/N\mathbf{Z})^\times
\to\mathbf{C}^\times,$ satisfying $\chi(-1)=(-1)^k$.  By definition,
$\chi(d)=0$ if $(d,N)>1,$ where $(\cdot,\cdot)$ is the g.c.d\null.  The
notation~${\bf1}_N$ for the principal character $\text{mod $N$},$
satisfying ${\bf1}_N(d)=1$ if~$(d,N)=1,$ will be used.  The trivial
character of period~$1$ will be denoted~$\bf1$.

In terms of~$q,$ the Dedekind eta function equals
$q^{1/24}\prod_{n=1}^\infty(1-q)^n$.  On~$\Gamma(1),$ it transforms
as~\cite{Rademacher73}
\begin{equation}
\label{eq:etatransf}
  \eta(\tfrac{a\tau+b}{c\tau+d}) =\left\{
  \begin{array}{ll}
    \left(\frac{d}c\right)\zeta_{24}^{3(1-c)+bd(1-c^2)+c(a+d)}[-{\rm i}(c\tau+d)]^{1/2}\,\eta(\tau), & \quad c{\rm\ odd},\\
    \left(\frac{c}d\right)\zeta_{24}^{3d+ac(1-d^2)+d(b-c)}[-{\rm i}(c\tau+d)]^{1/2}\,\eta(\tau), & \quad d{\rm\ odd},\\
  \end{array}
\right.
\end{equation}
if $c>0,$ where $\zeta_{24}:=\exp(2\pi{\rm i}/24),$ and the Jacobi symbol
is taken to satisfy $\bigl(\tfrac{c}{-d}\bigr) = \bigl(\tfrac{c}{d}\bigr)$.
Fine's notation $[\delta]$ for the function $\tau\mapsto\eta(\delta \tau)$
on~$\mathfrak{H}^*$ will be used, so that, e.g., $\Delta=[1]^{24}$.  At any
cusp $s=\tfrac{a}{d}\in\mathbf{Q}\cup\{{\rm i}\infty\}$ (in~lowest terms,
with $\tfrac10$ signifying~${\rm i}\infty$), the order of vanishing
of~$\eta(\delta\tau),$ denoted $\ord_s([\delta]),$ is given by a well-known
formula stated in Ref.~\cite{Martin97},
\begin{equation}
  \label{eq:newguy1}
  \ord_s([\delta]) = \frac1{24}(\delta,d)^2\!/\delta.
\end{equation}
Here, $\ord_s(\cdot)$~is computed with respect to a local parameter on the
quotient curve $X(1)=\Gamma(1)\setminus{\mathfrak{H}}^*\!,\,$ such as the
Klein--Weber $j$-invariant (which equals ${E_4}^3/\Delta =
12^3{E_4}^3/({E_4}^3-{E_6}^2)$ and is a~Hauptmodul for~$\Gamma(1)$).  As
usual, $\ord_{{\rm i}\infty}(f)$ is the lowest power of~$q$ in the Fourier
expansion of~$f$.

If $f$ is a modular form on~$\Gamma\!,$ its order of vanishing at a cusp
$s\in{\mathfrak{H}}^*\!,\,$ computed with respect to a local parameter
for~$\Gamma$ (i.e., on the quotient curve
$X=\Gamma\setminus{\mathfrak{H}}^*$) is
\begin{equation}
\label{eq:order}
  \Ord_{s,\Gamma}(f):=h_\Gamma(s)\cdot \ord_s(f),
\end{equation}
Here, $h_\Gamma(s)$ is the multiplicity with which the image of~$s$ in~$X$
is mapped to~$X(1),$ i.e., the width of the cusp~$s$.  If
$s\in\mathfrak{H}^*$~is not a cusp but rather a quadratic or cubic elliptic
fixed point of~$\Gamma$ (implying that $s\in\mathfrak{H}$), then
by~definition $s$~will be mapped doubly, resp.\ triply to~$X$.  In this
case,
\begin{equation}
  \label{eq:newguy2}
  \ord_s(f) = (2,\text{ resp.\ }3)\cdot \Ord_{s,\Gamma}(f),
\end{equation}
where $\ord_s(f)$~is the order of vanishing of~$f$ at the point
$s\in\mathfrak{H}$ in the conventional sense of analytic functions.  If
$f$~has no~poles and is single-valued on~$\mathfrak{H},$ i.e., has
no~branch points, then this order must be a non-negative integer.

In the case $\Gamma=\Gamma_0(N),$ the inequivalent cusps~$\tau=\tfrac{a}d$
on~$\mathfrak{H}^*$ may be taken to be the fractions
$\tfrac{a}{d}\in\mathbf{Q}$ with $d\divides N,$ $1\le a\le N,$ and with $a$
reduced modulo $(d,N/d)$ while remaining coprime to~$d$.  (E.g., the cusps
of $\Gamma_0(N)$ would be $\tfrac11,\tfrac12$ if~$N=2$; $\tfrac11,\tfrac13$
if~$N=3$; and $\tfrac11,\tfrac12,\tfrac14$ if~$N=4$.  Note that
$\tfrac11\sim0$ and $\tfrac1N\sim {\rm i}\infty$ under~$\Gamma_0(N)$.)  If
this convention is adhered~to, then each inequivalent cusp~$\tfrac{a}{d}$
will have width $h_{\Gamma_0(N)}(\tfrac{a}{d})=e_{d,N}:=N/d(d,N/d)$.

\begin{mydefinition}
\label{def:abc}
  $\mathscr{A}_r,\mathscr{B}_r,\mathscr{C}_r,$ $r=4,3,2,$ are certain
  functions on~$\mathfrak{H}^*\!,$ defined to have the eta-product
  representations
  \begin{alignat*}{2}
    \mathscr{A}_4 &=(2^6\cdot[2]^{24}+[1]^{24})^{1/4}/\,[1]^2[2]^2, \\
    \mathscr{B}_4 &= [1]^4/\,[2]^2, &\qquad&
    \mathscr{C}_4 = 2^{3/2}\cdot[2]^4/\,[1]^2; \\*[\jot]
    \mathscr{A}_3 &=(3^3\cdot[3]^{12}+[1]^{12})^{1/3}/\,[1][3], \\
    \mathscr{B}_3 &= [1]^3/\,[3], &\qquad&
    \mathscr{C}_3 = 3\cdot[3]^3/\,[1]; \\*[\jot]
    \mathscr{A}_2 &=(2^4\cdot[4]^8+[1]^8)^{1/2}/\,[2]^2, \\
    \mathscr{B}_2 &= [1]^4/\,[2]^2, &\qquad&
    \mathscr{C}_2 = 2^2\cdot[4]^4/\,[2]^2,
  \end{alignat*}
  so that by definition, ${\mathscr{A}_r}^r = {\mathscr{B}_r}^r +
  {\mathscr{C}_r}^r$.  At the infinite cusp (i.e., at~$q=0$), each
  $\mathscr{A}_r$ and~$\mathscr{B}_r$ equals unity, and each
  $\mathscr{C}_r$~vanishes.  The~$\mathscr{A}_r,$ defined as roots of
  single-valued modular forms, are potentially multivalued, but it will be
  shown that $\mathscr{A}_2,\mathscr{A}_3$ are single-valued.  One notes
  that $\mathscr{B}_4 = \mathscr{B}_2$ and $\mathscr{C}_4(q) =
  2^{-1/2}\cdot\mathscr{C}_2(q^{1/2})$.
\end{mydefinition}

\begin{myremark*}
  Connections to theta functions, such as Jacobi's theta-nulls
$\vartheta_2,\vartheta_3,\vartheta_4,$ will be discussed
in~\S\,\ref{sec:theta}.  (Also, see the Appendix.)  For the moment, observe
that by theta identities first proved by Euler, or alternatively by the
Jacobi triple product formula, $\mathscr{A}_2,\mathscr{B}_2,\mathscr{C}_2$
equal ${\vartheta_3}^2\!,\,{\vartheta_4}^2\!,\,{\vartheta_2}^2\!$.
Similarly, ${\mathscr{A}_4}^2={\vartheta_2}^4+{\vartheta_3}^4\!,\,$
$\mathscr{B}_4={\vartheta_4}^2\!,\,$ and
$\mathscr{C}_4=\sqrt{2\mathscr{A}_2\mathscr{C}_2}= 2^{1/2}\cdot
{\vartheta_2}{\vartheta_3}$.
\end{myremark*}

\begin{myproposition}
\label{prop:modularprop}
  $\mathscr{A}_r,\mathscr{B}_r,\mathscr{C}_r,$ $r=4,3,2$ are weight\/-$1$
  modular forms on the subgroups\/ $\Gamma_0(N),$ $N=2,3,4,$ respectively,
  with each being single-valued save for\/~$\mathscr{A}_4,$ the square of
  which is single-valued.  Each has exactly one equivalence class of zeroes
  on $\mathfrak{H}^*\!,\,$ at which its order of vanishing is\/~$1/r$
  (computed with respect to a local parameter for\/~$\Gamma_0(N)$), located
  as stated in Theorem\/~{\rm\ref{thm:main}}.  Under the Fricke involution
  $W_N:\tau\mapsto-1/N\tau$ for\/~$\Gamma_0(N),$
  $\mathscr{B}_r$~and\/~$\mathscr{C}_r$ are interchanged in the sense
  that\/ ${\mathscr{B}_r}^2|{W_N}= -{\mathscr{C}_r}^2,$ and\/
  ${\mathscr{A}_r}^2$ is negated.  There is an alternative, explicitly
  single-valued representation for\/ $\mathscr{A}_2,$ namely\/
  $\mathscr{A}_2 = [2]^{10}/\,[1]^4[4]^4$.
\end{myproposition}
\begin{proof}
  It follows from~$(\ref{eq:newguy1})$ that for each~$r,$ $\ord_{{\rm
  i}\infty}(\mathscr{B}_r)=0,$ $\ord_{{\rm i}\infty}(\mathscr{C}_r)=1/r,$
  and $\ord_0(\mathscr{C}_r)=0$; and for $r=4,3,2,$ that
  $\ord_0(\mathscr{B}_r)=1/9,1/8,1/9$.  Also,
  $\ord_{1/2}(\mathscr{B}_2)=\ord_{1/2}(\mathscr{C}_2)=0$.

  The cusps $\tau=0,{\rm i}\infty$ of $\Gamma_0(2),\Gamma_0(3)$ have widths
  $2,1$ and~$3,1,$ and the cusps $\tau=0,\tfrac12,{\rm i}\infty$
  of~$\Gamma_0(4)$ have widths $4,1,1$.  It follows from~(\ref{eq:order})
  that the order of $\mathscr{B}_r,\mathscr{C}_r$ at each cusp is zero,
  except at $\tau=0,{\rm i}\infty$ respectively, where the big-O order in
  each case equals~$1/r,$ as claimed.

  To prove the claim about the zeroes of~$\mathscr{A}_r,$ note that
  $t_2=2^{12}\cdot[2]^{24}/\,[1]^{24},$
  $t_3=3^{6}\cdot[3]^{12}/\,[1]^{12},$ $t_4=2^{8}\cdot[4]^{8}/\,[1]^{8}$
  are Hauptmoduls for $\Gamma_0(2),\Gamma_0(3),\Gamma_0(4),$ i.e., rational
  parameters for the associated quotient curves $X_0(N)$.  Each vanishes at
  the cusp $\tau={\rm i}\infty$ and has a pole at the cusp~$\tau=0$.
  (See~\cite{Maier12}; the normalization factors are unimportant here.)  By
  construction, $\mathscr{A}_4/\mathscr{B}_4=(1+t_2/2^6)^{1/4}\!,\,$
  $\mathscr{A}_3/\mathscr{B}_3=(1+t_3/3^3)^{1/3}\!,\,$ and
  $\mathscr{A}_2/\mathscr{B}_2=(1+t_4/2^4)^{1/2}\!.\,$ Hence,
  \begin{equation}
  \Ord_{0,\Gamma_0(N)}(\mathscr{A}_r)=\Ord_{0,\Gamma_0(N)}(\mathscr{B}_r)+\Ord_{0,\Gamma_0(N)}(\mathscr{A}_r/\mathscr{B}_r)=1/r-1/r=0,
  \end{equation}
  i.e., each~$\mathscr{A}_r$ must be regular and nonzero at the cusp
  $\tau=0$.  Also, each of these quotients $\mathscr{A}_r/\mathscr{B}_r$ is
  zero at the third fixed point; see~\cite[Table~2]{Maier12}.  It follows
  that $\mathscr{A}_r$~must have big-O order at the third fixed point equal
  to~$1/r$.  The third fixed point is quadratic, resp.\ cubic, for $r=4,$
  resp.\ $r=3$; hence by~(\ref{eq:newguy2}), the small-o order of vanishing
  there will be $2\cdot(1/4)=1/2,$ resp.\ $3\cdot(1/3)=1$.  One concludes
  that $\mathscr{A}_4$ has quadratic branch points on~$\mathfrak{H},$ but
  its square and~$\mathscr{A}_3$ are single-valued.

  The statements about the Fricke involution follow readily from the
  transformation law $\eta(-1/\tau)=(-{\rm i}\tau)^{1/2}\eta(\tau)$ and the
  definitions of~$\mathscr{A}_r,\mathscr{B}_r,\mathscr{C}_r$.  To prove
  that $\mathscr{A}_2 = [2]^{10}/\,[1]^4[4]^4,$ observe that $\mathscr{A}_2
  /\{ [2]^{10}/\,[1]^4[4]^4\}$ has zero order of vanishing at each of the
  three inequivalent cusps of~$\Gamma_0(4)$. 
\end{proof}

\begin{myproposition}\ 
  \begin{enumerate}
  \item On $\Gamma_0(2),$ ${{\mathscr{A}}_4}^2$ and
    ${{\mathscr{A}}_4}^4,\allowbreak{{\mathscr{B}}_4}^4,\allowbreak{{\mathscr{C}}_4}^4$
    have trivial character\/~${\bf 1}_2(d),$ which takes\/
    $d\equiv1\pmod2$ to\/~$1$.
  \item On $\Gamma_0(3),$ ${{\mathscr{A}}_3}$ and
    ${{\mathscr{A}}_3}^3,{{\mathscr{B}}_3}^3,{{\mathscr{C}}_3}^3$ have
    quadratic character\/
    $\chi_{-3}(d):=\bigl(\tfrac{-3}{d}\bigr)=\bigl(\tfrac{d}{3}\bigr),$ which
    takes\/ $d\equiv1,2\pmod 3$ to\/~$1,-1,$ and\/ ${{\mathscr{A}}_3}^2$
    has trivial character\/ ${\bf 1}_3(d),$ which takes\/
    $d\equiv1,2\pmod3$ to\/~$1$.    
  \item On $\Gamma_0(4),$ ${{\mathscr{A}}_2}$ has quadratic character\/
    $\chi_{-4}(d):=\bigl(\tfrac{-4}{d}\bigr),$ which takes\/
    $d\equiv1,3\pmod 4$ to\/~$1,-1,$ and\/
    ${{\mathscr{A}}_2}^2,\allowbreak{{\mathscr{B}}_2}^2,\allowbreak{{\mathscr{C}}_2}^2$
    have trivial character\/~${\bf 1}_4(d),$ which takes\/
    $d\equiv1,3\pmod4$ to\/~$1$.
  \end{enumerate}
\label{prop:chars}
\end{myproposition}

\begin{proof}
  To prove each statement, verify it on a generating set for the specified
  subgroup, using the transformation law~(\ref{eq:etatransf}).  For
  example, $\Gamma_0(3)$ has (minimal) generating set
  $\pm\left(\begin{smallmatrix}1&1\\0&1\end{smallmatrix}\right),$
  $\pm\left(\begin{smallmatrix}1&1\\-3&-2\end{smallmatrix}\right),$ and for
  each of the associated maps $\tau\mapsto\tfrac{a\tau+b}{c\tau+d},$ the
  power of~$\zeta_{24}$ appearing in the transformation law
  for~${\mathscr{B}_3}^3\!,\,$ deduced from~(\ref{eq:etatransf}), is
  consistent with the Dirichlet character~$\chi_{-3}$.  The same is true
  for~${\mathscr{C}_3}^3$; hence for~${\mathscr{A}_3}^3$ as~well, since
  ${\mathscr{A}_3}^3={\mathscr{B}_3}^3+{\mathscr{C}_3}^3$.  Hence, the
  claim involving
  ${\mathscr{A}_3}^3\!,\,{\mathscr{B}_3}^3\!,\,{\mathscr{C}_3}^3$ is
  proved.  Further details are left to the reader.  
\end{proof}

The formulas for $\dim\,{{M}}_k(\Gamma_0(N))$ and
$\dim\,{{S}}_k(\Gamma_0(N)),$ the dimensions of the vector spaces of all
modular forms and of cusp forms on~$\Gamma_0(N)$ of weight~$k,$ with
trivial character, are well known~\cite{Cohen77,Diamond2005}.  For
$\Gamma_0(2),\allowbreak\Gamma_0(3),\allowbreak\Gamma_0(4),$ the spaces
${M}_2,{M}_4$ have dimensions $1,2;\,\allowbreak1,2;\,\allowbreak2,3$
respectively; and there are no cusp forms of weight $2$ or~$4$.  Also,
$\dim\,{M}_6(\Gamma_0(2))=2,$ and there are no~cusp forms of weight~$6$
on~$\Gamma_0(2)$.  Similarly,
${M}_1(\Gamma_0(3),\chi_{-3}),\allowbreak{M}_3(\Gamma_0(3),\chi_{-3}),\allowbreak{M}_5(\Gamma_0(3),\chi_{-3})$
have dimensions $1,2,\allowbreak2,$ and
${M}_1(\Gamma_0(4),\chi_{-4}),\allowbreak{M}_3(\Gamma_0(4),\chi_{-4})$ have
dimensions $1,2,$ cusp forms being absent in all cases.  In the absence of
cusp forms, all modular forms in the preceding spaces are combinations of
Eisenstein series.

\begin{myproposition}
  The following spanning relations hold.
  \begin{enumerate}
  \item $M_2(\Gamma_0(2))=\bigl\langle{{\mathscr{A}}_4}^2\bigr\rangle,$
    $M_4(\Gamma_0(2))=\bigl\langle{{\mathscr{A}}_4}^4\!,\,{{\mathscr{B}}_4}^4\bigr\rangle,$
    $M_6(\Gamma_0(2))=\bigl\langle{{\mathscr{A}}_4}^6\!,\,{{\mathscr{A}}_4}^2{{\mathscr{B}}_4}^4\bigr\rangle$.
  \item $M_1(\Gamma_0(3),\chi_{-3})=\bigl\langle{{\mathscr{A}}_3}\bigr\rangle,$
    $M_2(\Gamma_0(3))=\bigl\langle{{\mathscr{A}}_3}^2\bigr\rangle,$
    $M_3(\Gamma_0(3),\chi_{-3})=\bigl\langle{{\mathscr{A}}_3}^3\!,\,{{\mathscr{B}}_3}^3\bigr\rangle,$
    $M_4(\Gamma_0(3))=\bigl\langle{{\mathscr{A}}_3}^4\!,\,{{\mathscr{A}}_3}{{\mathscr{B}}_3}^3\bigr\rangle,$
    $M_5(\Gamma_0(3),\chi_{-3})=\bigl\langle{{\mathscr{A}}_3}^5\!,\,{{\mathscr{A}}_3}^2{{\mathscr{B}}_3}^3\bigr\rangle$.
  \item $M_1(\Gamma_0(4),\chi_{-4})=\bigl\langle{{\mathscr{A}}_2}\bigr\rangle,$
    $M_2(\Gamma_0(4))=\bigl\langle{{\mathscr{A}}_2}^2\!,\,{{\mathscr{B}}_2}^2\bigr\rangle,$
    \hfil\break
    $M_3(\Gamma_0(4),\chi_{-4})=\bigl\langle{{{\mathscr{A}}_2}^3\!,\,{{\mathscr{A}}_2}{{\mathscr{B}}_2}^2}\bigr\rangle,$
    $M_4(\Gamma_0(4))=\bigl\langle{{\mathscr{A}}_2}^4\!,\,{{\mathscr{B}}_2}^4\!,\,{{\mathscr{A}}_2}^2{{\mathscr{B}}_2}^2\bigr\rangle$.
  \end{enumerate}
\label{prop:spanning}
\end{myproposition}
\begin{proof}
  Immediate, by Proposition~\ref{prop:chars} and dimension considerations.  
\end{proof}

Let $M_{\text{even}}(\Gamma)$ denote the graded ring of even-weight modular
forms on~$\Gamma$.  By exploiting the valence formula one can prove the
following generalization.

\begin{myproposition}\ 
\begin{enumerate}
    \item $M_{\rm{even}}(\Gamma_0(2)) = \mathbf{C}[{{\mathscr{A}}_4}^2\!,\,{{\mathscr{B}_4}^4}]\hphantom{{\mathscr{A}}_3}=\mathbf{C}[{{\mathscr{A}}_4}^2\!,\,{{\mathscr{B}_4}^4}-{{\mathscr{C}_4}^4}]$.
    \item $M_{\rm{even}}(\Gamma_0(3)) = \mathbf{C}[{{\mathscr{A}}_3}^2\!,\,{{\mathscr{A}_3\mathscr{B}_3}^3}]=\mathbf{C}[{{\mathscr{A}}_3}^2\!,\,{\mathscr{A}_3({\mathscr{B}_3}^3-{\mathscr{C}_3}^3)}]$.
    \item $M_{\rm{even}}(\Gamma_0(4)) = \mathbf{C}[{{\mathscr{A}}_2}^2\!,\,{{\mathscr{B}_2}^2}]\hphantom{{\mathscr{A}}_3}=\mathbf{C}[{{\mathscr{A}}_2}^2\!,\,{{\mathscr{B}_2}^2}-{{\mathscr{C}_2}^2}]$.
\end{enumerate}
\end{myproposition}

\smallskip
In the sequel, some standard Eisenstein machinery will be used.
(Cf.~\cite[Thms.\ 4.5.2,\ 4.6.2,\ 4.8.1]{Diamond2005}.)  Let a subgroup
$\Gamma_0(N),$ $N\ge2,$ and an integer weight ${k\ge1}$ be specified.  Let
a Dirichlet character
$\chi\colon(\mathbf{Z}/N\mathbf{Z})^\times\to\mathbf{C}^\times,$ extended
to~$\mathbf{Z},$ satisfying $\chi(-1)=(-1)^k,$ also be given.  The
conductor (primitive period) of~$\chi$ will divide~$N$.  For each pair
$\psi,\phi$ of Dirichlet characters, the conductors $u,v$ of which satisfy
$uv\divides N$ and for which $\psi\phi=\chi$ (the equality being one of
characters mod~$N$), there is an Eisenstein series $E_k^{\psi,\phi}\in
M_k(\Gamma_0(N),\chi),$ namely
\begin{equation}
  E_k^{\psi,\phi}(q) := 
\left\{
\begin{alignedat}{3}
&1 + \tfrac{2}{L(1-k,\phi)}\,\cdot \hat{E}_k({\bf1},\phi),&\qquad&\text{if}&\quad&\psi={\bf1},\\
&{2}\cdot\hat{E}_k({\psi},\phi),&\qquad&\text{if}&\quad&\psi\neq{\bf1},
\end{alignedat}
\right.
\end{equation}
where
\begin{equation}
\begin{split}
\hat{E}_k(\psi,\phi) &:=
  \sum_{n=1}^\infty \Biggl[
  \sum_{0<d|n}\psi(n/d)\phi(d)\,d^{k-1}
  \Biggr]q^n\\
&\hphantom{:}= \sum_{e=1}^\infty \sum_{d=1}^\infty \psi(e)\phi(d)\,d^{k-1}\, q^{ed},
\end{split}
\end{equation}
and the $L$-series value $L(1-k,\phi)$ lies in the extension
of~$\mathbf{Q}$ by the values of~$\phi$.  In the case when
$\chi=\mathbf{1}_N,$ the principal character mod~$N,$ these Eisenstein
series are of the form $E_k^{\psi,\psi^{-1}}\!,$ where $\psi$~ranges over
the characters $\text{mod $N$}$ with conductor~$u,$ subject to $u^2\divides
N$.  The subcase $\psi={\bf 1}$ is special: $E_k^{{\bf1},{\bf1}}$ reduces
to~$E_k,$ the $k$th Eisenstein series on~$\Gamma(1)$.

If $k\ge3,$ the collection of all $E_k^{\psi,\phi}(q^\ell),$ where
$E_k^{\psi,\phi}$ is of the above form and ${0<\ell\divides N/(uv)},$ is a
basis for $M_k(\Gamma_0(N),\chi)/S_k(\Gamma_0(N),\chi)$.  For instance, if
${\chi=\bf1}_N$ then these series are of the form
$E_k^{\psi,{\psi}^{-1}}(q^\ell)$ with $0<\ell\divides N/u^2,$ and are
equinumerous with the cusps of~$\Gamma_0(N),$ of which there are
$\sum_{0<d|N}\varphi((d,N/d))$ in~all.  (Here, $(\cdot,\cdot)$
and~$\varphi(\cdot)$ are the g.c.d.\ and totient functions.)  But
if~${k\le2},$ the preceding statements must be modified.  When $k\le2,$ the
Eisenstein series $E_k^{\psi,\phi}(q^\ell)$ are quasi-modular but
in~general are not modular, so the quotient
$M_k(\Gamma_0(N),\chi)/S_k(\Gamma_0(N),\chi)$ is a proper subspace of their
span.

\begin{table}
\caption{}
\begin{center}
{\small
\begin{tabular}{|l|}
\hline
$\begin{alignedat}{4}
&M_2(\Gamma_0(2)):&\quad&{{\mathscr{A}}_4}^2&\quad&=1+24\sum\sigma_1(n;\,0,1)q^n&\quad&=2E_2(q^2)-E_2(q)\\
\hline
&M_4(\Gamma_0(2)):&\quad&{{\mathscr{A}}_4}^4&\quad&=1+24\sum\sigma_3(n;\,3,2)q^n&\quad&=\tfrac15\bigl[4\,E_4(q^2)+E_4(q)\bigr]\\
&{}               &\quad&{{\mathscr{B}}_4}^4&\quad&=1-16\sum\sigma_3(n;\,-1,1)q^n&\quad&=\tfrac{1}{15}\bigl[16\,E_4(q^2)-E_4(q)\bigr]\\
&{}               &\quad&{{\mathscr{C}}_4}^4&\quad&=8\sum\sigma_3(n;\,7,8)q^n&\quad&=-\,\tfrac{4}{15}\bigl[E_4(q^2)-E_4(q)\bigr]\\
\hline
&M_6(\Gamma_0(2)):&\quad&{{\mathscr{A}}_4}^6&\quad&=1+18\sum\sigma_5(n;\,3,4)q^n&\quad&=\tfrac17\bigl[8\,E_6(q^2)-E_6(q)\bigr]\\
&{}               &\quad&{{\mathscr{A}}_4}^2{{\mathscr{B}}_4}^4&\quad&=1+8\sum\sigma_5(n;\,-1,1)q^n&\quad&=\tfrac1{63}\bigl[64\,E_6(q^2)-E_6(q)\bigr]\\
&{}               &\quad&{{\mathscr{A}}_4}^2{{\mathscr{C}}_4}^4&\quad&=2\sum\sigma_5(n;\,31,32)q^n&\quad&=\tfrac{8}{63}\bigl[E_6(q^2)-E_6(q)\bigr]
\end{alignedat}$
\\
\hline
\hline
$\begin{alignedat}{4}
&M_1(\Gamma_0(3),\chi_{-3}):&\quad&{{\mathscr{A}}_3}&\quad&=1+6\sum\sigma_0(n;\,0,1,-1)q^n&\quad&=E_1^{{\bf1},\chi_{-3}}(q)\\
\hline
&M_2(\Gamma_0(3)):&\quad&{{\mathscr{A}}_3}^2&\quad&=1+12\sum\sigma_1(n;\,0,1,1)q^n&\quad&=\tfrac12\bigl[3\,E_2(q^3)-E_2(q)\bigr]\\
\hline
&M_3(\Gamma_0(3),\chi_{-3}):&\quad&{{\mathscr{A}}_3}^3&\quad&={{\mathscr{B}}_3}^3+{{\mathscr{C}}_3}^3\text{ (see below)}&\quad&=E_3^{{\bf1},\chi_{-3}}(q)+\tfrac{27}{2}\,E_3^{\chi_{-3},{\bf1}}(q)\\
&                        &\quad&{{\mathscr{B}}_3}^3&\quad&=1-9\sum\sigma_2(n;\,0,1,-1)q^n&\quad&=E_3^{{\bf1},\chi_{-3}}(q)\\
&                        &\quad&{{\mathscr{C}}_3}^3&\quad&=27\sum\sigma^{\mathrm{c}}_2(n;\,0,1,-1)q^n&\quad&=\tfrac{27}{2}\,E_3^{\chi_{-3},{\bf1}}(q)\\
\hline
&M_4(\Gamma_0(3)):&\quad&{{\mathscr{A}}_3}^4&\quad&=1+8\sum\sigma_3(n;\,4,3,3)q^n&\quad&=\tfrac1{10}\bigl[9\,E_4(q^3)+E_4(q)\bigr]\\
&                &\quad&{{\mathscr{A}}_3}{{\mathscr{B}}_3}^3&\quad&=1-3\sum\sigma_3(n;\,-2,1,1)q^n&\quad&=\tfrac1{80}\bigl[81\,E_4(q^3)-E_4(q)\bigr]\\
&                &\quad&{{\mathscr{A}}_3}{{\mathscr{C}}_3}^3&\quad&=\sum\sigma_3(n;\,26,27,27)q^n&\quad&=-\,\tfrac9{80}\bigl[E_4(q^3)-E_4(q)\bigr]\\
\hline
&M_5(\Gamma_0(3),\chi_{-3}):&\quad&{{\mathscr{A}}_3}^5&\quad&={{\mathscr{A}}_3}^2{{\mathscr{B}}_3}^3+{{\mathscr{A}}_3}^2{{\mathscr{C}}_3}^3\text{ (see below)}&\quad&=E_5^{{\bf1},\chi_{-3}}(q)+\tfrac{27}{2}\,E_5^{\chi_{-3},{\bf1}}(q)\\
&                        &\quad&{{\mathscr{A}}_3}^2{{\mathscr{B}}_3}^3&\quad&=1+3\sum\sigma_4(n;\,0,1,-1)q^n&\quad&=E_5^{{\bf1},\chi_{-3}}(q)\\
&                        &\quad&{{\mathscr{A}}_3}^2{{\mathscr{C}}_3}^3&\quad&=27\sum\sigma^{\mathrm{c}}_4(n;\,0,1,-1)q^n&\quad&=\tfrac{27}{2}\,E_5^{\chi_{-3},{\bf1}}(q)
\end{alignedat}$
\\
\hline
\hline
$\begin{alignedat}{4}
&M_1(\Gamma_0(4),\chi_{-4}):&\quad&{{\mathscr{A}}_2}&\quad&=1+4\sum\sigma_0(n;\,0,1,0,-1)q^n&\quad&=E_1^{{\bf1},\chi_{-4}}(q)\\
\hline
&M_2(\Gamma_0(4)):&\quad&{{\mathscr{A}}_2}^2&\quad&=1+8\sum\sigma_1(n;\,0,1,1,1)q^n&\quad&=\tfrac13\bigl[4\,E_2(q^4)-E_2(q)\bigr]\\
&                 &\quad&{{\mathscr{B}}_2}^2&\quad&=1-8\sum\sigma_1(n;\,0,1,-2,1)q^n&\quad&=\tfrac13\bigl[8\,E_2(q^4)-6\,E_2(q^2)+E_2(q)\bigr]\\
&                 &\quad&{{\mathscr{C}}_2}^2&\quad&=8\sum\sigma_1(n;\,0,2,-1,2)q^n&\quad&=-\,\tfrac23\bigl[2\,E_2(q^4)-3\,E_2(q^2)+E_2(q)\bigr]\\
\hline
&M_3(\Gamma_0(4),\chi_{-4}):&\quad&{{\mathscr{A}}_2}^3&\quad&={{\mathscr{A}}_2}{{\mathscr{B}}_2}^2+{{\mathscr{A}}_2}{{\mathscr{C}}_2}^2\text{ (see below)}&\quad&=E_3^{{\bf1},\chi_{-4}}(q)+8\,E_3^{\chi_{-4},{\bf1}}(q)\\
&                        &\quad&{{\mathscr{A}}_2}{{\mathscr{B}}_2}^2&\quad&=1-4\sum\sigma_2(n;\,0,1,0,-1)q^n&\quad&=E_3^{{\bf1},\chi_{-4}}(q)\\
&                        &\quad&{{\mathscr{A}}_2}{{\mathscr{C}}_2}^2&\quad&=16\sum\sigma^{\mathrm{c}}_2(n;\,0,1,0,-1)q^n&\quad&=8\,E_3^{\chi_{-4},{\bf1}}(q)\\
\hline
&M_4(\Gamma_0(4)):&\quad&{{\mathscr{A}}_2}^4&\quad&=1+4\sum\sigma_3(n;\,4,4,3,4)q^n&\quad&=\tfrac1{15}\bigl[16\,E_4(q^4)-2\,E_4(q^2)+E_4(q)\bigr]\\
&                 &\quad&{{\mathscr{B}}_2}^4&\quad&=1-16\sum\sigma_3(n;\,-1,1)q^n&\quad&=\tfrac1{15}\bigl[16\,E_4(q^2)-E_4(q)\bigr]\\
&                 &\quad&{{\mathscr{C}}_2}^4&\quad&=4\sum\sigma_3(n;\,7,0,8,0)q^n&\quad&=-\,\tfrac{16}{15}\bigl[E_4(q^4)-E_4(q^2)\bigr]\\
&                 &\quad&{{\mathscr{A}}_2}^2{{\mathscr{B}}_2}^2&\quad&=1-2\sum\sigma_3(n;\,-1,0,1,0)q^n&\quad&=\tfrac{1}{15}\bigl[16\,E_4(q^4)-E_4(q^2)\bigr]\\
&                 &\quad&{{\mathscr{A}}_2}^2{{\mathscr{C}}_2}^2&\quad&=2\sum\sigma_3(n;\,7,8)q^n&\quad&=-\,\tfrac{1}{15}\bigl[E_4(q^2)-E_4(q)\bigr]\\
&                 &\quad&{{\mathscr{B}}_2}^2{{\mathscr{C}}_2}^2&\quad&=2\sum\sigma_3(n;\,-7,8,-9,8)q^n&\quad&=\tfrac{1}{15}\bigl[16\,E_4(q^4)-17\,E_4(q^2)+E_4(q)\bigr]
\end{alignedat}$
\\
\hline
\end{tabular}
\bigskip
}
\end{center}
\label{tab:1}
\end{table}

\begin{myproposition}
 One has the Eisenstein series and divisor-function representations shown
 in Table\/~{\rm\ref{tab:1}}, for monomials
 in\/~$\mathscr{A}_r,\mathscr{B}_r,$ $r=4,3,2,$ with multiplier systems of
 quadratic Dirichlet-character type.  (The ones
 involving\/~${\mathscr{C}}_r$ are included for completeness; they follow
 from ${{\mathscr{A}}_r}^r = {{\mathscr{B}}_r}^r + {{\mathscr{C}}_r}^r$.)
\end{myproposition}

\begin{proof}
  For each of the monomials in Proposition~\ref{prop:spanning}, by working
  out the first few coefficients in its $q$-expansion one determines the
  Eisenstein representation given in the rightmost column, and hence the
  full $q$-expansion.  This is a matter of linear algebra, since
  \begin{alignat*}{2}
    M_k(\Gamma_0(2))&\subseteq \langle E_k(q^2),E_k(q) \rangle,&\qquad& k=2,4,6,\\
    M_k(\Gamma_0(3))&\subseteq \langle E_k(q^3),E_k(q) \rangle,&\qquad& k=2,4,\\
    M_k(\Gamma_0(4))&\subseteq \langle E_k(q^4),E_k(q^2),E_k(q) \rangle,&\qquad& k=2,4,\\*
    M_k(\Gamma_0(3),\chi_{-3})&\subseteq \langle E_k^{{\bf1},\chi_{-3}}\!,E_k^{\chi_{-3},{\bf1}} \rangle,&\qquad& k=1,3,5,\\
    M_k(\Gamma_0(4),\chi_{-4})&\subseteq \langle E_k^{{\bf1},\chi_{-4}}\!,E_k^{\chi_{-4},{\bf1}} \rangle,&\qquad& k=1,3,
  \end{alignat*}
where $\subseteq$ signifies $\subset$ if $k\le2$ and $=$ if $k>2$.  The
$L$-series values $L(1-k,\chi_{-3})$ and $L(1-k,\chi_{-4})$ are computed from
$L(1-k,\phi)=-B_{k,\phi}/k$ and the generalized Bernoulli formula for any
Dirichlet character~$\phi$ to the modulus~$N$~\cite{Leopoldt58},
\begin{displaymath}
  \sum_{k=0}^\infty B_{k,\phi}\frac{x^k}{k!} =
  \frac{x}{e^{Nx}-1}\sum_{a=0}^{N-1}\phi(a)e^{ax}. 
\qedhere
\end{displaymath}
\end{proof}

\begin{myremark*}
  Each $q$-expansion in Table~\ref{tab:1} of a modular form of even
  weight~$k$ can alternatively be written in~terms of a
  $\sigma^{\mathrm{c}}_{k-1}$ conjugate divisor function, rather than a
  $\sigma_{k-1}$ divisor function.  For instance,
  \begin{subequations}
  \label{eq:KZ}
  \begin{alignat}{3}
    {{\mathscr{A}}_4}^2 &= {\vartheta_2}^4+{\vartheta_3}^4& &=1+24\sum_{n=1}^\infty\sigma_1(n;\,0,1)q^n&
    &=1+24\sum_{n=1}^\infty\sigma^{\mathrm{c}}_1(n;\,-1, 1)q^n,\\
    {{\mathscr{C}}_4}^4 &= 4\,{\vartheta_2}^4{\vartheta_3}^4& &=8\sum_{n=1}^\infty\sigma_3(n;\,7,8)q^n&
    &=64\sum_{n=1}^\infty\sigma^{\mathrm{c}}_3(n;\,0, 1)q^n.
  \end{alignat}
  \end{subequations}
  Using~(\ref{eq:KZ}ab), one can check that
  ${{\mathscr{A}}_4}^2,{{\mathscr{C}}_4}^4/64$ are identical to the forms
  $C,D$ used by Kaneko and Koike~\cite{Kaneko2004} as generators of~$M_{\rm
  even}(\Gamma_0(2))$.
\end{myremark*}

\begin{myremark*}
  The modular form
  $E_1^{\bf1,\chi_{-3}}(q)=1+6\sum_{n=1}^\infty\sigma_0(n;0,1,-1)q^n$
  figured in Wiles' proof of the Modularity Theorem; for a sketch,
  see~\cite[Ex.~9.6.4]{Diamond2005}.  Table~\ref{tab:1} reveals that
  this modular form is identical to~$\mathscr{A}_3,$ the Borweins' cubic
  theta function in the spirit of Ramanujan.  This observation may be new.
\end{myremark*}

\begin{myremark*}
  Each representation in Table~\ref{tab:1} can be rewritten as a Lambert
  series identity.  Of~the resulting identities, several were recorded by
  Ramanujan and have been given non-modular proofs by Berndt and
  others~\cite{Berndt99,Berndt2006}.
\end{myremark*}

\begin{myremark*}
  The difficulty in extending Table~\ref{tab:1} to higher-degree monomials
  in the triples $\mathscr{A}_r,\mathscr{B}_r,\mathscr{C}_r,$ i.e., in
  deriving simple expressions for their Fourier coefficients in~terms of
  divisor functions, is of~course that one begins to encounter cusp forms.
  To some extent one can work around this.  For instance, Van der
  Pol~\cite{vanderPol54} expressed the coefficients of
  ${\mathscr{A}_2}^{12},{\mathscr{B}_2}^{12},{\mathscr{C}_2}^{12},$ i.e.,
  ${\vartheta_3}^{24},{\vartheta_4}^{24},{\vartheta_2}^{24},$ with the aid
  of Ramanujan's tau function.  Recently Hahn~\cite[Thm.~2.1]{Hahn2008},
  for each even $k\ge4,$ worked~out the combination of the basis monomials
  $\{{{\mathscr{A}}_4}^{2a}{{\mathscr{B}}_4}^{4b}\!,\,\ 2a+4b=k\}$
  of~$M_k(\Gamma_0(2)),$ i.e., the theta polynomials $\{
  ({\vartheta_2}^4+{\vartheta_3}^4){\vphantom{\vartheta_4}}^{a}{\vartheta_4}^{8b}\!,\,\
  2a+4b=k\},$ which equals
  \begin{displaymath}
    E_k^{{\rm i}\infty}(q):=
    \tfrac1{2^k-1}\bigl[2^kE_k(q^2)-E_k(q)\bigr] =
    1+\frac{2k}{(2^k-1)B_k}\sum_{n=1}^\infty\sigma_{k-1}(n;\,-1,1)q^n.
  \end{displaymath}
  This is a weight-$k$ Eisenstein form on~$\Gamma_0(2)$ which vanishes at
  the cusp~${\tau=0}$ and is nonzero at $\tau={\rm i}\infty,$
  like~${\mathscr{B}_4}^k\!$.  In~effect, her combination of monomials
  (unlike the single monomial ${\mathscr{B}_4}^k$ for even~$k\ge6$) has no
  cusp-form component, and therefore has a $q$-expansion with coefficients
  expressible in~terms of divisor functions.
\end{myremark*}

\begin{table}
\caption{}
\begin{center}
{\small
\begin{tabular}{|l|}
\hline
$\begin{alignedat}{3}
&\text{$\Gamma_0(2),$ $k=1$}:&\quad&{{\mathscr{B}}_4}&\quad&=1-4\sum\sigma_0(n;\,0,1,-2,-1,0,1,2,-1)q^n\\
&                            &\quad&{{\mathscr{C}}_4}&\quad&=2^{3/2} \sum\sigma_0(n;\,0,1,-1,-1,0,1,1,-1)q^{n/4}\\
\hline
&\text{\hphantom{$\Gamma_0(2),$} $k=2$}:&\quad&{{\mathscr{B}}_4}^2&\quad&=1-8\sum\sigma_1(n;\,0,1,-2,1)q^n\\
&                            &\quad&{{\mathscr{C}}_4}^2&\quad&=4 \sum\sigma_1(n;\,0,2,-1,2)q^{n/2}
\end{alignedat}$
\\
\hline
\hline
$\begin{alignedat}{3}
&\text{$\Gamma_0(3),$ $k=1$}:&\quad&{{\mathscr{B}}_3}&\quad&=1-3\sum\sigma_0(n;\,0,1,-1,-3,1,-1,3,1,-1)q^n\\
&                            &\quad&{{\mathscr{C}}_3}&\quad&=3 \sum\sigma_0(n;\,0,1,-1,-1,1,-1,1,1,-1)q^{n/3}
\end{alignedat}$
\\
\hline
\hline
$\begin{alignedat}{3}
&\text{$\Gamma_0(4),$ $k=1$}:&\quad&{{\mathscr{B}}_2}&\quad&=1-4\sum\sigma_0(n;\,0,1,-2,-1,0,1,2,-1)q^n\\
&                            &\quad&{{\mathscr{C}}_2}&\quad&=4 \sum\sigma_0(n;\,0,1,-1,-1,0,1,1,-1)q^{n/2}
\\
\hline
&\text{\hphantom{$\Gamma_0(4),$} $k=3$}:&\quad&{{\mathscr{B}}_2}^3&\quad&=1+2\sum\sigma_2(n;\,0,2,-1,-2,0,2,1,-2)q^n\\
& & & & & \qquad{}-16\sum\sigma^{\mathrm{c}}_2(n;\,0,1,-8,-1,0,1,8,-1)q^n\\
&                            &\quad&{{\mathscr{C}}_2}^3&\quad&= \sum\sigma_2(n;\,0,-4,1,4,0,-4,-1,4)q^{n/2}\\
& & & & & \qquad{}+4\sum\sigma^{\mathrm{c}}_2(n;\,0,1,-4,-1,0,1,4,-1)q^{n/2}\\
\end{alignedat}$
\\
\hline
\end{tabular}
\bigskip
}
\end{center}
\label{tab:2}
\end{table}

\begin{myproposition}
  One has the supplementary\/ $q$-expansions shown in
  Table\/~{\rm\ref{tab:2}}, for certain powers
  of\/~$\mathscr{B}_r,\mathscr{C}_r,$ $r=4,3,2,$ the
  multiplier systems of which are not of Dirichlet-character type.  (In
  each, $k$~denotes the weight.)
\end{myproposition}

\begin{proof}
  Each $q$-expansion comes from an Eisenstein representation computed by
  linear algebra, like those of Table~\ref{tab:1}.  The starting points are
  \begin{align}
    {\mathscr{B}}_3(q),{\mathscr{C}}_3(q^3) \in M_1(\Gamma_0(9),\chi_{-3}),\\
    {\mathscr{B}}_2(q),{\mathscr{C}}_2(q^2) \in M_1(\Gamma_0(8),\chi_{-4}),
  \end{align}
  which follow from the definitions of ${\mathscr{B}}_3,{\mathscr{C}}_3$
  and ${\mathscr{B}}_2,{\mathscr{C}}_2,$ like Proposition~\ref{prop:chars}.
  (The statements about ${\mathscr{C}}_3(q^3),{\mathscr{C}}_2(q^2)$ here
  are equivalent to ${\mathscr{C}}_3\in M_1(\Gamma(3),\chi_{-3})$ and
  ${\mathscr{C}}_2\in M_1(\Gamma_0(4)\cap\Gamma(2),\chi_{-4})$.)  To derive
  the given expansions of ${\mathscr{B}}_4,{\mathscr{C}}_4$ and
  ${{\mathscr{B}}_4}^2\!,{{\mathscr{C}}_4}^2\!,$ one simply uses the facts
  that ${\mathscr{B}}_4={\mathscr{B}}_2$ and ${\mathscr{C}}_4(q)=
  2^{-1/2}\cdot{\mathscr{C}}_2(q^{1/2})$.  (The latter fact incidentally
  implies that ${\mathscr{C}}_4\in M_1(\Gamma(4),\chi_{-4})$.)

  The $k=1$ expansions in Table~\ref{tab:2} have previously been been
  derived by non-modular methods, in~\cite[\S\S\,32 and~33]{Fine88}
  and~\cite{Borwein91}.  The final two expansions, of weight-$3$ forms
  on~$\Gamma_0(4),$ may possibly be classical (since $\mathscr{B}_2 =
  {\vartheta_4}^2$ and $\mathscr{C}_2 = {\vartheta_2}^2$), but are more
  likely to be new.  They come from
  \begin{align}
    {{\mathscr{B}}_2}^3(q) &= -E_3^{{\bf1},\chi_{-4}}(q) +
    2\,E_3^{{\bf1},\chi_{-4}}(q^2) -8\,E_3^{\chi_{-4},{\bf1}}(q) +
    64\,E_3^{\chi_{-4},{\bf1}}(q^2) ,\label{eq:basedon1}\\ 
    {{\mathscr{C}}_2}^3(q^2) &=
    E_3^{{\bf1},\chi_{-4}}(q) - \,E_3^{{\bf1},\chi_{-4}}(q^2)
    +2\,E_3^{\chi_{-4},{\bf1}}(q) -8\,E_3^{\chi_{-4},{\bf1}}(q^2),\label{eq:basedon2}
  \end{align}
  in which the four $E_3$'s span $M_3(\Gamma_0(8),\chi_{-4})$.
\end{proof}

The (formal!)\ weight-$1$ modular
form~$\mathscr{A}_4=\sqrt{{\vartheta_2}^4+{\vartheta_3}^4}$
on~$\Gamma_0(2)$ does not fit into the preceding Eisenstein framework,
since it is multivalued.  This is why $\mathscr{A}_4$ and its odd powers
are not expanded in Table \ref{tab:1} or~\ref{tab:2}.  A~bit of
computation yields
\begin{equation}
  \mathscr{A}_4 = 1 + 12\bigl[q - 5\,q^2 + 64\,q^3 - 917\,q^4 + 14850\,q^5 +\cdots\bigr],
\end{equation}
but there is no obvious arithmetical interpretation of the (integral,
see~\cite{Heninger2006}) coefficients of this $q$-expansion, any more than
there is for the $q$-expansions
\begin{subequations}
\begin{align}
{E_4}^{1/4} &= 1+ 60\bigl[q -81\,q^2 + 11008\,q^3 -1751057\,q^4 +  \cdots\bigr],\\
{E_6}^{1/6} &= 1- 84\bigl[q +243\,q^2 + 78784\,q^3 + 29826307\,q^4 +  \cdots\bigr]
\end{align}
\end{subequations}
of the multivalued weight-$1$ forms ${E_4}^{1/4}\!,\,{E_6}^{1/6}$
on~$\Gamma(1),$ introduced in~\S\,\ref{sec:intro}.  It should be noted that
the form ${\mathscr{A}_4}^2\in M_2(\Gamma_0(2))$ is the theta function of
the $D_4$~lattice.

The divisor-function representations of Tables \ref{tab:1}
and~\ref{tab:2} can be viewed as theta identities; including even the $r=3$
ones, since $\mathscr{A}_3,\mathscr{B}_3,\mathscr{C}_3$ too can be
expressed in~terms of $\vartheta_2,\vartheta_3,\vartheta_4$.
(See~\cite{Borwein91} and~\S\,\ref{sec:theta}, below.)  They imply,
\emph{inter alia},
\begin{mytheorem}
\label{thm:6squares}
  Let\/ $r_{2s}(n),$ $n\ge1,$ resp.\ $t_{2s}(n),$ $n\ge0,$ denote the
  number of ways of representing an integer\/ $n$ as the sum of\/ $2s$
  squares, resp.\ triangles.  (These terms signify\/ $m^2,$ resp.\
  $m(m+1)/2,$ with\/ $m$ ranging over\/ $\mathbf{Z}$.)  Then in~terms of
  divisor and conjugate divisor functions,
  \begin{align*}
    r_2(n) &= 4\,\sigma_0(n;\, 0,1,0,-1),\\
    r_4(n) &= 8\,\sigma_1(n;\,0,1,1,1) = 8\,\sigma^{\mathrm{c}}_1(n;\,-3,1,1,1),\\
    r_6(n) &= 16\,\sigma^{\mathrm{c}}_2(n;\,0,1,0,-1) - 4\,\sigma_2(n;\,0,1,0,-1),\\
    r_8(n) &= 4\,\sigma_3(n;\,4,4,3,4) = 16\,\sigma^{\mathrm{c}}_3(n;\,15,1,-1,1);\\
    \hphantom{r_8(n)} &{\hphantom{=}}& \\
    t_2(n) &= 4\,\sigma_0(4n+1;\,0,1,-1,-1,0,1,1,-1)\\
           &= 4\,\sigma_0(8n+2;\,0,1,0,-1),\\
    t_4(n) &= 8\,\sigma_1(2n+1;\,0,2,-1,2) = 16\,\sigma^{\mathrm{c}}_1(2n+1;\,0,1,-2,1)\\
           &= 16\,\sigma_1(2n+1;\,1) = 16\,\sigma^{\mathrm{c}}_1(2n+1;\,1),\\
    t_6(n) &= \sigma_2(4n+3;\,0,-4,1,4,0,-4,-1,4) + 4\,\sigma^{\mathrm{c}}_2(4n+3;\,0,1,-4,-1,0,1,4,-1)\\
           &= 8\,\sigma_2(4n+3;\,0,-1,0,1),\\
    t_8(n) &= 4\,\sigma_3(2n+2;\,7,0,8,0) = 256\,\sigma^{\mathrm{c}}_3(2n+2;\,0,0,1,0)\\
           &= 32\,\sigma_3(n+1;\,7,8) =256\,\sigma^{\mathrm{c}}_3(n+1;\,0,1).
  \end{align*}
\end{mytheorem}
\begin{proof}
  $\mathscr{A}_2={\vartheta_3}^2$ and $\vartheta_3(q)=\sum_{m\in\mathbf{Z}}
  q^{m^2};$ hence $r_{2s}(n)$ is the coefficient of $q^n$ in the
  $q$-expansion of~${\mathscr{A}_2}^s$.  Similarly,
  $\mathscr{C}_2={\vartheta_2}^2$ and
  $\vartheta_2(q)=q^{1/4}\sum_{m\in\mathbf{Z}} q^{m(m+1)};$ hence
  $t_{2s}(n)$ is the coefficient of $q^{2n}$ in the $q$-expansion
  of~$\bigl(q^{-1/2}{\mathscr{C}_2}\bigr)^s$.

  Each $r_{2s}(n)$ formula is taken directly from Table \ref{tab:1}
  or~\ref{tab:2}, and if possible, rewritten in an alternative form based
  on a conjugate divisor function.  The same is true of the first line of
  each of the $t_{2s}(n)$ formulas.  The second, simpler lines of the
  latter follow by elementary arithmetic arguments.
\end{proof}

Theorem~\ref{thm:6squares} is a restatement of Jacobi's Two, Four, Six, and
Eight Squares Theorems, and the known formulas for
$t_2,t_4,t_6,t_8$~\cite{Ono95}.  But the present modular proof of the
formulas for $r_{6}(n),t_6(n),$ in~particular, is illuminating.  (For the
history of these difficult formulas, see~\cite[p.~80]{Milne2002}.)  The
present proof, unlike previous arithmetical or elliptic ones, makes it
clear for the first time how the two terms in the rather awkward formula
for~$r_6(n)$ come from
$E_3^{\chi_{-4},{\bf1}}\!,\,E_3^{{\bf1},\chi_{-4}}\in
M_3(\Gamma_0(4),\chi_{-4})$.  In contrast, a modular derivation of the
seemingly simple formula for~$t_6(n)$ has already been given by Ono et
al.~\cite{Ono95}; but the present derivation, based on
Eq.~(\ref{eq:basedon2}), reveals its complicated underpinnings.

Difficulties arise in extending any Eisenstein approach to $s>4,$ of
course.  As Rankin~\cite{Rankin65} showed, the power~${\vartheta_3}^{2s}$
(i.e.,~${\mathscr{A}_2}^s$) for each $s>4$ has a nonzero cusp-form
component.

\section{Proof of Theorem~\ref{thm:main}(2)}
\label{sec:mainthm2}

Using the results obtained in the last section, one can derive the
differential systems of Theorem~\ref{thm:main}(2) as an exercise in linear
algebra, as follows.

The definition of quasi-modular form used here is standard.
On~$\mathfrak{H}^*\!,\,$ a~holomorphic function~$f$ is quasi-modular of
weight~$2$ and depth~$\le1$ on a subgroup $\Gamma<\Gamma(1),$ with trivial
multiplier system, if
\begin{equation}
\label{eq:qmlaw}
f\left(\frac{a\tau+b}{c\tau+d}\right) =(c\tau+d)^2 f(\tau) +(s/2\pi{\rm
  i})\,c(c\tau+d),
\end{equation}
for all
$\pm\left(\begin{smallmatrix}a&b\\c&d\end{smallmatrix}\right)\in\Gamma$ and
  some $s\in\mathbf{C}$.  One writes $f\in M_2^{\le 1}(\Gamma)$.  The
  constant~$s$ is called the coefficient of affinity of~$f$.

\begin{mylemma}
  If\/ $\mathscr{F}\in M_k(\Gamma,\hat\chi_{\mathscr{F}})$ and\/
  $\mathscr{G}\in M_\ell(\Gamma,\hat\chi_{\mathscr{G}}),$ i.e.,
  $\mathscr{F},\mathscr{G}$ are modular forms on\/~$\Gamma$ with multiplier
  systems not required to be of Dirichlet-character type, and\/
  $\mathscr{F}$~vanishes only at cusps, then
  \begin{enumerate}
    \item $\mathscr{E}:=\mathscr{F}'/\mathscr{F} \in M_2^{\le1}(\Gamma),$
      and\/ $\mathscr{E}$ has coefficient of affinity\/ $k$.
    \item $k\,\mathscr{E}'-\mathscr{E}\cdot\mathscr{E} \in M_4(\Gamma)$.
    \item $k\,\mathscr{G}'-\ell\,\mathscr{E}\cdot\mathscr{G}\in
    M_{\ell+2}(\Gamma,\hat\chi_{\mathscr{G}})$.
  \end{enumerate}
\label{lem:1}
\end{mylemma}
\begin{proof}
  By differentiation of the transformation laws for
  $\mathscr{F},\mathscr{G}$ and~$\mathscr{E}$.  
\end{proof}

\begin{proof}[Proof of Theorem~\ref{thm:main}(2)]
Given $\mathscr{A}_r,\mathscr{B}_r,\mathscr{C}_r,$ $r=4,3,2,$ as in
Definition~\ref{def:abc}, define the function~$\mathscr{E}_r$ of
Theorem~\ref{thm:main} as~$({\mathscr{C}_r}^r)'/{\mathscr{C}_r}^r$.  By
part~(1) of the lemma, it is quasi-modular of weight~$2$ and depth~$\le1$
on~$\Gamma_0(2),\allowbreak\Gamma_0(3),\allowbreak\Gamma_0(4),$
respectively, with trivial multiplier system.  The space of such
quasi-modular forms is spanned respectively by $E_2(q^2),\allowbreak
E_2(q)$; by $E_2(q^3),E_2(q)$; and by $E_2(q^4),E_2(q^2),E_2(q)$.  By
working~out the first few Fourier coefficients of $\mathscr{C}_r$
and~$({\mathscr{C}_r}^r)'/{\mathscr{C}_r}^r,$ and comparing them with those
of these basis functions, one derives the Eisenstein and divisor-function
representations of~$\mathscr{E}_r$ stated in the theorem.

By part~(2) of the lemma,
$r\,\mathscr{E}_r'-\mathscr{E}_r\cdot\mathscr{E}_r$ must lie in
$M_4(\Gamma_0(2)),\allowbreak M_4(\Gamma_0(3)),\allowbreak
M_4(\Gamma_0(4))$ for $r=4,3,2$.  By working out the first two Fourier
coefficients of $r\,\mathscr{E}_r'-\mathscr{E}_r\cdot\mathscr{E}_r,$ and
comparing them with the $q$-expansions of the basis monomials of these
vector spaces, given in Table~\ref{tab:1}, one proves that in each case
this form equals
$-{\mathscr{A}_r}\!{\vphantom{\mathscr{A}_r}}^{4-r}{\mathscr{B}_r}^r,$ as
claimed.

A single example (the $r=3$ case) will suffice.  By direct computation,
\begin{equation}
  3\,\mathscr{E}_3'-\mathscr{E}_3\cdot\mathscr{E}_3 = -1+3\,q+\dots,
\end{equation}
and according to the table, $M_4(\Gamma_0(3))$ is spanned by
\begin{subequations}
\begin{align}
  {{\mathscr{A}}_3}^4&=1+24\,q+\dots,\\
  {{\mathscr{A}}_3}{{\mathscr{B}}_3}^3&=1-3\,q+\dots.
\end{align}
\end{subequations}
The identification of $3\mathscr{E}_3'-\mathscr{E}_3\cdot\mathscr{E}_3$
with $-{\mathscr{A}}_3{{\mathscr{B}}_3}^3$ is justified by the agreement to
first order in~$q$.

By part~(3) of the lemma,
$({\mathscr{A}_r}^r)'-\mathscr{E}_r\cdot{\mathscr{A}_r}^r$ must lie in the
spaces $M_6(\Gamma_0(2)),\allowbreak M_5(\Gamma_0(3),\chi_{-3}),\allowbreak
M_4(\Gamma_0(4)),$ for $r=4,3,2$.  By expanding in~$q$ again, and comparing
coefficients with the $q$-expansions of the spanning monomials listed in
Table~\ref{tab:1}, one proves that this form equals
$-{\mathscr{A}_r}^2{\mathscr{B}_r}^r,$ as claimed.  The details are
elementary.
\end{proof}

One can derive the generalized Chazy equations of Theorem~\ref{thm:main}(1)
by eliminating $\mathscr{A}_r,\mathscr{B}_r,\mathscr{C}_r$ from the
differential systems satisfied by
$\mathscr{A}_r,\mathscr{B}_r,\mathscr{C}_r;\mathscr{E}_r$.  But the
computations are undesirably lengthy, especially for~$r=3$.  Alternative,
more structured proofs of Theorem~\ref{thm:main}(1) will be given
in~\S\,\ref{sec:mainthm1}.

\section{Elliptic integral and differential theta identities}
\label{sec:theta}

This section is a digression, in which the systems of
Theorem~\ref{thm:main}(2) are employed to derive an elliptic integral
transformation formula and differential identities involving Jacobi's
theta-nulls.  The latter are defined on~$\mathfrak{H}\ni\tau,$ i.e.,
on~$\left|q\right|<1,$ by
\begin{equation*}
  \begin{aligned}
  \vartheta_2(q)&=\sum_{m\in\mathbf{Z}} q^{(m+\tfrac12)^2}\\
  &=2\cdot[4]^2/\,[2],
  \end{aligned}
\qquad
  \begin{aligned}
  \vartheta_3(q)&=\sum_{m\in\mathbf{Z}} q^{m^2}\\
  &=[2]^{5}/\,[1]^2[4]^2,
  \end{aligned}
\qquad
  \begin{aligned}
  \vartheta_4(q)&=\sum_{m\in\mathbf{Z}} (-1)^mq^{m^2}\\
  &=[1]^2/\,[2],
  \end{aligned}
\end{equation*}
the given eta representations following from classical $q$-series
identities.  Each~$\vartheta_i$ is a weight-$\tfrac12$ modular form
on~$\Gamma_0(4)$ with a non-Dirichlet multiplier
system~\cite[\S\,81]{Rademacher73}.  They satisfy ${\vartheta_3}^4 =
{\vartheta_2}^4 + {\vartheta_4}^4$.  As noted in~\S\,\ref{sec:prelims},
$\mathscr{A}_2,\mathscr{B}_2,\mathscr{C}_2$ equal
${\vartheta_3}^2\!,\,{\vartheta_4}^2\!,\,{\vartheta_2}^2\!,\,$ and
moreover~\cite{Borwein91}, e.g., $\mathscr{A}_3(q) =
\vartheta_3(q)\vartheta_3(q^3) + \vartheta_2(q)\vartheta_2(q^3)$.

The theta-nulls $\vartheta_2,\vartheta_3,\vartheta_4$ vanish respectively
at $q=0,-1,1,$ i.e., at the points $\tau={\rm i}\infty,1/2,0,$ which are
the inequivalent cusps of~$\Gamma_0(4)$.  Informally, each~$\vartheta_i$
has a simple zero at the respective cusp, and is nonzero and regular
elsewhere.  This does not mean that in the conventional analytic sense,
$\vartheta_i$~is bounded as either of the other two cusps is approached.
For instance, $\vartheta_3(q)\to\infty$ logarithmically as~$q\to1^-,$ i.e.,
as~$\tau\to0$ along the positive imaginary axis.  Having zero order of
vanishing at a finite cusp does not preclude a logarithmic divergence.

The reader is cautioned that in the classical literature, and in the
applied mathematics literature to this day, the argument of
each~$\vartheta_i$ is taken to be $q_2:=\sqrt{q}=\exp(\pi{\rm i}\tau)$
rather than $q=\exp(2\pi{\rm i}\tau)$.  Using $q_2$ rather than~$q$ is
equivalent to viewing the theta-nulls as modular forms on~$\Gamma(2)$
rather than~$\Gamma_0(4),$ since the two subgroups of~$\Gamma(1)$ are
conjugates under the $2$-isogeny $\tau\mapsto2\tau$ in~${\it
PSL}(2,\mathbf{R})$.  In this article the $\Gamma_0(4)$ convention is
adhered~to.

The following is a brief review of how theta-nulls arise from elliptic
integrals.  Consider the parametric family~$\mathcal{E}$ of elliptic plane
curves $E_\alpha/\mathbf{C}$ defined by the equation $y^2=(1-x^2)(1-\alpha
x^2),$ where $\alpha\in\mathbf{P}^1(\mathbf{C})\setminus\{0,1,\infty\}$.
The (first) complete elliptic integral $\mathsf{K}=\mathsf{K}(\alpha)$ is
defined by
\begin{equation}
  \mathsf{K}(\alpha)=\frac12\int_{0}^1 x^{-1/2}(1-x)^{-1/2}(1-\alpha x)^{-1/2}\,dx,
\end{equation}
which makes sense if $0\le\alpha<1,$ and can be continued to a holomorphic
function on~$\mathbf{P}^1(\mathbf{C})_\alpha,$ slit between $\alpha=1$
and~$\alpha=\infty$ to ensure single-valuedness.  The fundamental periods
of the curve~$E_\alpha$ are proportional to $\mathsf{K}(\alpha),{\rm
  i}\,\mathsf{K}(1-\alpha),$ so its period ratio
$\tau=\tau_1/\tau_2\in\mathfrak{H}$ is ${\rm
  i}\,\mathsf{K}(1-\alpha)/\mathsf{K}(\alpha)$.  Since
$\mathsf{K}(0)=\pi/2,$ it is convenient to normalize by defining
$\hat{\mathsf{K}}=\mathsf{K}/(\pi/2)$.

One can show (e.g., by comparing $q$-series) that if $\mathsf{K}$~is
regarded as a function of the nome $q=\exp(2\pi{\rm i}\tau),$ i.e.,
$\hat{\mathsf{K}}=\hat K(q),$ then $\hat K(q)$ equals~${\vartheta_3}^2(q),$
which is holomorphic and single-valued on~$\mathfrak{H}^*$.  The reason for
this equality is that in modern language, $\mathcal{E}$~is the elliptic
family attached to~$\Gamma_0(4)$.  The parameter~$\alpha$ can also be
viewed as a function of~$q,$ i.e., as a $\Gamma_0(4)$-stable holomorphic
function on~$\mathfrak{H},$ with a zero at~$\tau={\rm i}\infty$ and a pole
at~$\tau=0$: it~is a Hauptmodul for~$\Gamma_0(4)$.

So, $\hat K={\vartheta_3}^2=\mathscr{A}_2$.  By Table~\ref{tab:1}, $\hat
K\in M_1(\Gamma_0(4),\chi_{-4}),$ and $\hat K$~has the Eisenstein series
representation
\begin{equation}
  \label{eq:Kdef}
  \hat K(q) = 1 + 4\sum_{n=1}^\infty \sigma_0(n;\,0,1,0,-1)q^n =
  E_1^{{\bf1},\chi_{-4}}(q).
\end{equation}
This expansion is well known, as is the presence of the
character~$\chi_{-4}$ in the transformation law of~$\hat K$ under
$\tau\mapsto\tfrac{a\tau+b}{c\tau+d}$ with
$\pm\left(\begin{smallmatrix}a&b\\c&d\end{smallmatrix}\right)\in\Gamma_0(4)$.
But, analogous expansions and transformation properties for the {\em
second\/} complete elliptic integral are not.  This function
$\mathsf{E}=\mathsf{E}(\alpha)$ is defined locally (on $\alpha\le0<1$) by
\begin{equation}
  \label{eq:Edef}
  \mathsf{E}(\alpha)=\frac12\int_{0}^1 x^{-1/2}(1-x)^{-1/2}(1-\alpha x)^{1/2}\,dx.
\end{equation}
Since $\mathsf{E}(0)=\pi/2$ also, one normalizes by letting
$\hat{\mathsf{E}}=\mathsf{E}/(\pi/2)$.  As with~$\hat{\mathsf{K}},$
$\hat{\mathsf{E}}$~can be viewed as~$\hat E(q),$ a holomorphic and
single-valued function on~$\mathfrak{H}$.  It is a classical result
(see~\cite[p.~218]{Enneper1890} and~\cite[\S\,31]{Glaisher1885a}) that
\begin{equation}
\begin{split}
  \hat K(q)\hat E(q) &= 1+8\sum_{n=1}^\infty \frac{q^{2n}}{(1+q^{2n})^2}\\
  &= 1 + 4\sum_{n=1}^\infty \sigma_1(n;-1,0,1,0)q^n
  = \tfrac13\bigl[4\,E_2(q^4)-E_2(q^2)\bigr].
\end{split}
\label{eq:split}
\end{equation}
Hence $\hat E=(\hat K\hat E)/\hat K,$ i.e.,
\begin{equation}
\hat E(q) =\frac{1+4\sum_{n=1}^\infty \sigma_1(n;\,-1,0,1,0)q^n}
{1+4\sum_{n=1}^\infty \sigma_0(n;\,0,1,0,-1)q^n}
=\frac{4\,E_2(q^4)-E_2(q^2)}{3\,E_1^{{\bf 1},\chi_{-4}}(q)}
\end{equation}
(cf.\ \cite[\S\,38]{Glaisher1885a}).  Remarkably, the divisor-function
representation of~(\ref{eq:split}) is identical to that of the
quasi-modular form $\mathscr{E}_2\in M_2^{\le1}(\Gamma_0(4)),$ given in
Theorem~\ref{thm:main}.  So,
\begin{equation}
  \mathscr{E}_2 = ({\mathscr{C}_2}^2)'/{\mathscr{C}_2}^2
  = ({\vartheta_2}^4)'/{\vartheta_2}^4 =\hat K \hat E
\end{equation}
and $\hat K \hat E \in M_2^{\le1}(\Gamma_0(4))$.  Also, one can write $\hat
E = {\mathscr{E}}_2 / {\mathscr{A}}_2$.

\begin{myproposition}
\label{prop:e}
  The forms $\hat K,$ $\hat E,$ $\hat K\hat E$ have the transformation
  laws
  \begin{align*}
    \hat K(q_1) &= \chi_{-4}(d)(c\tau+d)\hat K(q),\\
    \hat E(q_1) &= \chi_{-4}(d)\left[(c\tau+d)\hat E(q) + (\pi{\rm i})^{-1}c\hat K(q)^{-1}\right],\\
    \hat K\hat E(q_1) &=(c\tau+d)^2\,\hat K\hat E(q) + (\pi{\rm i})^{-1}c(c\tau+d),
  \end{align*}
for all
$\pm\left(\begin{smallmatrix}a&b\\c&d\end{smallmatrix}\right)\in\Gamma_0(4)$.
Here, $q=\exp(2\pi{\rm i}\tau),$ $q_1=\exp(2\pi{\rm i}\tau_1)$ with\/
$\tau_1=\tfrac{a\tau+b}{c\tau+d}$.
\end{myproposition}

\begin{proof}
  $\hat K={\mathscr{A}_2}\in M_1(\Gamma_0(4),\chi_{-4}),$ hence its law is
  known.  The quasi-modular law for $\hat K\hat
  E=\mathscr{E}_2=({\mathscr{C}_2}^2)'/({\mathscr{C}_2}^2),$ of the
  type~(\ref{eq:qmlaw}), follows from Lemma~\ref{lem:1}(1).  Taking the
  quotient yields the law for~$\hat E$.
\end{proof}

\begin{myremark*}
  This transformation law under $\Gamma_0(4)$ for the (normalized) second
  complete elliptic integral $\hat E(q)$ is arguably the most informative
  obtained to~date.  Tricomi \cite[Chap.~IV, \S\,2]{Tricomi51} has some
  related results, but it is difficult to compare them, since he (i)~used
  homogeneous modular forms, i.e., functions of $\tau_1,\tau_2$ rather
  than~$\tau,$ (ii)~worked in~terms of
  $\mathsf{K}(\alpha),\mathsf{E}(\alpha),$ and especially, (iii)~treated
  only $\left(\begin{smallmatrix}a&b\\c&d\end{smallmatrix}\right)=
    \left(\begin{smallmatrix}1&1\\0&1\end{smallmatrix}\right),
      \left(\begin{smallmatrix}0&-1\\1&0\end{smallmatrix}\right),$ which
        are generators of~$\Gamma(1)$ rather of than $\Gamma_0(4)$
        (or~$\Gamma(2)$).
\end{myremark*}

\begin{myremark*}
  One can similarly define quasi-modular forms $\hat K\hat G,\hat K\hat
  I\in M_2^{\le1}(\Gamma_0(4))$ by
  \begin{subequations}
  \begin{align}
  ({\mathscr{A}_2}^2)'/{\mathscr{A}_2}^2
    &= ({\vartheta_3}^4)'/{\vartheta_3}^4 =:\hat K \hat G,\\
  ({\mathscr{B}_2}^2)'/{\mathscr{B}_2}^2
    &= ({\vartheta_4}^4)'/{\vartheta_4}^4 =:\hat K \hat I
  \end{align}
  \end{subequations}
(cf.\ Glaisher~\cite{Glaisher1885a}), and work~out the transformation laws
  of $\hat G,\hat I$.  Each of $\hat G,\hat I$ has a representation as a
  complete elliptic integral, analogous to~(\ref{eq:Edef}) for~$\hat E$.
\end{myremark*}

\begin{myproposition}
\label{prop:KE}
  The theta-nulls\/ $\vartheta_2,\vartheta_3,\vartheta_4,$ together
  with\/~$\hat K\hat E,$ satisfy a differential system on the
  half-plane\/~$\mathfrak{H},$ namely
  \begin{align*}
    4\,\vartheta_2'/\vartheta_2 &= \hat K\hat E, & 2(\hat K\hat E)'&=(\hat
    K\hat E)^2 - {\vartheta_3}^4{\vartheta_4}^4,\\
    4\,\vartheta_3'/\vartheta_3 &= \hat K\hat E - {\vartheta_4}^4,\\
    4\,\vartheta_4'/\vartheta_4 &= \hat K\hat E - {\vartheta_3}^4,
  \end{align*}
  where\/ ${}'$ signifies\/ $q\,{\rm d}/{\rm d}q = (2\pi{\rm i})^{-1}{\rm
    d}/{\rm d}\tau$.
\end{myproposition}

\begin{proof}
  Substitute
  ${\vartheta_3}^2\!,\,{\vartheta_4}^2\!,\,{\vartheta_2}^2\!,\,\hat K\hat
  E$ for $\mathscr{A}_2,\mathscr{B}_2,\mathscr{C}_2,\mathscr{E}_2$ in the
  $r=2$ system of Theorem~\ref{thm:main}(2).  
\end{proof}

\begin{myremark*}
  This system of coupled ODEs may be new, though it can be deduced from
  identities of Glaisher and of the Borweins~\cite[\S\,2.3]{Borwein87}.
  For $i=2,3,4,$ one can derive from~it a nonlinear third-order ODE
  satisfied by~$\vartheta_i,$ by eliminating the other three dependent
  variables.  For each~$\vartheta_i$ this turns~out to be
  \begin{equation}
    \label{eq:jacobi3}
    (\vartheta^2\vartheta'''-15\,\vartheta\vartheta'\vartheta''+30\,\vartheta'^3)^2
    + 32 (\vartheta\vartheta''-3\,\vartheta'^2)^3 = \vartheta^{10}(\vartheta\vartheta''-3\,\vartheta'^2)^2.
  \end{equation}
  This is the 1847 equation of Jacobi~\cite{Jacobi1848}, which was
  mentioned in~\S\,\ref{sec:intro}.  His derivation used differentiation
  with respect to Hauptmoduls for $\Gamma_0(4)$ (his $k^2$ and~$k'^2$).

  For an easy proof that each of $\vartheta_2,\vartheta_3,\vartheta_4$ must
  satisfy the \emph{same} third-order ODE, reason as follows.  First,
  work~out the differential systems for
  $\vartheta_2,\vartheta_3,\vartheta_4;\hat K\hat G$ and
  $\vartheta_2,\vartheta_3,\vartheta_4;\hat K\hat I$ that are analogues of
  the system for $\vartheta_2,\vartheta_3,\vartheta_4;\hat K\hat E$ in
  Proposition~\ref{prop:KE}.  Then notice that up~to cyclic permutations of
  the ordered pairs $(\vartheta_2,\hat K\hat
  E),\allowbreak(\vartheta_3,\hat K\hat G),\allowbreak(\vartheta_4,\hat
  K\hat I),$ the three systems are the same.  Hence, eliminating all
  dependent variables except a single~$\vartheta_i$ must yield the same
  equation, irrespective of~$i$.

  Brezhnev~\cite[\S\,7]{Brezhnev2006} has recently derived a different but
  related differential system, symmetric and elegant, in which the
  dependent variables are $\vartheta_2,\vartheta_3,\vartheta_4,$ and
  (in~the notation used here) the element $(\hat K\hat E+\hat K\hat G+\hat
  K\hat I)(q)$ of~$M_2^{\le1}(\Gamma_0(4)),$ which by~examination is
  proportional to~$E_2(q^2)$.  His system can be obtained by averaging
  together the three preceding ones; and this averaging ensures symmetry.
\end{myremark*}

\section{Proofs of Theorem~\ref{thm:main}(1); Hypergeometric identities}
\label{sec:mainthm1}

Direct derivations of the generalized Chazy equations of
Theorem~\ref{thm:main}(1) will now be given.  They will not employ, except
superficially, the differential systems satisfied by the weight-$1$ modular
forms $\mathscr{A}_r,\mathscr{B}_r,\mathscr{C}_r$.

Two proofs of Theorem~\ref{thm:main}(1) are supplied.  The first is an
explicitly modular, linear-algebraic one. It is modeled after Resnikoff's
proof~\cite{Resnikoff65} of Eq.~(\ref{eq:jacobi3}), Jacobi's nonlinear
third-order ODE (for~$\vartheta=\vartheta_3$).  Equation~(\ref{eq:chazy}),
the Chazy equation satisfied by $u=(2\pi{\rm i}/12)E_2,$ can be proved
similarly.

The second proof employs analytic manipulations of Picard--Fuchs equations,
and relies on results of~\cite{Maier12}.  It is based on a sort of
nonlinear hypergeometric identity, stated as
Proposition~\ref{prop:special}, which holds for certain very special
parameter values that appear in Picard--Fuchs equations attached
to~$\Gamma_0(N),$ $N=2,3,4$.  Remarkably, this identity has an extension to
\emph{all} parameter values, namely Theorem~\ref{thm:general}.

Rankin~\cite{Rankin76} gives an altogether different sort of proof of the
Chazy equation~(\ref{eq:chazyb}), based on elementary arithmetic methods.
One may speculate that the generalized Chazy equations can also be derived
by such methods.

\subsection{A modular proof of Theorem~\ref{thm:main}(1)}
\label{subsec:51}

Define $\mathscr{A}_r,\mathscr{B}_r,\mathscr{C}_r$ as
in~\S\,\ref{sec:prelims}, and let
$\mathscr{E}_r=({\mathscr{C}_r}^r)'/{\mathscr{C}_r}^r,$ as in the proof of
Theorem~\ref{thm:main}(2).  For $r=4,3,2,$ $\mathscr{E}_r$~is quasi-modular
of weight~$2$ and depth~$\le1$ on~$\Gamma_0(N),$ $N=2,3,4,$ respectively.
For $k=4,6,8,\dotsc,$ define ${\mathfrak{u}}_k^{(r)}$ by
\begin{equation*}
{\mathfrak{u}}_4^{(r)}=
r\mathscr{E}_r'-\mathscr{E}_r\cdot\mathscr{E}_r,
\quad\qquad
{\mathfrak{u}}_{k+2}^{(r)}={\mathfrak{u}}_k^{(r)}{}' - (k/r)\mathscr{E}_r\cdot
{\mathfrak{u}}_k^{(r)}.
\end{equation*}
By Lemma~\ref{lem:1}, ${\mathfrak{u}}_k^{(r)}\in M_k(\Gamma_0(N))$.  If
$u:=(2\pi{\rm i}/r)\mathscr{E}_r$ and $u_k$~is defined in~terms of~$u$ as
in~\S\,\ref{sec:intro}, then one has that $u_k=(2\pi{\rm
i})^{k/2}\,{\mathfrak{u}}_k/r^2$ for all~$k$.

By the last differential equation in Theorem~\ref{thm:main}(2c),
${\mathfrak{u}}_4^{(r)}$ equals
$-{\mathscr{A}_r}\!{\vphantom{\mathscr{A}_r}}^{4-r}{\mathscr{B}_r}^r\!$.
By Theorem~\ref{thm:main}(2b), $\Ord_0(\mathscr{A}_r)=0$ and
$\Ord_{0,\Gamma_0(N)}(\mathscr{B}_r)=1/r$; hence one has that
${\Ord_{0,\Gamma_0(N)}({\mathfrak{u}}_4^{(r)})=1}$.  It is evident that
$\Ord_{0,\Gamma_0(N)}({\mathfrak{u}}_k^{(r)})\ge1$ for~$k\ge4$.

According to the valence formula~\cite[Chap.~V]{Schoeneberg74}, the total
number of zeroes of a nonzero element $f\in M_k(\Gamma_0(N)),$ counted with
respect to local parameters, is equal to $(k/12)[\Gamma(1):\Gamma_0(N)]$.
It follows that if at any $s\in\mathfrak{H}^*\!,\,$ it is the case that
$\Ord_{s,\Gamma_0(N)}(f) > (k/12)[\Gamma(1):\Gamma_0(N)],$ then $f=0$.
Here, the subgroup index $[\Gamma(1):\Gamma_0(N)]$ equals $3,4,6$
when~$N=2,3,4$.

In the following analyses, the superscript ${}^{(r)}$ will be omitted for
readability.

\begin{itemize}
\item $r=4,$ $\Gamma_0(N)=\Gamma_0(2)$.  One sets $k=12,$ i.e., uses linear
algebra on~$M_{12}(\Gamma_0(2))$.  For each
$g\in{\mathfrak{V}}=\{{\mathfrak{u}}_4{\mathfrak{u}}_8,
{{\mathfrak{u}}_6}^2\!,\, {{\mathfrak{u}}_4}^3\},$ it is the case that
$g\in M_{12}(\Gamma_0(2))$ and ${\Ord_{0,\Gamma_0(2)}(g)\ge2}$.  There is a
linear combination~$f$ of the three monomials in~${\mathfrak{V}}$ for which
$\Ord_{0,\Gamma_0(2)}(f)\ge4$.  But if $f\in M_{12}(\Gamma_0(2))$ vanishes
with order greater than $(k/12)[\Gamma(1):\Gamma_0(2)]=3,$ then $f=0$.

This combination can be found by direct computation, using $q$-series (even
though $q$-series are expansions at the infinite cusp, not at~$\tau=0$).
To~$O(q^1),$
\begin{align*}
{\mathfrak{u}}_4{\mathfrak{u}}_8 &= \tfrac32 - 8\,q+\dots,\\
{{\mathfrak{u}}_6}^2 &= 1+16\,q+\dots,\\
{{\mathfrak{u}}_4}^3 &= -1+48\,q+\dots.
\end{align*}
There is a unique combination (up~to scalar multiples) that is zero to
this order, and must therefore vanish identically; namely,
$2{\mathfrak{u}}_4{\mathfrak{u}}_8 - 2{{\mathfrak{u}}_6}^2 + {{\mathfrak{u}}_4}^3$.
Its vanishing is equivalent to $u_4u_8 - u_6^2 + 8u_4^3=0$.
  \item $r=3,$ $\Gamma_0(N)=\Gamma_0(3)$.  One sets $k=20,$ i.e., uses
linear algebra on~$M_{20}(\Gamma_0(3))$.  For each
$g\in{\mathfrak{V}}=\{{\mathfrak{u}}_4{{\mathfrak{u}}_8}^2\!,\,
{{\mathfrak{u}}_6}^2{{\mathfrak{u}}_8},
{{\mathfrak{u}}_4}^3{{\mathfrak{u}}_8},
{{\mathfrak{u}}_4}^2{{\mathfrak{u}}_6}^2\!,\, {{\mathfrak{u}}_4}^5\},$ it
is the case that $g\in M_{20}(\Gamma_0(3))$ and
${\Ord_{0,\Gamma_0(3)}(g)\ge3}$.  There is a linear combination~$f$ of the
five monomials in~${\mathfrak{V}}$ for which $\Ord_{0,\Gamma_0(3)}(g)\ge7$.
But if $f\in M_{20}(\Gamma_0(3))$ vanishes with order greater than
$(k/12)[\Gamma(1):\Gamma_0(3)]=20/3,$ then $f=0$.

As in the $r=4$ case, this combination can be found by a direct computation
(a~tedious one).  To~$O(q^3),$
\begin{align*}
{\mathfrak{u}}_4{{\mathfrak{u}}_8}^2 &= -\tfrac{64}9 -\tfrac{112}3\,q + 23\,q^2 - \tfrac{7123}3\,q^3+\dots,\\
{{\mathfrak{u}}_6}^2{{\mathfrak{u}}_8} &= -\tfrac{128}{27} - \tfrac{368}9\,q -\tfrac{944}3\,q^2 + \tfrac{7381}{9}\,q^3+\dots,\\
{{\mathfrak{u}}_4}^3{{\mathfrak{u}}_8} &= \tfrac83 - 13\,q - 201\,q^2 + 2075\,q^3+\dots,\\
{{\mathfrak{u}}_4}^2{{\mathfrak{u}}_6}^2 &= \tfrac{16}9 -\tfrac83\,q - 71\,q^2 -\tfrac{2654}3\,q^3+\dots,\\
{{\mathfrak{u}}_4}^5 &= -1+15\,q + 45\,q^2 - 2145\,q^3 +\dots.
\end{align*}
There is a unique combination (up~to scalar multiples) that is zero to this
order, and therefore must vanish identically; namely,
\begin{equation*}
9\,{\mathfrak{u}}_4{\mathfrak{u}}_8^2 - 9\,{\mathfrak{u}}_6^2{\mathfrak{u}}_8 + 24\,{\mathfrak{u}}_4^3{\mathfrak{u}}_8 -15\,{\mathfrak{u}}_4^2{\mathfrak{u}}_6^2 + 16\,{\mathfrak{u}}_4^5.
\end{equation*}
Its vanishing is equivalent to $u_4u_8^2 - u_6^2u_8 + 24u_4^3u_8
-15u_4^2u_6^2 + 144u_4^5=0$.

\item $r=2,$ $\Gamma_0(N)=\Gamma_0(4)$.  No linear algebra is needed, since
as noted in the statement of Theorem~\ref{thm:main}(1),
$\mathscr{E}_2(q)=\mathscr{E}_4(q^2)$.  By comparing
\begin{align*}
{\mathfrak{u}}_4^{(4)}&=
4\mathscr{E}_4'-\mathscr{E}_4\cdot\mathscr{E}_4
&
{\mathfrak{u}}_{k+2}^{(4)}&={\mathfrak{u}}_k^{(4)}{}' - (k/4)\mathscr{E}_4\cdot
{\mathfrak{u}}_k^{(4)},\\
{\mathfrak{u}}_4^{(2)}&=
2\mathscr{E}_2'-\mathscr{E}_2\cdot\mathscr{E}_2
&
{\mathfrak{u}}_{k+2}^{(2)}&={\mathfrak{u}}_k^{(2)}{}' - (k/2)\mathscr{E}_r\cdot
{\mathfrak{u}}_k^{(2)},
\end{align*}
one deduces that
\begin{equation*}
  {\mathfrak{u}}_k^{(2)}(q) = 2^{(k-4)/2}\, {\mathfrak{u}}_k^{(4)}(q^2).
\end{equation*}
But in the $r=4$ case,
\begin{equation*}
  2\,{\mathfrak{u}}_4{\mathfrak{u}}_8 - 2\,{{\mathfrak{u}}_6}^2 + {{\mathfrak{u}}_4}^3 =0
\end{equation*}
(see the treatment above).  Hence, for $r=2,$
\begin{equation*}
  {\mathfrak{u}}_4{\mathfrak{u}}_8 - {{\mathfrak{u}}_6}^2 + 2\,{{\mathfrak{u}}_4}^3 =0.
\end{equation*}
This is equivalent to $u_4u_8 - u_6^2+8u_4^3=0$.
\end{itemize}

\subsection{A hypergeometric proof of Theorem~\ref{thm:main}(1)}
\label{subsec:52}

This proof is in the spirit of Jacobi, since it employs Hauptmoduls and
derivatives with respect to them.  It uses the results of~\cite{Maier12},
which were inspired by the following standard theorem on subgroup actions
of ${\it PSL}(2,\mathbf{R})$ \cite[\S\,44, Thm.~15]{Ford51}.

\begin{mytheorem}
\label{thm:automorphic}
  Let\/ $\Gamma<{\it PSL}(2,\mathbf{R})$ be a Fuchsian group of M\"obius
  transformations of\/~$\mathfrak{H}$ {\rm(}of~the first kind\/{\rm)} that
  has a Hauptmodul\/ $t=t(\tau),$ i.e., a non-constant simple automorphic
  function with a single simple zero on a fundamental region of\/~$\Gamma$.
  Then\/ $\tau$~can be expressed as a {\rm(}multivalued\/{\rm)} function
  of\/~$t$ as\/~$f_1/f_2,$ a~ratio of independent solutions\/ $f_1,f_2$ of
  some second-order differential equation
  \begin{equation}
    \label{eq:fuchsian}
    \mathcal{L}^{(\Gamma)}f:=\left[{D_t}^2 + P(t)\cdot D_t + Q(t)\right]f=0
  \end{equation}
  on\/~$\mathbf{P}^1(\mathbf{C})_t,$ in which\/ $P,Q\in\mathbf{C}(t)$.
\end{mytheorem}

Equation~(\ref{eq:fuchsian}) is called a Picard--Fuchs equation (the term
being historically most accurate when $\Gamma<\Gamma(1)$).  It is an ODE on
the genus-zero curve $\mathbf{P}^1(\mathbf{C})_t,$ which is essentially the
fundamental region of~$\Gamma$ with boundary identifications, i.e., the
(compactified) quotient of~$\mathfrak{H}$ by~$\Gamma$.  It follows from a
second theorem on automorphic functions~\cite[\S\,110, {Thm}.~6]{Ford51}
that Eq.~(\ref{eq:fuchsian}) must be a `Fuchsian' ODE, i.e., all its
singular points on $\mathbf{P}^1(\mathbf{C})_t$ must be regular.  These
points are bijective with the vertices of the fundamental region
of~$\Gamma$.  The difference of the two characteristic exponents of the
operator~$\mathcal{L}^{(\Gamma)}$ will be~$0$ at a cusp, and $1/n$ at an
order-$n$ elliptic fixed point.  That is, it will be the reciprocal of the
order of the associated stabilizing subgroup.

The Picard--Fuchs equation has solution space $\mathbf{C}f_1 \oplus
\mathbf{C}f_2,$ i.e., $(\mathbf{C}\tau\oplus\mathbf{C})f_2$.  It will
shortly be of interest to determine whether the logarithmic derivative
$u:=\dot f_2/f_2$ also satisfies an ODE, in this case with respect to
$\tau\in\mathfrak{H}\!$.  (As~always, the dot signifies differentiation
with respect to~$\tau$.)  For this, the following will be useful.  Let
$u_k,$ $k=4,6,\dotsc,$ be defined as in Theorem~\ref{thm:main}, i.e.,
$u_4=\dot u-u^2$ and $u_{k+2}=\dot u_k - kuu_k,$ and let differentiation
with respect to~$t$ be denoted by a subscripted~$t$.

\begin{mylemma}
\label{lem:useful}
  One can write\/ $u_k=\hat u_k \dot t^{k/2},$ where the sequence\/ $\hat
  u_4,\hat u_6,\dotsc$ follows from\/ $\hat u_4=-Q$ and\/ $\hat
  u_{k+2}=(\hat u_k)_t+(k/2)P\hat u_k$.  Thus
  \begin{align*}
    u_4 &= -Q\,\dot t^2,\\
    u_6 &= -(Q_t+2\,PQ)\dot t^3,\\
    u_8 &= -(Q_{tt}+5\,PQ_t+2\,QP_t + 6\,P^2Q)\dot t^4.
  \end{align*}
\end{mylemma}

\begin{proof}
  $\dot t=({\rm d}\tau/{\rm d}t)^{-1}= 1/(f_1/f_2)_t=f_2^2/w,$ where
  $w=w(f_1,f_2)$ is the Wronskian.  Similarly, $\ddot t=P\dot t^2+2u\dot t$
  comes from~(\ref{eq:fuchsian}) by differential calculus.  The recurrence
  $\hat u_{k+2}=(\hat u_k)_t+(k/2)P\hat u_k$ comes from $u_{k+2}=\dot u_k -
  kuu_k$ by substituting ${\rm d}/{\rm d}\tau=\dot t\,D_t,$ and exploiting
  these facts.
\end{proof}

Picard--Fuchs differential operators ${\mathcal L}^{(\Gamma)}$ that illustrate
Theorem~\ref{thm:automorphic} were obtained in~\cite{Maier12} for the
groups $\Gamma=\Gamma_0(N),$ $N=2,3,4,$ among others.  The corresponding
Hauptmoduls $t=t_N=t_N(\tau)$ were chosen to be
\begin{equation*}
  t_2(\tau)=2^{12}\cdot[2]^{24}/\,[1]^{24},\qquad
  t_3(\tau)=3^{6}\cdot[3]^{12}/\,[1]^{12},\qquad
  t_4(\tau)=2^{8}\cdot[4]^{8}/\,[1]^{8},
\end{equation*}
where the prefactors are of arithmetical significance but are not important
here (they could equally well be set equal to unity).  In each of these
three cases, $t_N=0$ corresponds to the infinite cusp, and $t_N=\infty$ to
the cusp $\tau=0$.  For $N=2,3,4,$ the operators
$\mathcal{L}_N={D_{{t_N}}}\!{\vphantom{D}}^2+P(t_N)\cdot D_{t_N}+Q(t_N)$
were computed to be
\begin{subequations}
\label{eq:ldefs}
  \begin{align}
    \mathcal{L}_2&={D_{{t_2}}}\!{\vphantom{D}}^2+\left[\frac{1}{t_2}+\frac{1}{2(t_2+64)}\right]D_{t_2}+\frac{1}{16\,t_2(t_2+64)},\\
    \mathcal{L}_3&={D_{{t_3}}}\!{\vphantom{D}}^2+\left[\frac{1}{t_3}+\frac{2}{3(t_3+27)}\right]D_{t_3}+\frac{1}{9\,t_3(t_3+27)},\\
    \mathcal{L}_4&={D_{{t_4}}}\!{\vphantom{D}}^2+\left[\frac{1}{t_4}+\frac{1}{t_4+16\mathstrut}\right]D_{t_4}+\frac{1}{4\,t_4(t_4+16)}.
  \end{align}
\end{subequations}
Each is a Gauss hypergeometric operator, up~to a scaling of the independent
variable.  That is, each has three (regular) singular points, located at
$t=t_N^*,\infty,0,$ where $t_N^*$ (respectively $-64,-27,-16$) is the third
fixed point of~$\Gamma_0(N)$ on the quotient curve
$X_0(N)=\Gamma_0(N)\setminus\mathfrak{H}^*\cong\mathbf{P}^1(\mathbf{C})_{t_N}$.
It is respectively a quadratic elliptic point, a cubic one, and a third
cusp (the image of $\tau=1/2$), as mentioned in~\S\,\ref{sec:intro}.

For each~$N,$ there is a solution $f=h_N(t_N)$ of $\mathcal{L}_Nf=0$ that
is holomorphic and equal to unity at~$t_N=0$.  It was shown
in~\cite{Maier12} that if ${\mathfrak{h}}_N(\tau):=(h_N\circ t_N)(\tau),$
then
\begin{equation*}
  {\mathfrak{h}}_2(\tau)=[1]^4/\,[2]^2,\qquad
  {\mathfrak{h}}_3(\tau)=[1]^3/\,[3],\qquad 
  {\mathfrak{h}}_4(\tau)=[1]^4/\,[2]^2.
\end{equation*}
That is, the holomorphic local solution of~(\ref{eq:fuchsian}) at the
infinite cusp, in each case, can be continued to a weight-$1$ modular form
on~$\mathfrak{H}^*$.  In~fact,
${\mathfrak{h}}_2,{\mathfrak{h}}_3,{\mathfrak{h}}_4$ are respectively equal
to ${\mathscr{B}}_4,{\mathscr{B}}_3,{\mathscr{B}}_2$ in the notation of the
present article (see Definition~\ref{def:abc}).

For $N=2,3,4,$ a weight-$1$ modular form $\bar {\mathfrak{h}}_N(\tau) =
(\bar{h}_N\circ t_N) (\tau)$ that vanishes at the infinite cusp, and has
zero order of vanishing at the cusp~$\tau=0,$ is obtained by multiplying
$\mathfrak{h}_N(\tau)$ by an appropriate power of the
Hauptmodul~$t_N(\tau)$.  Let
\begin{equation*}
  \bar{h}_2(t_2) = 2^{-3/2}\, t_2^{1/4}{h}_2(t_2),\quad
  \bar{h}_3(t_3) = 3^{-1}\, t_3^{1/3}{h}_3(t_3),\quad
  \bar{h}_4(t_4) = 2^{-2}\, t_4^{1/2}{h}_4(t_4).
\end{equation*}
Then by Definition~\ref{def:abc},
$\bar{\mathfrak{h}}_2(\tau),\bar{\mathfrak{h}}_3(\tau),\bar{\mathfrak{h}}_4(\tau)$
are identical to~$\mathscr{C}_4,\mathscr{C}_3,\mathscr{C}_2$.  It follows
by changing (dependent) variables in the equations $\mathcal{L}_Nh_N=0$
that $\bar h_N$~satisfies the slightly modified Picard--Fuchs equation
$\bar{\mathcal{L}}_N{\bar h}_N=0,$ where
\begin{subequations}
\label{eq:lbardefs}
  \begin{align}
\bar{\mathcal{L}}_2&= {D_{{t_2}}}\!{\vphantom{D}}^2 + \left[ \frac1{2t_2} + \frac1{2({t_2}+64)}\right]D_{t_2} + \frac4{t_2^2(t_2+64)},\\
\bar{\mathcal{L}}_3&= {D_{{t_3}}}\!{\vphantom{D}}^2 + \left[ \frac1{3t_3} + \frac2{3(t_3+27)}\right]D_{t_3} + \frac3{t_3^2(t_3+27)},\\
\bar{\mathcal{L}}_4&= {D_{{t_4}}}\!{\vphantom{D}}^2 + \left[ \frac0{t_4} + \frac1{t_4+16\mathstrut}\right]D_{t_4} + \frac4{t_4^2(t_4+16)}.
  \end{align}
\end{subequations}
The fixed points of each $\Gamma_0(N)$ on the corresponding quotient
$X_0(N)\cong\mathbf{P}^1(\mathbf{C})_{t_N}$ are visible
in~(\ref{eq:lbardefs}), just as in~(\ref{eq:ldefs}).

Each modified equation $\bar{\mathcal{L}}_N f=0$ is of the form
$\mathcal{L}_{\alpha,\beta,\gamma}f=0,$ where
\begin{equation}
\label{eq:gen2F1}
\mathcal{L}_{\alpha,\beta,\gamma}:= {D_t}^2 + \left[ \frac{\alpha+\beta}t +
  \frac{1-\alpha}{t-t^*}\right]D_t + \frac{[\gamma^2 - (1-\alpha-\beta)^2]t^*}{4\,t^2(t-t^*)}.
\end{equation}
The operator $\mathcal{L}_{\alpha,\beta,\gamma}$ is the general
second-order Fuchsian operator on~$\mathbf{P}^1(\mathbf{C})_t$ that has
singular points at~$t=t^*,\infty,0$ with respective exponent differences
$\alpha,\beta,\gamma,$ and with one exponent at each of $t=t^*,\infty$
constrained to be zero.  It is of hypergeometric but not
Gauss-hypergeometric type.  The solutions of
$\mathcal{L}_{\alpha,\beta,\gamma}f=0$ include
\begin{equation}
  t^{(1-\alpha-\beta-\gamma)/2}\:{}_2F_1\left(\tfrac{1-\alpha-\beta-\gamma}2,\tfrac{1-\alpha+\beta-\gamma}2;\,1-\gamma\,;\:t/t^*\right),
\end{equation}
which is the local solution at~$t=0$ associated to the exponent
${(1-\alpha-\beta-\gamma)/2}$.  (Here, ${}_2F_1(\lambda,\mu;\nu;x)$ is the
Gauss hypergeometric function, defined and single-valued on the disk
$\left|x\right|<1$.)  This is the representation of the form $\bar
{\mathfrak{h}}_N=\mathscr{C}_r$ as a (multivalued) function of the
Hauptmodul~$t=t_N$.  For $N=2,3,4,$ the parameters $(\alpha,\beta,\gamma),$
which are the reciprocals of the orders of elements of~$\Gamma_0(N)$ that
stabilize the corresponding fixed points, are respectively $(\frac12,0,0),$
$(\frac13,0,0),$ $(0,0,0)$.

By Theorem~\ref{thm:main}(2), $u=(2\pi{\rm i}/r)\mathscr{E}_r$ equals
$\dot{\mathscr{C}}_r/\mathscr{C}_r$.  Hence, if one takes the operator
$\mathcal{L}^{(\Gamma)}=\mathcal{L}^{(\Gamma_0(N))}$ to equal
$\bar{\mathcal{L}}_N$ rather than~${\mathcal{L}}_N,$ the function~$u$ of
Theorem~\ref{thm:main}(1) will agree with the function $u=\dot f_2/f_2,$ in
the notation of Theorem~\ref{thm:automorphic} and Lemma~\ref{lem:useful}.
Therefore Theorem~\ref{thm:main}(1) will follow immediately from

\begin{myproposition}
  \label{prop:special}
Let\/ $\tau=f_1/f_2,$ a ratio of independent solutions of the
hypergeometric ODE\/ $\mathcal{L}_{\alpha,\beta,\gamma}f=0$
on\/~$\mathbf{P}^1(\mathbf{C})_t$.  Let\/ $u=\dot f_2/f_2,$ and let\/
$u_k,$ $k=4,6,\dotsc$ be defined by\/ $u_4=\dot u-u^2$ and\/ $u_{k+2}=\dot
u_k-kuu_k$.  Then\/ $u,$ regarded as a function of\/ $\tau,$ will satisfy a
nonlinear third-order ODE: the generalized Chazy equation
\begin{alignat*}{3}
  &u_4u_8 - u_6^2 + 8\,u_4^3=0,&\qquad&\text{if}&\quad&(\alpha,\beta,\gamma)=(\tfrac12,0,0),\\
  &u_4u_8^2 - u_6^2u_8 + 24\,u_4^3u_8 -15\,u_4^2u_6^2 + 144\,u_4^5=0,&\qquad&\text{if}&\quad&(\alpha,\beta,\gamma)=(\tfrac13,0,0),\\
  &u_4u_8 - u_6^2 + 8\,u_4^3=0,&\qquad&\text{if}&\quad&(\alpha,\beta,\gamma)=(0,0,0).
\end{alignat*}
\end{myproposition}
\begin{proof}
  By direct computation, using the expressions of Lemma~\ref{lem:useful}
  for $u_4,u_6,u_8$ in~terms of the coefficient functions
  $P,Q\in\mathbf{C}(t),$ which can be read~off from the
  formula~(\ref{eq:gen2F1}) for~$\mathcal{L}_{\alpha,\beta,\gamma}$.
\end{proof}

So, in each of the cases $N=2,3,4,$ i.e., $r=4,3,2,$ the quasi-modular form
$u=\dot{\mathscr{C}}_r/\mathscr{C}_r$ satisfies a generalized Chazy
equation.  This proof of Theorem~\ref{thm:main}(1) is more analytic than
the proof given in~\S\,\ref{subsec:51}, and less explicitly modular.

The reader may wonder whether this second, alternative proof was necessary.
It required extra machinery, such as the Picard--Fuchs equations
$\mathcal{L}_Nh_N=0$ and ${\bar{\mathcal{L}}_N{\bar h}_N=0},$ and the
Hauptmoduls~$t_N,$ the $q$-expansions of which are relatively complicated
and are not discussed here.  Also, Proposition~\ref{prop:special} is
restricted to very special triples of parameter values.

In fact, Proposition~\ref{prop:special} is the tip of an interesting
iceberg.  The following is its extension to arbitrary triples
$(\alpha,\beta,\gamma)$.

\begin{mytheorem}
  \label{thm:general}
  Let\/ $\tau=f_1/f_2,$ a ratio of independent solutions of the
  hypergeometric ODE\/ $\mathcal{L}_{\alpha,\beta,\gamma}f=0$
  on\/~$\mathbf{P}^1(\mathbf{C})_t$.  Let\/ $u=\dot f_2/f_2,$ and let\/
  $u_k,$ $k=4,6,\dotsc$ be defined by\/ $u_4=\dot u-u^2$ and\/
  $u_{k+2}=\dot u_k-kuu_k$.  Then\/ $u,$ regarded as a function of\/
  $\tau,$ will satisfy a nonlinear third-order ODE: a generalized Chazy
  equation
  \begin{equation*}
    C_{88}\,u_4^2u_8^2 + C_{86}\,u_4u_6^2u_8 + C_{84}\,u_4^4u_8 + C_{66}\,u_6^4 +
    C_{64}\,u_4^3u_6^2 + C_{44}\,u_4^6 = 0
  \end{equation*}
  with coefficients symmetric under\/ $\alpha\leftrightarrow\beta,$ namely
  \begin{align*}
    C_{88} &= (2\alpha-1)(2\beta-1)(\alpha+\beta-\gamma-1)^2(\alpha+\beta+\gamma-1)^2    ,\\
    C_{86} &= -
\bigl[(2\alpha-1)(3\beta-1)+(3\alpha-1)(2\beta-1)\bigr]
(\alpha+\beta-\gamma-1)^2(\alpha+\beta+\gamma-1)^2  ,\\
    C_{84} &= -16(2\alpha-1)(2\beta-1)(\alpha+\beta-1)(\alpha+\beta-\gamma-1)(\alpha+\beta+\gamma-1)  ,\\
    C_{66} &= (3\alpha-1)(3\beta-1)(\alpha+\beta-\gamma-1)^2(\alpha+\beta+\gamma-1)^2 ,\\
    C_{64} &= 4
\left[ 2(2\alpha-1)^2(3\beta-1) + 2(3\alpha-1)(2\beta-1)^2 - 3(\alpha-\beta)^2 \right]
\\
           &\qquad\qquad\qquad\qquad\qquad\qquad\qquad{}\times(\alpha+\beta-\gamma-1)(\alpha+\beta+\gamma-1) ,\\
    C_{44} &= 64(2\alpha-1)(2\beta-1)(\alpha+\beta-1)^2.
  \end{align*}
\end{mytheorem}
\begin{proof}
  With the aid of a computer algebra system, eliminate~$t$ from the
  expressions for $u_8/u_4^2$ and $u_6^2/u_4^3$ that follow from
  Lemma~\ref{lem:useful}.  As in the proof of
  Proposition~\ref{prop:special}, $P,Q\in\mathbf{C}(t)$ come from
  Eq.~(\ref{eq:gen2F1}).
\end{proof}

Theorem~\ref{thm:general} is a \emph{nonlinear hypergeometric identity},
relevant even to hypergeometric equations without modular applications.  It
belongs to the theory of special functions, but as one sees, in a loose
sense it is a relation of linear dependence among modular forms of
weight~$24$ (since each monomial has that weight).

Rational exponent differences $\alpha,\beta,\gamma\,$ that are not members
of $\{0,\tfrac12,\tfrac13\}$ occur in the theory of automorphic functions
on subgroups of~${\it PSL}(2,\mathbf{R})$ that are not subgroups of
$\Gamma(1)={\it PSL}(2,\mathbf{Z})$.  This will be illustrated in the next
section.

\section{Differential systems for weight-$1$ forms}
\label{sec:final}

The systems of Theorem~\ref{thm:main}(2), satisfied by triples
$\mathscr{A}_r,\mathscr{B}_r,\mathscr{C}_r,$ will now be greatly
generalized.  Associated to any first-kind Fuchsian subgroup $\Gamma<{\it
PSL}(2,\mathbf{R})$ that is a triangle group, i.e., that has a hyperbolic
triangular fundamental domain in~$\mathfrak{H}^*$ and a Hauptmodul, there
are weight-$1$ modular forms $\mathcal{A}\!,\mathcal{B},\mathcal{C}$
(possibly multivalued) that vanish respectively at the three vertices.  The
forms will satisfy a system of coupled nonlinear first-order ODEs with
independent variable~$\tau$.

The key result on this is Theorem~\ref{thm:triangular}, which is proved by
hypergeometric manipulations related to those of Ohyama~\cite{Ohyama96}.
It deals with the case when $\Gamma$~has at~least one cusp, which without
loss of generality can be taken to be $\tau={\rm i}\infty$.
In~\S\S\,\ref{subsec:final2} and~\ref{subsec:final3}, as illustrations, the
triangle subgroups $\Gamma<{\it PSL}(2,\mathbf{R})$ that are commensurable
with $\Gamma(1)={\it PSL}(2,\mathbf{Z})$ are examined.  These include
$\Gamma(1)$ itself; $\Gamma_0(N),$ $N=2,3$; the Fricke extensions
$\Gamma_0^+(N),$ $N=2,3$; the index-$2$ subgroup $\Gamma^2$ of~$\Gamma(1)$;
and two others~\cite{Takeuchi77}.  For each such~$\Gamma\!,\,$ the forms
$\mathcal{A}\!,\mathcal{B},\mathcal{C}$ are worked~out explicitly, as are
the hypergeometric representation of~$\mathcal{A},$ the differential system
the forms satisfy, and the generalized Chazy equation that the system
implies.

\subsection{Hypergeometric manipulations}
\label{subsec:final1}

Let $\Gamma<{\it PSL}(2,\mathbf{R})$ be a Fuchsian subgroup (of the first
kind), regarded as a group of M\"obius transformations of
$\mathfrak{H}^*\ni\tau$.  If it has a Hauptmodul $t=t(\tau),$ i.e., is of
genus zero, then
$\Gamma\setminus\mathfrak{H}^*\cong\mathbf{P}^1(\mathbf{C})_t$ and
$\tau$~can be expressed as a ratio $\tau_1/\tau_2$ of two solutions
$f=\tau_1,\tau_2$ of a Picard--Fuchs equation $\mathcal{L}^{(\Gamma)}f=0$
on~$\mathbf{P}^1(\mathbf{C})_t,$ as stated in
Theorem~\ref{thm:automorphic}.  Its solution space will be
$\mathbf{C}\tau_1 \oplus \mathbf{C}\tau_2 =
(\mathbf{C}\tau\oplus\mathbf{C})\tau_2$.  The solution~$\tau_2$ can be
viewed as a weight-$1$ form on~$\Gamma$.  This follows by considering the
homogeneous counterpart $\overline\Gamma<{\it SL}(2,\mathbf{R})$
to~$\Gamma\!,\,$ which acts on vectors~$\binom{\tau_1}{\tau_2},$ and the
associated homogeneous forms, which are functions of~$\tau_1,\tau_2$.

Suppose that $\Gamma$ is a triangle group, i.e., is of genus zero with a
triangular fundamental domain and hence with three inequivalent classes of
fixed points on~$\mathfrak{H}^*\!,\,$ say classes A,B,C\null.  They
correspond to three conjugacy classes of stabilizing elements
of~$\Gamma\!,\,$ either elliptic or parabolic.  The group $\Gamma$ is
specified up~to conjugacy by their orders, i.e., by a signature
$(n_{\mathcal{A}},n_{\mathcal{B}},n_{\mathcal{C}})$ such as the signature
$(3,2,\infty)$ of~$\Gamma(1)$.  It will be assumed that at~least one of
these orders is~$\infty,$ i.e., at least one of the three classes is
parabolic, corresponding to a cusp.  Without loss of generality one can
take $n_{\mathcal{C}}=\infty,$ and the cusp to be $\tau={\rm i}\infty$.
This infinite cusp will be fixed by some $\tau\mapsto\tau+\upsilon,$ where
by definition, $\upsilon\in\mathbf{R}^+$ is its width.

By a M\"obius transformation the Hauptmodul~$t$ can be redefined, if
necessary, so that $t=0$ at the infinite cusp, and $t=t^*,\infty,$ for some
$t^*\in\mathbf{C}\setminus\{0\},$ on the fixed points in the
classes~A,B\null.  The Picard--Fuchs equation, being hypergeometric, will
then have $t=t^*,\infty,0$ as its (regular) singular points.  Their
respective exponent differences $(\alpha,\beta,\gamma\,)$ will equal
$(1/n_{\mathcal{A}},1/n_{\mathcal{B}},1/n_{\mathcal{C}}\!)$.  These are
vertex angles in~terms of $\pi$~radians, and necessarily
$\alpha+\beta+\gamma<1$.  (If the convention of the last paragraph is
adopted then $n_{\mathcal{C}}=\infty$ and $\gamma=1/\infty=0,$ but for
reasons of symmetry $n_{\mathcal{C}}$ and~$\gamma$ will be kept as
independent parameters.)

Fuchs's relation on characteristic exponents implies that the six exponents
of any second-order ODE of hypergeometric type must sum to unity.  This
leaves two degrees of freedom in the choice of~$\mathcal{L}^{(\Gamma)},$ as
in~\S\,\ref{subsec:52}.  (Cf.\ ${\mathcal{L}}_N$
vs.~$\bar{\mathcal{L}}_N$.)  Let ${\mathcal{L}}^{(\Gamma)}$ be chosen to
have exponents
\begin{equation}
  \{0,\alpha\}\text{ at }t=t^*,\quad
  \{0,\beta\}\text{ at }t=\infty,\quad
  \{\tfrac{1-\alpha-\beta-\gamma}2,\tfrac{1-\alpha-\beta+\gamma}2\}\text{ at }t=0,
\end{equation}
i.e., so that there is a zero exponent at each of $t=t^*,\infty$.  With this
choice,
\begin{equation}
\label{eq:Labc}
\mathcal{L}^{(\Gamma)}=\mathcal{L}_{\alpha,\beta,\gamma}:= {D_t}^2 + \left[ \frac{\alpha+\beta}t +
  \frac{1-\alpha}{t-t^*}\right]D_t + \frac{[\gamma^2 - (1-\alpha-\beta)^2]t^*}{4\,t^2(t-t^*)},
\end{equation}
defined as in (\ref{eq:gen2F1}).  Let the local solution of
$\mathcal{L}^{(\Gamma)}f=\mathcal{L}_{\alpha,\beta,\gamma}f=0$ at the
singular point $t=0$ (i.e., the infinite cusp), associated to the exponent
$(1-\alpha-\beta-\gamma)/2,$ be denoted $C=C(t)$.  Then the lifted function
$\mathcal{C}(\tau):=C(t(\tau))$ will be a weight-$1$ form on~$\Gamma\!,\,$
which vanishes at cusps in class~C because $\alpha+\beta+\gamma<1$.

Also, define (potentially multivalued) functions
$\mathcal{A}(\tau),\mathcal{B}(\tau)$ that vanish at the fixed points in
the classes~A,B, at which $t=t^*,$ resp.\ $t=\infty,$ by
\begin{subequations}
\label{eq:varpropto}
\begin{align}
\mathcal{A} & = \left[(t-t^*)/t\right]^{1/\rho} \mathcal{C} ,\\
\mathcal{B} & = \left[-t^*/t\right]^{1/\rho} \mathcal{C} ,
\end{align}
\end{subequations}
where
\begin{equation}
\label{eq:rhodef}
  \rho := \frac2{1-\alpha-\beta-\gamma} =
  \frac2{1-n_{\mathcal{A}}^{-1}-n_{\mathcal{B}}^{-1}-n_{\mathcal{C}}^{-1}}.
\end{equation}
With these definitions,
\begin{equation}
  {\mathcal{A}}^\rho = {\mathcal{B}}^\rho + {\mathcal{C}}^\rho.
\end{equation}
The corresponding quotients
\begin{subequations}
  \label{eq:gathered}
  \begin{gather}
    t/(t-t^*)=\mathcal{C}^\rho/\mathcal{A}^\rho,\\
    t/t^*=-\mathcal{C}^\rho/\mathcal{B}^\rho,
  \end{gather}
\end{subequations}
are normalized Hauptmoduls, the respective values of which on the classes
A,B,C are $\infty,1,0$ and~$1,\infty,0$.  The $\Gamma$-specific quantity
$\rho\in\mathbf{Q}^+$ of~(\ref{eq:rhodef}) generalizes the `signature'~$r$
that parametrizes Ramanujan's alternative theories of elliptic integrals.
(Recall that $r=4,3,2$ correspond to
$\Gamma=\Gamma_0(2),\allowbreak\Gamma_0(3),\Gamma_0(4),$ i.e., to
$(n_{\mathcal{A}},n_{\mathcal{B}},n_{\mathcal{C}})=\allowbreak(2,\infty,\infty),\allowbreak(3,\infty,\infty),\allowbreak(\infty,\infty,\infty)$.)

The functions $\mathcal{A}(\tau),\mathcal{B}(\tau)$ could also be defined
as $A(t(\tau))$ and $B(t(\tau)),$ where $A(t),B(t)$ are solutions of
Picard--Fuchs equations having appropriately modified exponents, but the
same exponent \emph{differences}
as~${\mathcal{L}}^{(\Gamma)}_{\alpha,\beta,\gamma},$ i.e.,
$\alpha,\beta,\gamma$.  (Cf.\ the relation between
$\mathcal{L}_N,\bar{\mathcal{L}}_N$ in~\S\,\ref{subsec:52}.)  Their
respective exponents would be
\begin{subequations}
\begin{align}
  &\{\tfrac{1-\alpha-\beta-\gamma}2,\tfrac{1+\alpha-\beta-\gamma}2\}\text{ at }t=t^*,\quad
  \{0,\beta\}\text{ at }t=\infty,\quad
  \{0,\gamma\}\text{ at }t=0,\\
  &\{0,\alpha\}\text{ at }t=t^*,\quad
  \{\tfrac{1-\alpha-\beta-\gamma}2,\tfrac{1-\alpha+\beta-\gamma}2\}\text{ at }t=\infty,\quad
  \{0,\gamma\}\text{ at }t=0.
\end{align}
\end{subequations}
Gauss-hypergeometric representations of
$\mathcal{A}\!,\mathcal{B},\mathcal{C}$ in~terms of~${}_2F_1$ are
\begin{subequations}
\label{eq:rewritable}
  \begin{align}
\mathcal{A}(\tau) & =  A(t(\tau)) = {}_2F_1\bigl(\tfrac{1-\alpha-\beta-\gamma}2, \tfrac{1+\alpha-\beta-\gamma}2;\,1-\gamma\,;\;t(\tau)/\left[t(\tau)-t^*\right]\bigr),\\
\mathcal{B}(\tau) & =  B(t(\tau)) = {}_2F_1\bigl(\tfrac{1-\alpha-\beta-\gamma}2, \tfrac{1-\alpha+\beta-\gamma}2;\,1-\gamma\,;\;t(\tau)/t^*),\\
\mathcal{C}(\tau) & =  C(t(\tau)) = (-t/t^*)^{1/\rho}{}_2F_1\bigl(\tfrac{1-\alpha-\beta-\gamma}2, \tfrac{1-\alpha+\beta-\gamma}2;\,1-\gamma\,;\;t(\tau)/t^*),\label{eq:rewritablec}
  \end{align}
\end{subequations}
in which the normalizations of $\mathcal{A}\!,\mathcal{B},\mathcal{C},$ not
previously specified, have been set by requiring that
$\mathcal{A}\!,\mathcal{B}$ equal unity at the infinite cusp, at~which
$t=0$.  The parameters and arguments of the ${}_2F_1$'s are determined by
${}_2F_1(\lambda,\mu;\nu; x)$ having exponents $\{0,1-\nu\}$ at~$x=0,$
$\{0,\nu-\lambda-\mu\}$ at~$x=1,$ and $\{\lambda,\mu\}$ at~$x=\infty$.  

The representations
\begin{subequations}
\label{eq:rewritable2}
  \begin{align}
\mathcal{A} & = {}_2F_1\bigl(\tfrac{1-\alpha-\beta-\gamma}2,
\tfrac{1+\alpha-\beta-\gamma}2;\,1-\gamma\,;\;\mathcal{C}^\rho/\mathcal{A}^\rho\bigr),\\
\mathcal{B} & = {}_2F_1\bigl(\tfrac{1-\alpha-\beta-\gamma}2,
\tfrac{1-\alpha+\beta-\gamma}2;\,1-\gamma\,;\;-\mathcal{C}^\rho/\mathcal{B}^\rho\bigr)
  \end{align}
\end{subequations}
follow from (\ref{eq:rewritable}ab) with the aid of~(\ref{eq:gathered}ab).
The identities (\ref{eq:rewritable2}ab) are equivalent: they are related by
Pfaff's transformation of~${}_2F_1$.  The function
${}_2F_1(\lambda,\mu;\nu;x)$ is defined on the disk $\left| x\right|<1,$ so
(\ref{eq:rewritable2}ab) hold in a neighborhood of the infinite cusp, at
which $\mathcal{C}=0$.  In~fact each will hold near any cusp in the
class~C, if an appropriate constant of proportionality is included.
Similarly, there follow
\begin{subequations}
\label{eq:rewritable3}
  \begin{align}
\mathcal{C} & \varpropto {}_2F_1\bigl(\tfrac{1-\alpha-\beta-\gamma}2,
\tfrac{1-\alpha-\beta+\gamma}2;\,1-\alpha\,;\;\mathcal{A}^\rho/\mathcal{C}^\rho\bigr),\\
\mathcal{C} & \varpropto {}_2F_1\bigl(\tfrac{1-\alpha-\beta-\gamma}2,
\tfrac{1-\alpha-\beta+\gamma}2;\,1-\beta\,;\;-\mathcal{B}^\rho/\mathcal{C}^\rho\bigr),
  \end{align}
\end{subequations}
which hold near any fixed point in the class~A, resp.~B, with the constant
of proportionality dependent on the choice of fixed point.

As defined, $\mathcal{A}\!,\mathcal{B},\mathcal{C}$ are formally weight-$1$
modular forms on~$\Gamma\!,\,$ with some multiplier systems; and they
vanish respectively on the classes A,B,C of fixed points of~$\Gamma$
on~$\mathfrak{H}^*$.  However, if A, resp.~B comprises \emph{elliptic}
fixed points, then $\mathcal{A}\!,$ resp.~$\mathcal{B}$ may be a
multivalued function of~$\tau$.  This is because of the fractional powers
in their definitions~(\ref{eq:varpropto}ab).

The test for multivaluedness on~$\mathfrak{H}$ is as follows.  By
construction, each of $\mathcal{A}\!,\mathcal{B},\mathcal{C}$ has order of
vanishing (computed with respect to a local parameter for~$\Gamma\!,\,$
e.g.,~$t$) equal to~$1/\rho$.  Fixed points in classes A,B,C are mapped to
$\Gamma\setminus\mathfrak{H}^*\cong\mathbf{P}^1(\mathbf{C})_t$ with
multiplicities $n_\mathcal{A},n_\mathcal{B},n_\mathcal{C}$.  If
$n_\mathcal{A}<\infty,$ resp.\ $n_\mathcal{B}<\infty,$ signaling
ellipticity, then the order of vanishing of $\mathcal{A}\!,\mathcal{B}$ at
the associated elliptic points on~$\mathfrak{H},$ in classes~A,B, will be
$n_\mathcal{A}/\rho,n_\mathcal{B}/\rho\in\mathbf{Q}^+$.  If this is not an
integer then $\mathcal{A},$ resp.~$\mathcal{B}$ will be multivalued, i.e.,
the $k$'th root of a true modular form (of weight~$k$), where $k$ equals
the numerator of the fraction $\rho/n_\mathcal{A}=\rho\alpha,$ resp.\
$\rho/n_\mathcal{B}=\rho\beta,$ expressed in lowest terms.

The generalization of Theorem~\ref{thm:main}(2) to arbitrary triangle
groups can now be stated and proved.  As always, ${}'$ signifies $q\,{\rm
d}/{\rm d}q = (2\pi{\rm i})^{-1}{\rm d}/{\rm d}\tau$.

\begin{mytheorem}
\label{thm:triangular}
Let\/\, $\Gamma<{\it PSL}(2,\mathbf{R})$ be a triangle group with
signature\/ $(n_\mathcal{A}, n_\mathcal{B}, n_\mathcal{C}\!)$ and exponent
parameters\/ $(\alpha,\beta,\gamma)=(1/n_\mathcal{A}, 1/n_\mathcal{B},
1/n_\mathcal{C}),$ with\/ $\alpha+\beta+\gamma<1,$ and let $\rho :=
\frac2{1-\alpha-\beta-\gamma}$.  Assume that\/ $\gamma=0,$ i.e., that the
third vertex is a cusp, and define the formal weight\/-$1$ modular forms\/
$\mathcal{A}\!,\mathcal{B},\mathcal{C}$ as above, vanishing at the
classes\/ {\rm A,B,C} of fixed points of\/ $\Gamma$ on\/ $\mathfrak{H}^*$
(the last containing the infinite cusp, of width\/ $\upsilon$), and
satisfying\/ $\mathcal{A}^\rho = \mathcal{B}^\rho + \mathcal{C}^{\rho}$.
Then ${\mathcal{A}}^\rho,{\mathcal{B}}^\rho,{\mathcal{C}}^\rho,$ along
with\/ the weight\/-$2,$ depth\/-${\le1}$ quasi-modular form\/
$\mathcal{E}:=\upsilon ({\mathcal{C}}^\rho)'/({\mathcal{C}}^\rho)$ that is
associated with class\/ {\rm C,} satisfy the coupled system of nonlinear
first-order equations
\begin{subequations}
\label{eq:last}
\begin{align}
	  \upsilon({\mathcal{A}}^\rho)' &= \mathcal{E}\cdot {\mathcal{A}}^\rho - {\mathcal{A}}^{\rho(1-\alpha)}{\mathcal{B}^{\rho(1-\beta)}},\label{eq:lasta}\\
	  \upsilon({\mathcal{B}}^\rho)' &= \mathcal{E}\cdot {\mathcal{B}}^\rho - {\mathcal{A}}^{\rho(1-\alpha)}{\mathcal{B}^{\rho(1-\beta)}},\label{eq:lastb}\\
	  \upsilon({\mathcal{C}}^\rho)' &= \mathcal{E}\cdot {\mathcal{C}}^{\rho\vphantom{()}},\label{eq:lastc}\\
	  \upsilon\rho\,\mathcal{E}'&= \mathcal{E}\cdot \mathcal{E} - {\mathcal{A}}^{\rho(1-2\alpha)}{\mathcal{B}^{\rho(1-2\beta)}},\label{eq:lastd}
\end{align}
\end{subequations}
from which a generalized Chazy equation\/ ${C}_{p}$ for\/
$u=(2\pi{\rm i}/\upsilon\rho)\mathcal{E},$ parametrized by\/
$\alpha,\beta;\gamma\,$ as in Theorem\/ {\rm\ref{thm:general}}, can be
derived by elimination.  (The third equation says that\/
$u=\dot{\mathcal{C}}/\mathcal{C}$.)
\end{mytheorem}
\begin{proof}
Equation~(\ref{eq:lastc}) is true by definition, and (\ref{eq:lasta}) is
implied by (\ref{eq:lastb}),\allowbreak(\ref{eq:lastc}) and
$\mathcal{A}^\rho = \mathcal{B}^\rho + \mathcal{C}^{\rho}$.  It remains to
prove (\ref{eq:lastb}),(\ref{eq:lastd}).  They will come from a useful
formula for $\dot t = dt/d\tau,$ deduced as follows.

The solution space of the Picard--Fuchs equation
$\mathcal{L}_{\alpha,\beta,\gamma}f=0$ is $\mathbf{C}f_1\oplus\mathbf{C}f_2
:=\allowbreak \mathbf{C}\tau\, C\oplus\mathbf{C}\,C =
\left(\mathbf{C}\tau\oplus\mathbf{C}\right)C,$ where $\tau$~is viewed as a
(multivalued) function on the quotient curve $\mathbf{P}^1(\mathbf{C})_t,$
and $C$~is defined in~terms of ${}_2F_1$ by~(\ref{eq:rewritablec}).  From
the expression for~$\mathcal{L}_{\alpha,\beta,\gamma}$ given
in~(\ref{eq:Labc}), the Wronskian $w=w(f_1,f_2)=w(\tau C,C)$ must equal a
multiple of $1/t^{\alpha+\beta}(t-t^*)^{1-\alpha}$.  The constant of
proportionality can be calculated by taking the $t\to0$ limit, in which the
infinite cusp is approached.  In this limit,
\begin{equation}
  \tau\sim(\upsilon/2\pi{\rm i})\log t,\qquad C\sim(-t/t^*)^{1/\rho},
\end{equation}
the former coming from $t\sim\text{const}\cdot q^{1/\upsilon}\!,$ which is
true since $\upsilon$ is~the width of the infinite cusp.  One readily
deduces that
\begin{equation}
  1/w(t) = \left[2\pi{\rm i}(-t^*)^{-\beta}/\upsilon\right]\cdot
  t^{\alpha+\beta} (t-t^*)^{1-\alpha}.
\end{equation}
As in the proof of Lemma~\ref{lem:useful}, $\dot t=f_2^2/w,$ i.e., $\dot
t=\mathcal{C}^2(\tau) / w(t(\tau))$ on~$\mathfrak{H}$.  Taking account
of~(\ref{eq:varpropto}ab), one can rewrite this in the useful form
\begin{equation}
\label{eq:dott}
\dot t = [2\pi{\rm i}\,(-t^*)/\upsilon]\,\mathcal{A}^{\rho(1-\alpha)}\mathcal{B}^{-\rho(1+\beta)}\mathcal{C}^\rho.  
\end{equation}
Now consider the logarithmic derivative of the equality $\mathcal{B}^\rho =
(-t^*/t)\mathcal{C}^\rho,$ i.e.,
\begin{equation}
  ({\mathcal{B}^\rho})'/\mathcal{B}^\rho
=   ({\mathcal{C}^\rho})'/\mathcal{C}^\rho
  - t'/t.
\end{equation}
By employing (\ref{eq:dott}) to expand $t'=\dot t/2\pi{\rm i},$ one obtains
Eq.~(\ref{eq:lastb}).  Equation (\ref{eq:lastd}) follows by similar
manipulations.
\end{proof}

Let $\mathcal{E}_{A}$ and $\mathcal{E}_{B}$ denote the weight-$2,$
depth-$\le1$ quasi-modular forms associated with classes A and~B, i.e.,
$\upsilon({\mathcal{A}}^\rho)'/{\mathcal{A}}^\rho$ and
$\upsilon({\mathcal{B}}^\rho)'/{\mathcal{B}}^\rho,$ just as the form
$\mathcal{E}_\mathrm{C}=\mathcal{E}=\upsilon({\mathcal{C}}^\rho)'/{\mathcal{C}}^\rho$
is associated with class~C\null.  Moreover, let
$u_\mathrm{A},u_\mathrm{B},u_\mathrm{C}$ denote the normalized forms
$(2\pi{\rm i}/\upsilon\rho)\mathcal{E}_\mathrm{A},$ $(2\pi{\rm
i}/\upsilon\rho)\mathcal{E}_\mathrm{B},$ $(2\pi{\rm
i}/\upsilon\rho)\mathcal{E}_\mathrm{C},$ so that
$u_\mathrm{A}=\dot{\mathcal{A}}/\mathcal{A},$
$u_\mathrm{B}=\dot{\mathcal{B}}/\mathcal{B},$
$u_\mathrm{C}=\dot{\mathcal{C}}/\mathcal{C}$.  Then a bit of calculus
applied to Eqs.\ (\ref{eq:last}abcd) yields

\begin{mycorollary}
  The weight\/-$2,$ depth\/-$\le1$ quasi-modular forms\/
  $u_\mathrm{A},u_\mathrm{B},u_\mathrm{C}$ satisfy
  \begin{align*}
    \dot u_\mathrm{A} &= u_\mathrm{A}^2 - (1+\rho\alpha)(u_\mathrm{A}-u_\mathrm{B})(u_\mathrm{A}-u_\mathrm{C}),\\
    \dot u_\mathrm{B} &= u_\mathrm{B}^2 - (1+\rho\beta)(u_\mathrm{B}-u_\mathrm{C})(u_\mathrm{B}-u_\mathrm{A}),\\
    \dot u_\mathrm{C} &= u_\mathrm{C}^2 - (1+\rho\gamma)(u_\mathrm{C}-u_\mathrm{A})(u_\mathrm{C}-u_\mathrm{B}).
  \end{align*}
\end{mycorollary}
\smallskip
This is a so-called generalized Darboux--Halphen (gDH) system of ODEs
\cite{Ablowitz2006,Chakravarty2010,Harnad2000}.  It is evident that the gDH
system with
$(\alpha,\beta,\gamma)=(1/n_\mathcal{A},1/n_\mathcal{B},1/n_\mathcal{C})$
and $\rho=\frac2{1-\alpha-\beta-\gamma}$ arises naturally from the unique
(up~to conjugacy) triangle subgroup of~${\it PSL}(2,\mathbf{R})$ with
signature $(n_\mathcal{A},n_\mathcal{B},n_\mathcal{C})$.  Examples of gDH
systems coming from modular subgroups have appeared in the literature;
e.g., the ones coming from the six (up~to conjugacy) triangle subgroups of
$\Gamma(1)={\it PSL}(2,\mathbf{Z}),$ which are incidentally the only gDH
systems for which some linear combination of
$u_\mathrm{A},u_\mathrm{B},u_\mathrm{C}$ satisfies the classical Chazy
equation~\cite{Chakravarty2010}.  However, the general statement is new.

\subsection{Triangle groups commensurable with $\Gamma(1)$}
\label{subsec:final2}

The triangle subgroups of ${\it PSL}(2,\mathbf{R})$ commensurable with
$\Gamma(1)={\it PSL}(2,\mathbf{Z})$ are well known.  (Subgroups
$\Gamma_1,\Gamma_2$ are said to be commensurable if $\Gamma_1\cap\Gamma_2$
is of finite index in both.)  Up~to conjugacy there are exactly
nine~\cite{Takeuchi77}, listed in Table~\ref{tab:3}.  Each is hyperbolic
with at~least one cusp.  They are of three types, and it will be shown that
to each type there is an associated differential system, parametrized
by~$\rho,$ which is satisfied by weight-$1$ forms
$\mathcal{A}\!,\,\mathcal{B},\mathcal{C}$.  For the first type the system
will be that of Theorem~\ref{thm:main}(2), in which $\rho$~equals the
signature of Ramanujan's elliptic theories (i.e., $\rho=r=4,3,2$).

Type~I comprises $\Gamma_0(N),$ $N=2,3,4,$ and Type~II comprises
$\Gamma(1)$ and the Fricke extensions $\Gamma^+_0(N),$ $N=2,3,$ which are
not subgroups of~$\Gamma(1)$.  Type~III comprises three groups that will be
called $\mathrm{2a}',\mathrm{4a}',\mathrm{6a}'$.  The group~$\mathrm{2a}'$
is the index-$2$ subgroup $\Gamma^2<\Gamma(1),$ but the latter two are not
subgroups of~$\Gamma(1)$.  These names are taken from Harnad and McKay
\cite{Harnad2000}\footnote{The primes on
$\mathrm{2a}',\mathrm{4a}',\mathrm{6a}'$ indicate a relation to the groups
labeled $\mathrm{2a},\mathrm{4a},\mathrm{6a}$ in the extended
Conway--Norton classification.  In the Conway--Norton notation used
in~\cite{Harnad2000}, the Type~I groups
$\Gamma_0(2),\Gamma_0(3),\Gamma_0(4)$ are referred~to as
$\mathrm{2B},\mathrm{3B},\mathrm{4C},$ and the Type~II groups
$\Gamma(1),\Gamma_0^+(2),\Gamma_0^+(3)$
as~$\mathrm{1A},\mathrm{2A},\mathrm{3A}$.}.  It is known that the
intersections of the groups $\mathrm{4a}',\mathrm{6a}'$ with~$\Gamma(1)$
are $\Gamma_0(4)\cap\Gamma(2),$ $\Gamma_0(3)\cap\Gamma(2),$ which are
conjugates of $\Gamma_0(8),\allowbreak\Gamma_0(12)$ under
$\tau\mapsto2\tau$.

For each triangle group~$\Gamma,$ the exponents
$(\alpha,\beta,\gamma\,)=(1/n_{\mathcal{A}},1/n_{\mathcal{B}},1/n_{\mathcal{C}}\!)$
are given; as is $\rho=2/(1-\alpha-\beta-\gamma),$ which subsumes the
signature of Ramanujan's theories.  For concreteness, generators
$a,b,c\in\Gamma$ of corresponding stabilizing subgroups are given as~well.
(What are given are homogeneous versions $\bar a,\bar b,\bar c\in{\it
GL}(2,\mathbf{R})$.  To~convert them to $a,b,c\in{\it PSL}(2,\mathbf{R}),$
satisfying $cab=\pm{I},$ divide each by its determinant and prepend~$\pm$.)
These generators are adapted from~\cite[Table~1]{Harnad2000}.  By
convention, $n_{\mathcal{C}}=\infty$ and $\gamma=0,$ and the corresponding
class~$\mathrm{C}$ of fixed points (cusps) includes the infinite cusp.  The
width of the infinite cusp is denoted~$\upsilon,$ as above.

The expressions for the formal (i.e., potentially multivalued) weight-$1$
modular forms $\mathcal{A}\!,\mathcal{B},\mathcal{C},$ which in each case
satisfy ${\mathcal{A}}^\rho = {\mathcal{B}}^\rho + {\mathcal{C}}^\rho\!,\,$
come as follows.  
Starting with the ones for $\Gamma$ of Type~I and for
$\Gamma=\Gamma(1),$ which have already been discussed, they are obtained as
consequences of the index-$2$ subgroup relations~\cite{Takeuchi77a}
\begin{gather}
\label{eq:gather21}
  \Gamma_0(2) <_2\, \Gamma_0^+(2), \qquad   \Gamma_0(3) <_2\, \Gamma_0^+(3); \\
  \mathrm{2a}' <_2\, \Gamma(1), \qquad \mathrm{4a}' <_2\, \Gamma^+_0(2),  \qquad \mathrm{6a}' <_2\, \Gamma^+_0(3).
\label{eq:gather22}
\end{gather}
For instance, suppose that $\tilde\Gamma<_2\Gamma$ and that the respective
signatures satisfy $2\tilde\rho=\rho$.  The goal is to relate the
associated triples
$\tilde{\mathcal{A}}\!,\tilde{\mathcal{B}},\tilde{\mathcal{C}}$ and
$\mathcal{A}\!,\mathcal{B},\mathcal{C}$.  Let the corresponding classes of
fixed points be $\tilde{\mathrm{A}},\tilde{\mathrm{B}},\tilde{\mathrm{C}}$
and $\mathrm{A},\mathrm{B},\mathrm{C}$.  Suppose that
$\tilde{\mathrm{B}},\tilde{\mathrm{C}}$ are cusp classes
under~$\tilde\Gamma,$ which merge into a single class~$\mathrm{C}$
under~$\Gamma,$ but that $\mathrm{A}=\tilde{\mathrm{A}}$.  This is
precisely what happens when
$(\tilde\Gamma,\Gamma)=(\Gamma_0(N),\Gamma_0^+(N))$ for $N=2,3$.

Under these assumptions, one will have $\mathcal{A}=\tilde{\mathcal{A}}$.
Expressions for $\mathcal{B},\mathcal{C}$ in~terms of
$\tilde{\mathcal{B}},\tilde{\mathcal{C}}$ come from a Hauptmodul relation.
Hauptmoduls for $\tilde\Gamma,\Gamma,$ i.e., rational parameters on the
quotient curves
$\mathfrak{H}^*\setminus\tilde\Gamma,\,\mathfrak{H}^*\setminus\Gamma,$ will
be
$\tilde\lambda=\tilde{\mathcal{C}}^{\tilde\rho}/\tilde{\mathcal{A}}^{\tilde\rho},$
resp.\ $\lambda={\mathcal{C}}^{\rho}/{\mathcal{A}}^{\rho}$.  The index-$2$
relation $\tilde\Gamma<_2\Gamma$ induces a double covering
${\mathfrak{H}^*\setminus\tilde\Gamma}\to\allowbreak\mathfrak{H}^*\setminus\Gamma,$
i.e.,
$\mathbf{P}^1(\mathbf{C})_{\tilde\lambda}\to\mathbf{P}^1(\mathbf{C})_{\lambda},$
i.e., a quadratic rational map $\tilde\lambda\mapsto\lambda$.  Since
$\tilde\lambda=1,0$ (corresponding to
classes~$\tilde{\mathrm{B}},\tilde{\mathrm{C}}$) must be mapped to
$\lambda=0$ (corresponding to class~$\mathrm{C}$), and
$\tilde\lambda=\infty$ (corresponding to class~$\tilde{\mathrm{A}}$) must
be mapped to $\lambda=\infty$ (corresponding to class~${\mathrm{A}}$), the
map must be
\begin{equation}
  \lambda = 4\,\tilde\lambda(1-\tilde \lambda) = 1 - (1-2\,\tilde\lambda)^2.
\end{equation}

\begin{landscape}
\ 
\vfill
\begin{table}[h]
\caption{For each triangle subgroup $\Gamma<{\it PSL}(2,\mathbf{R})$
  commensurable with $\Gamma(1)={\it PSL}(2,\mathbf{Z}),$ the basic data
  and the triple $\mathcal{A}\!,\mathcal{B},\mathcal{C}$ of (possibly
  multivalued) weight-$1$ modular forms, satisfying
  $\mathcal{A}^\rho=\mathcal{B}^\rho + \mathcal{C}^\rho$.  The nine
  subgroups are partitioned into Types I,II,III\null.  If
  $n_{\mathcal{A}}<\infty,$ resp.\ $n_{\mathcal{B}}<\infty,$ then the
  minimum power of $\mathcal{A},$ resp.~$\mathcal{B},$ which is
  single-valued, equals the numerator of $\rho\alpha,$ resp.~$\rho\beta,$
  expressed in lowest terms.
  The forms $\mathcal{A}\!,\mathcal{B}$ on~$\mathrm{2a}'=\Gamma^2$ can
  alternatively be written as $({\mathscr{B}_2}^2 -
  \bar\zeta_3{\mathscr{C}_2}^2)^{1/2}(q^{1/2})$ and $({\mathscr{B}_2}^2 -
  \zeta_3{\mathscr{C}_2}^2)^{1/2}(q^{1/2}),$ where $\zeta_3=\exp(2\pi{\rm
  i}/3)$; and the forms $\mathcal{A}\!,\mathcal{B}$ on~$\mathrm{4a}'$ as
  $\mathscr{A}_2\pm{\rm i}\,\mathscr{C}_2$.}
\begin{center}
{\small
\begin{tabular}{|l||c|c|c|c|c|c|l||l|l|l|}
\hline
$\Gamma$ & $(n_{\mathcal{A}},n_{\mathcal{A}},n_{\mathcal{A}})$ & $(\alpha,\beta,\gamma)$ &
$\rho$ & $\bar a$ & $\bar b$ & $\bar c$ & $\upsilon$ & $\mathcal{A}$ &
$\mathcal{B}$ & $\mathcal{C}$ \\
\hline
\hline
$\Gamma_0(2)$ & $(2,\infty,\infty)$ & $(\frac12,0,0)$ & 
$4$ & 
$\vphantom{{\Bigl[}^{1/6}}
\left(\begin{smallmatrix}1&-1\\2&-1\end{smallmatrix}\right)$ &
$\left(\begin{smallmatrix}-1&0\\2&-1\end{smallmatrix}\right)$ &
$\left(\begin{smallmatrix}1&1\\0&1\end{smallmatrix}\right)$ &
$1$ &
$\mathscr{A}_4$ &
$\mathscr{B}_4$ &
$\mathscr{C}_4$ \\
\hline
$\Gamma_0(3)$ & $(3,\infty,\infty)$ & $(\frac13,0,0)$ & 
$3$ & 
$\vphantom{{\Bigl[}^{1/6}}
\left(\begin{smallmatrix}1&-1\\3&-2\end{smallmatrix}\right)$ &
$\left(\begin{smallmatrix}-1&0\\3&-1\end{smallmatrix}\right)$ &
$\left(\begin{smallmatrix}1&1\\0&1\end{smallmatrix}\right)$ &
$1$ &
$\mathscr{A}_3$ &
$\mathscr{B}_3$ &
$\mathscr{C}_3$ \\
\hline
$\Gamma_0(4)$ & $(\infty,\infty,\infty)$ & $(0,0,0)$ & 
$2$ & 
$\vphantom{{\Bigl[}^{1/6}}
\left(\begin{smallmatrix}1&-1\\4&-3\end{smallmatrix}\right)$ &
$\left(\begin{smallmatrix}-1&0\\4&-1\end{smallmatrix}\right)$ &
$\left(\begin{smallmatrix}1&1\\0&1\end{smallmatrix}\right)$ &
$1$ &
$\mathscr{A}_2$ &
$\mathscr{B}_2$ &
$\mathscr{C}_2$ \\
\hline
\hline
$\Gamma(1)$ & $(3,2,\infty)$ & $(\frac13,\frac12,0)$ & 
$12$ & 
$\vphantom{{\Bigl[}^{1/6}}
\left(\begin{smallmatrix}0&-1\\1&-1\end{smallmatrix}\right)$ &
$\left(\begin{smallmatrix}0&-1\\1&0\end{smallmatrix}\right)$ &
$\left(\begin{smallmatrix}1&1\\0&1\end{smallmatrix}\right)$ &
$1$ &
${E_4}^{1/4}$ &
${E_6}^{1/6}$ &
$(12^3\Delta)^{1/12}$ \\
\hline
$\Gamma_0^+(2)$ & $(4,2,\infty)$ & $(\frac14,\frac12,0)$ & 
$8$ & 
$\vphantom{{\Bigl[}^{1/6}}
\left(\begin{smallmatrix}0&-1\\2&-2\end{smallmatrix}\right)$ &
$\left(\begin{smallmatrix}0&-1\\2&0\end{smallmatrix}\right)$ &
$\left(\begin{smallmatrix}1&1\\0&1\end{smallmatrix}\right)$ &
$1$ &
$\mathscr{A}_4$ &
$\left[{{\mathscr{B}_4}^4 -{\mathscr{C}_4}^4}\right]^{1/4}$ &
$2^{1/4}\sqrt{\mathscr{B}_4\mathscr{C}_4}$ \\
\hline
$\Gamma_0^+(3)$ & $(6,2,\infty)$ & $(\frac16,\frac12,0)$ & 
$6$ & 
$\vphantom{{\Bigl[}^{1/6}}
\left(\begin{smallmatrix}0&-1\\3&-3\end{smallmatrix}\right)$ &
$\left(\begin{smallmatrix}0&-1\\3&0\end{smallmatrix}\right)$ &
$\left(\begin{smallmatrix}1&1\\0&1\end{smallmatrix}\right)$ &
$1$ &
$\mathscr{A}_3$ &
$\left[{{\mathscr{B}_3}^3 -{\mathscr{C}_3}^3}\right]^{1/3}$ &
$2^{1/3}\sqrt{\mathscr{B}_3\mathscr{C}_3}$ \\
\hline
\hline
$\mathrm{2a}'=\Gamma^2$ & $(3,3,\infty)$ & $(\frac13,\frac13,0)$ & 
$6$ & 
$\vphantom{{\Bigl[}^{1/6}}
\left(\begin{smallmatrix}1&-3\\1&-2\end{smallmatrix}\right)$ &
$\left(\begin{smallmatrix}0&-1\\1&-1\end{smallmatrix}\right)$ &
$\left(\begin{smallmatrix}1&2\\0&1\end{smallmatrix}\right)$ &
$2$ &
$\left[E_6+{\rm i}\sqrt{12^3\Delta}\right]^{1/6}$ &
$\left[E_6-{\rm i}\sqrt{12^3\Delta}\right]^{1/6}$ &
$(2{\rm i})^{1/6}(12^3\Delta)^{1/12}$ \\
\hline
$\mathrm{4a}'$ & $(4,4,\infty)$ & $(\frac14,\frac14,0)$ & 
$4$ & 
$\vphantom{{\Bigl[}^{1/6}}
\left(\begin{smallmatrix}2&-5\\2&-4\end{smallmatrix}\right)$ &
$\left(\begin{smallmatrix}0&-1\\2&-2\end{smallmatrix}\right)$ &
$\left(\begin{smallmatrix}1&2\\0&1\end{smallmatrix}\right)$ &
$2$ &
$\left[{\mathscr{B}_4}^2 + {\rm i}\,{\mathscr{C}_4}^2 \right]^{1/2}$ &
$\left[{\mathscr{B}_4}^2 - {\rm i}\,{\mathscr{C}_4}^2 \right]^{1/2}$ &
$(4{\rm i})^{1/4}\sqrt{{\mathscr{B}}_4{\mathscr{C}}_4}$ \\
\hline
$\mathrm{6a}'$ & $(6,6,\infty)$ & $(\frac16,\frac16,0)$ & 
$3$ & 
$\vphantom{{\Bigl[}^{1/6}}
\left(\begin{smallmatrix}3&-7\\3&-6\end{smallmatrix}\right)$ &
$\left(\begin{smallmatrix}0&-1\\3&-3\end{smallmatrix}\right)$ &
$\left(\begin{smallmatrix}1&2\\0&1\end{smallmatrix}\right)$ &
$2$ &
$\left[{\mathscr{B}_3}^{3/2} + {\rm i}\,{\mathscr{C}_{3}}^{3/2} \right]^{2/3}$ &
$\left[{\mathscr{B}_3}^{3/2} - {\rm i}\,{\mathscr{C}_{3}}^{3/2} \right]^{2/3}$ &
$(4{\rm i})^{1/3}\sqrt{{\mathscr{B}}_3{\mathscr{C}_3}}$ \\
\hline
\end{tabular}
}
\end{center}
\label{tab:3}
\end{table}
\vfill
\ 
\end{landscape}

\noindent
Here, the proportionality constant (i.e., $4$) is determined by the
condition that ${\lambda=1},$ corresponding to the class~$\mathrm{B}$ of
fixed points under~$\Gamma,$ should be a critical value of the map.  If it
were not, then $\mathrm{B}$~would also be such a class under
$\tilde\Gamma$; which would violate the assumption that $\tilde\Gamma$ is a
triangle group, with only three such classes.

Using the above expressions for $\tilde\lambda,\lambda,$ and also the
identities $\tilde{\mathcal{A}}^{\tilde\rho} =
\tilde{\mathcal{B}}^{\tilde\rho} + \tilde{\mathcal{C}}^{\tilde\rho}\!,\,$
${\mathcal{A}}^\rho = {\mathcal{B}}^\rho + {\mathcal{C}}^\rho\!,\,$ with
$\rho=2\tilde\rho,$ one immediately obtains
\begin{subequations}
\begin{align}
\mathcal{A}&=\tilde{\mathcal{A}},\\
\mathcal{B}&= \sqrt[\tilde\rho]{\tilde{\mathcal{B}}^{\tilde\rho} - \tilde{\mathcal{C}}^{\tilde\rho}},\\
\mathcal{C}&=2^{1/\tilde\rho}\sqrt{\tilde{\mathcal{B}}\tilde{\mathcal{C}}}.
\end{align}
\end{subequations}
Applied to the pairs $(\tilde\Gamma,\Gamma)$ of~(\ref{eq:gather21}), these
yield the triples $\mathcal{A}\!,\mathcal{B},\mathcal{C}$ for
$\Gamma=\Gamma_0^+(2),\Gamma_0^+(3)$ that are shown in Table~\ref{tab:3}.
A~similar but `reversed' procedure, applied to the $(\tilde\Gamma,\Gamma)$
of~(\ref{eq:gather22}), allows the triples
$\mathcal{A}\!,\mathcal{B},\mathcal{C}$ for the Type-III groups
$\mathrm{2a}',\mathrm{4a}',\mathrm{6a}',$ to be computed in~terms of those
for the corresponding Type-II groups
$\Gamma(1),\Gamma_0^+(2),\Gamma_0^+(3)$.  The resulting triples are given
in the table.

Alternative representations for the forms $\mathcal{A}\!,\mathcal{B}$ on
the group $\mathrm{2a}'=\Gamma^2$ are supplied in the caption, and are
derived as follows.  Although these forms are not single-valued, their
squares ${\mathcal{A}}^{2},{\mathcal{B}}^2$ are single-valued, by the test
for single-valuedness mentioned immediately before
Theorem~\ref{thm:triangular} (and reproduced in the caption).  Each of
${\mathcal{A}}^{2},{\mathcal{B}}^2$ has a $\{1,\zeta_3,\zeta_3^2\}$-valued
multiplier system, and since $\Gamma^2$~has as index-$3$ subgroup the
principal modular subgroup~$\Gamma(2),$ each of them lies in
$M_2(\Gamma(2))$.  But $\Gamma(2)$~is is conjugated to~$\Gamma_0(4)$
by~$\tau\mapsto2\tau$.  Since ${{\mathscr{B}_2}}^2,{{\mathscr{C}_2}}^2$
span~$M_2(\Gamma_0(4)),$ the forms ${\mathcal{A}}^{2},{\mathcal{B}}^2$ must
be combinations of
${{\mathscr{B}_2}}^2(q^{1/2}),{{\mathscr{C}_2}}^2(q^{1/2}),$ i.e., of
${\vartheta_4}^4(q^{1/2}),{\vartheta_2}^4(q^{1/2})$.  The combinations are
easily worked~out by linear algebra, if one expands to second order in
$q_2=q^{1/2}$.

\subsection{Explicit systems and Chazy equations}
\label{subsec:final3}

For each of the nine (conjugacy classes of) triangle groups~$\Gamma$
commensurable with $\Gamma(1)={\it PSL}(2,\mathbf{Z}),$ the associated
differential system and generalized Chazy equation are computed below.
They come respectively from Theorems \ref{thm:triangular}
and~\ref{thm:general}.  For each~$\Gamma,$ a hypergeometric (i.e.,
elliptic-integral) representation of the corresponding weight-$1$
form~$\mathcal{A},$ coming from Eq.~(\ref{eq:rewritable2}a), is given
as~well.

As was explained in~\S\,\ref{subsec:final2}, these triangle subgroups are
of three types, denoted I,II,III\null.  From a classical-analytic rather
than a modular point of view, they differ in the dependence of the exponent
differences $(\alpha,\beta,\gamma)$ on the signature~$\rho$.

\begin{itemize}
  \item 
    Type I, for which $(\alpha,\beta,\gamma)=(1-\tfrac2\rho,0,0)$.  It
    comprises $\Gamma=\Gamma_0(2),\Gamma_0(3),\Gamma_0(4),$ for which
    $\rho=4,3,2$; in each case the infinite cusp has width $\upsilon=1$.
    For each of these groups the associated triple
    ${\mathcal{A}}_\rho,{\mathcal{B}}_\rho,{\mathcal{C}}_\rho$ of
    weight-$1$ forms equals
    ${\mathscr{A}}_\rho,{\mathscr{B}}_\rho,{\mathscr{C}}_\rho,$ and the
    weight-$2$ quasi-modular form ${\mathcal{E}}_\rho := \upsilon
    ({{\mathcal{C}}_\rho}^\rho)'/({{\mathcal{C}}_\rho}^\rho)$
    equals ${\mathscr{E}}_\rho$.  The system satisfied by
    ${\mathscr{A}}_\rho,{\mathscr{B}}_\rho,{\mathscr{C}}_\rho;{\mathscr{E}}_\rho$
    and the generalized Chazy equation satisfied by $u=(2\pi{\rm
    i}/\rho){\mathscr{E}}_\rho$ were given in Theorem~\ref{thm:main}.  

    By Eq.~(\ref{eq:rewritable2}a), the hypergeometric representation for
    ${\mathcal{A}}_\rho={\mathscr{A}}_\rho$ is ${\mathscr{A}}_\rho =
    \hat{\mathsf{K}}_\rho^{\mathrm{I}}(\lambda_\rho\!),$ where
    $\lambda_\rho:={{\mathscr{C}}_\rho}^\rho / {{\mathscr{A}}_\rho}^\rho$
    is a Hauptmodul for~$\Gamma$ and
  \begin{align}
    \hat{\mathsf{K}}_\rho^{\mathrm{I}}(\lambda_\rho\!) 
    :&= {}_2F_1(\tfrac1\rho,1-\tfrac1\rho;\,1;\,\lambda_\rho\!)\\
    &= \frac{\sin(\pi/\rho)}\pi \int_{0}^{1} x^{-1/\rho}  (1-x)^{-1+1/\rho}
    (1-\lambda_\rho x)^{-1/\rho}\,dx
    \nonumber
  \end{align}
    is the (normalized) Type-I complete elliptic integral.  These cases of
    Type~I correspond to Ramanujan's elliptic theories of signature~$\rho,$
    for $\rho=4,3,2$ (see~\cite{Berndt95,Borwein91}).  The classical
    (Jacobi) case is $\rho=2,$ and $\hat{\mathsf{K}}_2^{\mathrm{I}}$~is the
    (normalized) complete integral~$\hat{\mathsf{K}},$ which was introduced
    in~\S\,\ref{sec:theta}.
  \item
    Type II, for which
    $(\alpha,\beta,\gamma)=(\tfrac12-\tfrac2\rho,\tfrac12,0)$.  It
    comprises
    $\Gamma=\Gamma(1),\allowbreak\Gamma_0^+(2),\allowbreak\Gamma_0^+(3),$
    for which $\rho=12,8,6$; in each case the infinite cusp has width
    $\upsilon=1$.  The associated triples
    ${\mathcal{A}}_\rho,{\mathcal{B}}_\rho,{\mathcal{C}}_\rho$ of
    weight-$1$ forms are in Table~\ref{tab:3}.  By direct
    computation, the weight-$2$ quasi-modular forms
    $\mathcal{E}_\rho:=\upsilon
    ({\mathcal{C}_\rho}^\rho)'/({\mathcal{C}_\rho}^\rho)$ are
      \begin{alignat}{2}
	\mathcal{E}_{12}(q) &= \vphantom{\bigl[\bigr]}E_2(q)\label{eq:e12}\\
	                 &= 1 -24 \sum_{n=1}^\infty \sigma_1(n;1)q^n& &= 1 -24\sum_{n=1}^\infty \sigma^{\mathrm{c}}_1(n;1)q^n, \notag\\
	\mathcal{E}_8(q) &= \tfrac13\bigl[2\,E_2(q^2) + E_2(q)\bigr]\label{eq:e8}\\
	                 &= 1 - 8\sum_{n=1}^\infty \sigma_1(n;2,1)q^n& &=1 -8\sum_{n=1}^\infty \sigma^{\mathrm{c}}_1(n;3,1)q^n,\notag\\
	\mathcal{E}_6(q) &= \tfrac14\bigl[3\,E_2(q^3) + E_2(q)\bigr]\label{eq:e6}\\
	                 &= 1 -6\sum_{n=1}^\infty \sigma_1(n;2,1,1)q^n& &= 1 -6\sum_{n=1}^\infty \sigma^{\mathrm{c}}_1(n;4,1,1)q^n.\notag
      \end{alignat}
    They lie respectively in $M_2^{\le1}(\Gamma(1)),$
    $M_2^{\le1}(\Gamma_0^+(2)),$ $M_2^{\le1}(\Gamma_0^+(3))$.  By
    Theorem~\ref{thm:triangular}, the differential system parametrized
    by~$\rho$ is
    \begin{subequations}
      \label{eq:PQR2full}
      \begin{align}
	({{\mathcal{A}}_\rho}^{\rho})' &= {\mathcal{E}}_\rho\cdot {{\mathcal{A}}_\rho}^{\rho} -
	{{\mathcal{A}}_\rho}^{\rho/2+2} {{\mathcal{B}}_\rho}^{\rho/2},\\
	({{\mathcal{B}}_\rho}^{\rho})' &= \mathcal{E}_\rho\cdot {{\mathcal{B}}_\rho}^{\rho} -
	{{\mathcal{A}}_\rho}^{\rho/2+2} {{\mathcal{B}}_\rho}^{\rho/2},\\
	({{\mathcal{C}}_\rho}^{\rho})' &= \mathcal{E}_\rho\cdot {{\mathcal{C}}_\rho}^{\rho},\\
	\rho\,\mathcal{E}_\rho' &= \mathcal{E}_\rho\cdot\mathcal{E}_\rho - {{\mathcal{A}_\rho}}^4.
      \end{align}
    \end{subequations}
    This is an extension of Ramanujan's $P$--$Q$--$R$ system
    (\ref{eq:PQR2}abcd), to which it reduces when $\rho=12$
    and~$\Gamma=\Gamma(1)$.  

    The generalized Chazy equation satisfied by $u=(2\pi{\rm
    i}/\upsilon\rho){\mathcal{E}}_\rho=\dot{\mathcal{C}}_\rho/\mathcal{C}_\rho,$
    according to Theorem~\ref{thm:general}, is the nonlinear third-order
    ODE $C_{p_\rho},$ i.e., ${p_\rho}=0,$ in which the polynomial
    $p_\rho\in\mathbf{C}[u_4,u_6,u_8]$ is defined by
    \begin{align}
      p_{12} &= u_8 + 24\,u_4^2,\label{eq:newchazy1}\\
      p_8 &= 2\,u_4u_8 - u_6^2 + 32\,u_4^3  ,\label{eq:newchazy2}\\
      p_6 &= 4\,u_4u_8 - 3\,u_6^2 + 48\,u_4^3,\label{eq:newchazy3}
    \end{align}
    and $u_4,u_6,u_8$ were given in Definition~\ref{def:u468}.  The
    differential equation $C_{p_{12}}$ associated to~$\Gamma(1),$ coming
    from~(\ref{eq:newchazy1}), is the classical Chazy
    equation~(\ref{eq:chazy}) that is satisfied by $u=(2\pi{\rm i}/12)E_2$.
    The polynomials (\ref{eq:newchazy2}),(\ref{eq:newchazy3}) yield the
    generalized Chazy equations associated to
    $\Gamma_0^+(2),\allowbreak\Gamma_0^+(3),$ which are new.

    By Eq.~(\ref{eq:rewritable2}a), the hypergeometric representation for
    the weight-$1$ form ${\mathcal{A}}_\rho$ is ${\mathcal{A}}_\rho =
    \hat{\mathsf{K}}^{\mathrm{II}}_\rho(\lambda_\rho\!),$ in which
    $\lambda_\rho:={{\mathcal{C}}_\rho}^\rho / {{\mathcal{A}}_\rho}^\rho$
    is a Hauptmodul for~$\Gamma$ and
  \begin{subequations}
    \label{eq:lastminute}
  \begin{align}
    \hat{\mathsf{K}}_\rho^{\mathrm{II}}(\lambda_\rho\!) 
    :&= {}_2F_1(\tfrac1\rho,\tfrac12-\tfrac1\rho;\,1;\,\lambda_\rho\!)\\
    &= \frac{\cos(\pi/\rho)}\pi \int_{0}^{1} x^{-1/2-1/\rho}  (1-x)^{-1/2+1/\rho}
    (1-\lambda_\rho x)^{-1/\rho}\,dx    \notag
  \end{align}
    is the (normalized) Type-II complete elliptic integral.  Equivalently,
    \begin{equation}
    \left[\hat{\mathsf{K}}_\rho^{\mathrm{II}}(\lambda_\rho\!)\right]^2
     = {}_3F_2(\tfrac2\rho,\tfrac12,1-\tfrac2\rho;\,1,1;\,\lambda_\rho\!).
    \end{equation}
  \end{subequations}
    Such representations, when $\rho=12$ and~$\Gamma=\Gamma(1),$ are fairly
    well known.  By Table~\ref{tab:3}, the $\rho=12$ versions
    of~(\ref{eq:lastminute}ab) are
    \begin{subequations}
    \begin{align}
      \label{eq:112512}
      {E_4}^{1/4} &= {}_2F_1(\tfrac1{12},\tfrac5{12};\,1;\,12^3/j),\\
      {E_4}^{1/2} &= {}_3F_2(\tfrac1{6},\tfrac1{2},\tfrac56;\,1,1;\,12^3/j),
    \end{align}
    \end{subequations}
    where $j={E_4}^3/\Delta=12^3{E_4}^3/({E_4}^3-{E_6}^2)$ is the
    Klein--Weber invariant, the canonical Hauptmodul for~$\Gamma(1),$ so
    that $12^3/j = ({E_4}^3-{E_6}^2)/{E_4}^3$.  Equation (\ref{eq:112512})
    was known to Dedekind and was rediscovered by Stiller~\cite{Stiller88}.
    These identities hold in a neighborhood of the infinite cusp, at which
    $j=\infty$ and $12^3/j=0$.  In the same way,
    Eq.~(\ref{eq:rewritable2}b) yields
    \begin{equation}
      {E_6}^{1/6} = {}_2F_1\left(\tfrac1{12},\tfrac7{12};\,1;\,12^3/(12^3-j)\right).
    \end{equation}
    From~(\ref{eq:rewritable3}a) one also has
    \begin{subequations}
    \begin{align}
      \Delta^{1/12} & \varpropto {}_2F_1(\tfrac1{12},\tfrac1{12};\,\tfrac23;\,j/12^3),\\
      \Delta^{1/6} & \varpropto {}_3F_2(\tfrac1{6},\tfrac1{6},\tfrac16;\,\tfrac13,\tfrac23;\,j/12^3),
    \end{align}
    \end{subequations}
    which hold near any cubic elliptic fixed point, where $j/12^3=0$.
    (E.g., near $\tau=\zeta_3=\exp(2\pi{\rm i}/3)$.)  The constants of
    proportionality depend on the choice of fixed point.

    The $\rho=8,6$ representations, for
    $\Gamma=\Gamma_0^+(2),\Gamma_0^+(3),$ were derived by Zudilin
    \cite[Eqs.~(23bc)]{Zudilin2003}.  The corresponding differential
    systems that he obtained~\cite[Props.\ 6,7]{Zudilin2003} are equivalent
    to the $\rho=8,6$ cases of the system (\ref{eq:PQR2full}abcd), but are
    more complicated as they are not expressed in~terms of weight-$1$
    forms.
  \item
    Type III, for which
    $(\alpha,\beta,\gamma)=(\tfrac12-\tfrac1\rho,\tfrac12-\tfrac1\rho,0)$.
    It comprises
    $\Gamma=(\mathrm{2a}'=\Gamma^2),\allowbreak\mathrm{4a}',\allowbreak\mathrm{6a}',$
    for which $\rho=6,4,3$; in each case the infinite cusp has width
    $\upsilon=2$.  The associated triples
    ${\mathcal{A}}_\rho,{\mathcal{B}}_\rho,{\mathcal{C}}_\rho$ of
    weight-$1$ forms are in Table~\ref{tab:3}.  The weight-$2$
    quasi-modular forms $\mathcal{E}_\rho:=\upsilon
    ({\mathcal{C}_\rho}^\rho)'/({\mathcal{C}_\rho}^\rho),$ $\rho=6,4,3,$
    are identical to the Type-II forms
    $\mathcal{E}_{12},\mathcal{E}_{8},\mathcal{E}_{6},$ given in Eqs.\
    (\ref{eq:e12}),(\ref{eq:e8}),(\ref{eq:e6}).

    By Theorem~\ref{thm:triangular}, the differential system parametrized
    by~$\rho$ is
    \begin{subequations}
      \begin{align}
	2({{\mathcal{A}}_\rho}^{\rho})' &= {\mathcal{E}}_\rho\cdot {{\mathcal{A}}_\rho}^{\rho} -
	{{\mathcal{A}}_\rho}^{\rho/2+1} {{\mathcal{B}}_\rho}^{\rho/2+1},\\
	2({{\mathcal{B}}_\rho}^{\rho})' &= \mathcal{E}_\rho\cdot {{\mathcal{B}}_\rho}^{\rho} -
	{{\mathcal{A}}_\rho}^{\rho/2+1} {{\mathcal{B}}_\rho}^{\rho/2+1},\\
	2({{\mathcal{C}}_\rho}^{\rho})' &= \mathcal{E}_\rho\cdot {{\mathcal{C}}_\rho}^{\rho},\\
	2\rho\,\mathcal{E}_\rho' &= \mathcal{E}_\rho\cdot\mathcal{E}_\rho - {{\mathcal{A}_\rho}}^2{{\mathcal{B}_\rho}}^2.
      \end{align}
    \end{subequations}
    Although this system is significantly different from
    (\ref{eq:PQR2full}abcd), the system of Type~II, the resulting
    generalized Chazy equations~$C_{p_\rho},$ $\rho=6,4,3,$ are identical
    to the Type-II equations for $\rho=12,8,6$ (see
    (\ref{eq:newchazy1}),(\ref{eq:newchazy2}),(\ref{eq:newchazy3})).

    By Eq.~(\ref{eq:rewritable2}a), the hypergeometric representation for
    the weight-$1$ form ${\mathcal{A}}_\rho$ is ${\mathcal{A}}_\rho =
    \hat{\mathsf{K}}^{\mathrm{III}}_\rho(\lambda_\rho\!),$ in which
    $\lambda_\rho:={{\mathcal{C}}_\rho}^\rho / {{\mathcal{A}}_\rho}^\rho$
    is a Hauptmodul for~$\Gamma$ and
  \begin{align}
    \hat{\mathsf{K}}_\rho^{\mathrm{III}}(\lambda_\rho\!) 
    :&= {}_2F_1(\tfrac{1}{\rho},\tfrac12;\,1;\,\lambda_\rho\!)\\
    &= \frac1\pi \int_{0}^{1} x^{-1/2}  (1-x)^{-1/2}
    (1-\lambda_\rho x)^{-1/\rho}\,dx
    \nonumber
  \end{align}
    is the (normalized) Type-III complete elliptic integral.  These
    representations are new.  The case $\rho=4,$ i.e.,
    $\Gamma=\mathrm{4a}',$ is especially noteworthy.  It follows from the
    formulas 
    \begin{gather}
      {\mathcal{A}}_4,{\mathcal{B}}_4 = {\mathscr{A}}_2 \pm {\rm
      i}\,{\mathscr{C}}_2 = {\vartheta_3}^2 \pm {\rm i} {\vartheta_2}^2 \\
      {{\mathcal{A}}_4}^4 = {{\mathcal{B}}_4}^4 + {{\mathcal{C}}_4}^4
    \end{gather}
    that when $\rho=4,$ the equation ${\mathcal{A}}_\rho =
    \hat{\mathsf{K}}^{\mathrm{III}}_\rho(\lambda_\rho\!)$ specializes to
    \begin{equation}
      {\vartheta_3}^2 \pm{\rm i}{\vartheta_2}^2 = 
      {}_2F_1\left(\tfrac14,\tfrac12;\,1;\,
      \frac{({\vartheta_3}^2 \pm{\rm i}{\vartheta_2}^2)^4 - ({\vartheta_3}^2 \mp{\rm i}{\vartheta_2}^2)^4}{({\vartheta_3}^2 \pm{\rm i}{\vartheta_2}^2)^4}
      \right).
    \end{equation}
    It is unclear whether this remarkable theta identity has a non-modular
    proof.
\end{itemize}

\subsection{Discussion}
\label{subsec:final4}

The results of~\S\,\ref{subsec:final3} suggest that elliptic integrals of
Types II and~III (parametrized by the signature~$\rho$) deserve further
study, much like the elliptic integrals of Type~I, which are those of
Ramanujan's alternative theories~\cite{Berndt95}.  His theories fit into a
larger framework: one that is larger by a factor of three, at~least.

The new generalized Chazy equations $C_p$ are especially interesting, since
they open a `modular window' into the space of nonlinear third-order ODEs.
Each of the generalized Chazy equations derived in this article can be
integrated in closed form in~terms of modular functions.  This has
taxonomic ramifications.  The classical Chazy equation $u_8 + 24u_4^2=0$
has the \emph{Painlev\'e property}, in that its solutions have no~`movable
branch points'~\cite[\S\,7.1.5 and Ex.~6.5.14]{Ablowitz91}.  The nonlinear
third-order ODEs in~$u$ which have this property, and in~which $\dddot
u$~is polynomial in $\ddot u,\dot u,u$ and rational in~$x,$ were classified
by Chazy~\cite{Chazy11} into classes numbered I~through~XIII\null.  But
to~date, there has been no~extension of his scheme to third-order ODEs with
the property, in~which $\dddot u$~is non-polynomial but \emph{rational}
in~$\ddot u,\dot u,u$.

The equations $C_p$ of this article provide examples.  (For the defining
polynomials~$p\in\mathbf{C}[u_4,u_6,u_8],$ see Eqs.\
(\ref{eq:oldchazy1})--(\ref{eq:oldchazy3}) and
(\ref{eq:newchazy1})--(\ref{eq:newchazy3}).)  Each equation $p=0$ is a
nonlinear third-order ODE satisfied by $u=\dot{\mathcal{C}}/\mathcal{C},$
where $\mathcal{C}$~is the weight-$1$ modular form that vanishes on the
third class of fixed points of a triangle subgroup $\Gamma<{\it
PSL}(2,\mathbf{R})$ with signature
$(n_{\mathcal{A}},n_{\mathcal{B}};\,n_{\mathcal{C}})$.  Other
than~$C_{p_3},$ the rather complicated weight-$20$ ODE coming from
$\Gamma=\Gamma_0(3),$ these nonlinear ODEs lie in a single new class.  With
one seeming exception, each is of the form
\begin{equation}
\label{eq:newclass}
(M-2)u_4 u_8 - (M-3)u_6^2 + 8\,M\,u_4^3=0.
\end{equation}
Equation~(\ref{eq:newclass}) is of weight~$12$ unless $M=3,$ in which case
the $u_6^2$~term drops~out and it reduces to the classical Chazy equation,
of weight~$8$.  By Theorem~\ref{thm:general}, Eq.~(\ref{eq:newclass}) comes
from the triangle groups with signatures $(M,2;\infty)$ and $(M,M;\infty),$
i.e., with vertex angles $(\alpha,\beta;\gamma)$ (expressed in~terms of
$\pi$~radians) equal to $(\tfrac1M,\tfrac12;0)$ or~$(\tfrac1M,\tfrac1M;0)$.
The abovementioned seeming exception is the generalized Chazy equation
attached to $\Gamma_0(2)$ and~$\Gamma_0(4),$ which must be obtained from
Eq.~(\ref{eq:newclass}) by taking a formal $M\to\infty$ limit.  But such
limits are familiar from Chazy's analysis.  For instance, the classical
Chazy equation, which is attached to the groups $\Gamma(1)$ and~$\Gamma^2,$
with respective signatures $(3,2;\infty)$ and $(3,3;\infty),$ is also the
formal $N\to\infty$ limit of the Chazy-XII equation
\begin{equation}
\label{eq:chazyxii}
(N^2-36) u_8 + 24\,N^2\,u_4^2=0,
\end{equation}
which is of weight~$8$ and comes from the triangle groups with signatures
$(3,2;N)$ and $(3,3;N/2)$.  To date, the Chazy-XII class is the only one
that has been given a modular interpretation, e.g., expressed in~terms of
the forms $u_4,u_6,u_8$.

One expects that when an extended Chazy classification is finally
constructed by non-modular, classical-analytic techniques, the
equation~(\ref{eq:newclass}), parametrized by integer~$M,$ will belong to
an additional `modular' class.  Interestingly, it is the limiting
($N\to\infty$) case of a \emph{two-parameter} generalized Chazy equation,
\begin{equation}
  \label{eq:endit}
\left[(M-2)^2N^2 - 4M^2\right]\cdot\left[(M-2)u_4u_8 - (M-3)u_6^2 \right] +
8\, M(M - 2)^2 N^2 \,u_4^3 = 0,
\end{equation}
which is also of weight~$12$ (generically).  By Theorem~\ref{thm:general},
Eq.~(\ref{eq:endit}) comes from the triangle groups with signatures
$(M,2;N)$ and~$(M,M;N/2),$ i.e., with $(\alpha,\beta;\gamma)$ equal to
$(\frac1M,\frac12;\frac1N)$ or~$(\frac1M,\frac1M;\frac2N)$.  The
fundamental domains of the latter triangle groups are hyperbolic isosceles
triangles.  Equation~(\ref{eq:endit}) reduces to Eq.~(\ref{eq:chazyxii}),
the weight-$8$ Chazy--XII equation, when $M=3$.

In a similar way, one can derive a two-parameter generalized Chazy equation
extending~$C_{p_3},$ the weight-$20$ ODE that comes from~$\Gamma_0(3)$.
The extension, also of weight~$20$ (generically), is
\begin{sizemultline}{\small}
\label{eq:weight20}
\bigl[(2M-3)^2 N^2 - 9M^2\bigr]^2\cdot \bigl[(M-2)u_4u_8-(M-3)u_6^2\bigr]u_8\\
{}+12\,MN^2\bigl[(2M-3)^2N^2-9M^2\bigr]\cdot u_4^2\bigl[4(M-2)(2M-3)u_4u_8-(M-3)(5M-9)u_6^2\bigr]\\
{}+576\,M^2(M-2)(2M-3)^2N^4\cdot u_4^5=0.
\end{sizemultline}
By Theorem~\ref{thm:general}, this nonlinear ODE comes from the triangle
group with signature $(3,M;N),$ i.e., with $(\alpha,\beta;\gamma)$ equal to
$(\frac13,\frac1M;\frac1N)$.  It reduces to~$C_{p_3}$ when $M\to\infty$ and
$N\to\infty,$ to the Chazy--XII equation when $M=2,$ and to the classical
Chazy equation when $M=2$ and $N\to\infty,$ or when $M=3$ and~$N\to\infty$.

Writing $u_8,u_6,u_4$ in~terms of $\dddot u,\ddot u,\dot u,u,$ one sees
that like $C_{p_3}$ itself, the ODE (\ref{eq:weight20}), when $M\neq2,$
expresses $\dddot u$ as a degree-$2$ algebraic function of~$\ddot u,\dot
u,u,$ rather than a rational function.  It is of a more general type than
the ODE~(\ref{eq:endit}).

\section*{Appendix: Theta representations and AGM identities}
\renewcommand\thesection{A}
\setcounter{equation}{0}

The functions $\mathscr{A}_r,\mathscr{B}_r,\mathscr{C}_r,$ $r=2,3,4,$
satisfying ${\mathscr{A}_r}^r = {\mathscr{B}_r}^r + {\mathscr{C}_r}^r$ on
the disk $\left| q\right|<1,$ were originally defined by the
Borweins~\cite{Borwein91} as the sums of theta series of certain quadratic
forms, occurring in Ramanujan's theories of elliptic functions to
alternative bases~\cite{Berndt95}.  To~facilitate comparison, several of
their results are restated below in the notation of the present article.
They defined
$(\mathscr{A}_2,\mathscr{B}_2,\mathscr{C}_2)=({\vartheta_3}^2\!,\,
{\vartheta_4}^2\!,\,{\vartheta_2}^2),$
$({\mathscr{A}_4}^2,{\mathscr{B}_4}^2,{\mathscr{C}_4}^2)= ( {\vartheta_2}^4
+ {\vartheta_3}^4\!,\, {\vartheta_4}^4\!,\,
2{\vartheta_2}^2{\vartheta_3}^2),$ and also
\begin{subequations}
\begin{align}
  \mathscr{A}_3(q) &= \sum_{n,m\in\mathbf{Z}} q^{n^2+nm+m^2}   ,\\
  \mathscr{B}_3(q) &= \sum_{n,m\in\mathbf{Z}} \zeta_3^{n-m} q^{n^2+nm+m^2} ,\\
  \mathscr{C}_3(q) &= \sum_{n,m\in\mathbf{Z}} q^{(n+\frac13)^2+(n+\frac13)(m+\frac13)+(m+\frac13)^2},
\end{align}
\end{subequations}
where $\zeta_3$~is a primitive third root of unity.  Their AGM identities
include the quadratic signature-$2$ identities
\begin{subequations}
\begin{align}
  \mathscr{A}_2(q^2) &= \left[(\mathscr{A}_2 + \mathscr{B}_2)/2\right](q),\\
  \mathscr{B}_2(q^2) &= \sqrt{\mathscr{A}_2\mathscr{B}_2}\,(q),\\
  \mathscr{C}_2(q^2) &= \left[(\mathscr{A}_2 - \mathscr{B}_2)/2\right](q),
\end{align}
\end{subequations}
which originated with Jacobi, the quartic signature-$2$ identities
\begin{subequations}
\begin{align}
  \vartheta_3(q^4) &= \left[(\vartheta_3 + \vartheta_4)/2\right](q),\\
  \vartheta_4(q^4) &= \sqrt[4]{({\vartheta_3}^2+{\vartheta_4}^2)(\vartheta_3\vartheta_4)/2}\,\,(q),\\
  \vartheta_2(q^4) &= \left[(\vartheta_3 - \vartheta_4)/2\right](q),
\end{align}
\end{subequations}
the cubic signature-$3$ identities
\begin{subequations}
\begin{align}
  \mathscr{A}_3(q^3) &= \left[(\mathscr{A}_3 + 2\,\mathscr{B}_3)/3\right](q),\\
  \mathscr{B}_3(q^3) &=
  \sqrt[3]{({\mathscr{A}_3}^2+\mathscr{A}_3\mathscr{B}_3 +
  {\mathscr{B}_3}^2)\mathscr{B}_3 / 3}\:\,(q),\\
  \mathscr{C}_3(q^3) &= \left[(\mathscr{A}_3 - \mathscr{B}_3)/3\right](q),
\end{align}
\end{subequations}
and the quadratic signature-$4$ identities
\begin{subequations}
\begin{align}
  {\mathscr{A}_4}^2(q^2) &= \left[({\mathscr{A}_4}^2 + 3\,{\mathscr{B}_4}^2)/4\right](q),\\
  {\mathscr{B}_4}^2(q^2) &= \sqrt{({\mathscr{A}_4}^2+{\mathscr{B}_4}^2){\mathscr{B}_4}^2/2}\:\,(q),\\
  {\mathscr{C}_4}^2(q^2) &= \left[({\mathscr{A}_4}^2 - {\mathscr{B}_4}^2)/4\right](q).
\end{align}
\end{subequations}


\begin{thebibliography}{10}
\providecommand{\url}[1]{{#1}}
\providecommand{\urlprefix}{URL }
\expandafter\ifx\csname urlstyle\endcsname\relax
  \providecommand{\doi}[1]{DOI~\discretionary{}{}{}#1}\else
  \providecommand{\doi}{DOI~\discretionary{}{}{}\begingroup
  \urlstyle{rm}\Url}\fi

\bibitem{Ablowitz2006}
Ablowitz, M.J., Chakravarty, S., Hahn, H.: Integrable systems and modular forms
  of level {$2$}.
\newblock J.~Phys.~A \textbf{39}(50), 15,341--15,353 (2006)

\bibitem{Ablowitz91}
Ablowitz, M.J., Clarkson, P.A.: Solitons, Nonlinear Evolution Equations and
  Inverse Scattering.
\newblock No. 149 in London Mathematical Society Lecture Note Series. Cambridge
  Univ. Press, Cambridge, UK (1991)

\bibitem{Berndt95}
Berndt, B.C., Bhargava, S., Garvan, F.G.: Ramanujan's theories of elliptic
  functions to alternative bases.
\newblock Trans. Amer. Math. Soc. \textbf{347}(11), 4163--4244 (1995)

\bibitem{Berndt99}
Berndt, B.C.: Fragments by {R}amanujan on {L}ambert series.
\newblock In: S.~Kanemitsu, K.~Gyr (eds.) Number Theory and its Applications,
  no.~2 in Dev. Math., pp. 35--49. Kluwer, Dordrecht (1999)

\bibitem{Berndt2006}
Berndt, B.C.: Number Theory in the Spirit of {R}amanujan, \emph{Student
  Mathematical Library}, vol.~34.
\newblock American Mathematical Society (AMS), Providence, RI (2006)

\bibitem{Borwein87}
Borwein, J.M., Borwein, P.B.: Pi and the {AGM}: A Study in Analytic Number
  Theory and Computational Complexity, \emph{Canadian Mathematical Society
  Series of Monographs and Advanced Texts}, vol.~4.
\newblock Wiley, New York (1987)

\bibitem{Borwein91}
Borwein, J.M., Borwein, P.B.: A cubic counterpart of {J}acobi's identity and
  the {AGM}.
\newblock Trans. Amer. Math. Soc. \textbf{323}(2), 691--701 (1991)

\bibitem{Brezhnev2006}
Brezhnev, Y.V.: On functions of {J}acobi and {W}eierstrass {(I)} (2006),
\newblock preprint, available on-line as arXiv:math/0601371

\bibitem{Bureau87}
Bureau, F.J.: Sur des syst{\`e}mes diff{\'e}rentiels non lin{\'e}aires du
  troisi{\`e}me ordre et les {\'e}quations diff{\'e}rentielles non
  lin{\'e}aires associ{\'e}es.
\newblock Acad. Roy. Belg. Bull. Cl. Sci. (5) \textbf{73}(6--9), 335--353
  (1987)

\bibitem{Chakravarty2010}
Chakravarty, S., Ablowitz, M.J.: Parameterizations of the {C}hazy equation.
\newblock Stud. Appl. Math. \textbf{124}(2), 105--135 (2010)

\bibitem{Chazy11}
Chazy, J.: Sur les {\'e}quations diff{\'e}rentielles du troisi{\`e}me et
  d'ordre sup{\'e}rieur dont l'int{\'e}grale g{\'e}n{\'e}rale a ses points
  critiques fixes.
\newblock Acta Math. \textbf{34}, 317--385 (1911)

\bibitem{Clarkson96}
Clarkson, P.A., Olver, P.J.: Symmetry and the {Chazy} equation.
\newblock J.~Differential Equations \textbf{124}(1), 225--246 (1996)

\bibitem{Cohen77}
Cohen, H., Oesterl{\'e}, J.: Dimensions des espaces de formes modulaires.
\newblock In: Modular Functions of One Variable, {VI}, no. 627 in Lecture Notes
  in Mathematics, pp. 69--78. Springer-Verlag, New York/Berlin (1977)

\bibitem{Diamond2005}
Diamond, F., Shurman, J.: A First Course in Modular Forms.
\newblock Springer-Verlag, New York/Berlin (2005)

\bibitem{Ehrenpreis94}
Ehrenpreis, L.: Singularities, functional equations, and the circle method.
\newblock In: G.E. Andrews, D.M. Bressoud, L.A. Parson (eds.) The {R}ademacher
  Legacy to Mathematics, no. 166 in Contemporary Mathematics, pp. 35--80. Amer.
  Math. Soc., Providence, RI (1994)

\bibitem{Enneper1890}
Enneper, A.: Elliptische {F}unctionen: {T}heorie und {G}eschichte, 2nd edn.
\newblock Nebert-Verlag, Halle, Germany (1890)

\bibitem{Fine88}
Fine, N.J.: Basic Hypergeometric Series and Applications.
\newblock No.~27 in Mathematical Surveys and Monographs. Amer. Math. Soc.,
  Providence, RI (1988)

\bibitem{Ford51}
Ford, L.R.: Automorphic Functions, 2nd edn.
\newblock Chelsea Publishing Co., New York (1951)

\bibitem{Glaisher1885b}
Glaisher, J.W.L.: On certain sums of products of quantities depending upon the
  divisors of a number.
\newblock Messenger Math. \textbf{15}, 1--20 (1885)

\bibitem{Glaisher1885a}
Glaisher, J.W.L.: On the quantities {$K,E,J,G,K',E',J',G'$} in elliptic
  functions.
\newblock Quart. J.~Pure Appl. Math. \textbf{20}, 313--361 (1885)

\bibitem{Hahn2008}
Hahn, H.: Eisenstein series associated with {$\Gamma\sb 0(2)$}.
\newblock Ramanujan~J. \textbf{15}(2), 235--257 (2008)

\bibitem{Harnad2000}
Harnad, J., McKay, J.: Modular solutions to equations of generalized {H}alphen
  type.
\newblock Proc. Roy.~Soc. London Ser.~A \textbf{456}(1994), 261--294
(2000)

\bibitem{Heninger2006}
Heninger, N., Rains, E.M., Sloane, N.J.A.: On the integrality of {$n$}th roots
  of generating functions.
\newblock J.~Combin. Theory Ser.~A \textbf{113}(8), 1732--1745 (2006)

\bibitem{Jacobi1848}
Jacobi, C.G.J.: {\"U}ber die {D}ifferentialgleichung welcher die {R}eihen
  {$1\pm2q+2q^4\pm2q^9+{\rm etc.}$}, {$2q^{1/4}+2q^{9/4}+2q^{25/4}+{\rm etc.}$}
  {G}en{\"u}ge leisten.
\newblock J.~Reine Angew. Math. \textbf{36}, 97--112 (1848)

\bibitem{Kaneko2004}
Kaneko, M., Koike, M.: Quasimodular solutions of a differential equation of
  hypergeometric type.
\newblock In: K.~Hashimoto, K.~Miyake, H.~Nakamura (eds.) Galois Theory and
  Modular Forms, no.~11 in Dev. Math., pp. 329--336. Kluwer, Boston (2004)

\bibitem{Leopoldt58}
Leopoldt, H.W.: Eine {V}erallgemeinerung der {B}ernoullischen {Z}ahlen.
\newblock Abh. Math. Sem. Univ. Hamburg \textbf{22}, 131--140 (1958)

\bibitem{Maier12}
Maier, R.S.: On rationally parametrized modular equations.
\newblock J.~Ramanujan Math. Soc. \textbf{24}(1), 1--73 (2009)

\bibitem{Martin97}
Martin, Y., Ono, K.: Eta-quotients and elliptic curves.
\newblock Proc. Amer. Math. Soc. \textbf{125}(11), 3169--3176 (1997)

\bibitem{Milne2002}
Milne, S.C.: Infinite families of exact sums of squares formulas, {J}acobi
  elliptic functions, continued fractions, and {S}chur functions.
\newblock Ramanujan~J. \textbf{6}(1), 7--149 (2002)

\bibitem{Ohyama96}
Ohyama, Y.: Systems of nonlinear differential equations related to second order
  linear equations.
\newblock Osaka J.~Math. \textbf{33}(4), 927--949 (1996)

\bibitem{Ono95}
Ono, K., Robbins, S., Wahl, P.T.: On the representation of integers as sums of
  triangular numbers.
\newblock Aequationes Math. \textbf{50}(1--2), 73--94 (1995)

\bibitem{VanderPol51}
van~der Pol, B.: On a non-linear partial differential equation satisfied by the
  logarithm of the {J}acobian theta-functions, with arithmetical applications,
  {I},~{II}.
\newblock Indagationes Math. \textbf{13}, 261--271, 272--284 (1951)

\bibitem{vanderPol54}
van~der Pol, B.: The representation of numbers as sums of eight, sixteen and
  twenty-four squares.
\newblock Indagationes Math. \textbf{16}, 349--361 (1954)

\bibitem{Rademacher73}
Rademacher, H.: Topics in Analytic Number Theory, \emph{Die Grundlehren der
  mathematischen Wissenschaften}, vol. 169.
\newblock Springer-Verlag, New York/Berlin (1973)

\bibitem{Ramamani70}
Ramamani, V.: On some identities conjectured by {R}amanujan in his lithographed
  notes connected with partition theory and elliptic modular functions ---
  their proofs --- interconnection with various other topics in the theory of
  numbers and some generalizations.
\newblock Ph.D. thesis, University of Mysore, Mysore, India (1970)

\bibitem{Ramamani89}
Ramamani, V.: On some algebraic identities connected with {R}amanujan's work.
\newblock In: N.K. Thakare, K.C. Sharma, T.T. Raghunathan (eds.) Ramanujan
  International Symposium on Analysis (Pune, 1987), pp. 277--291. Macmillan of
  India, New Delhi (1989)

\bibitem{Rankin56}
Rankin, R.A.: The construction of automorphic forms from the derivatives of a
  given form.
\newblock J.~Indian Math. Soc. (N.S.) \textbf{20}, 103--116 (1956)

\bibitem{Rankin65}
Rankin, R.A.: Sums of squares and cusp forms.
\newblock Amer. J.~Math. \textbf{87}, 857--860 (1965)

\bibitem{Rankin76}
Rankin, R.A.: Elementary proofs of relations between {E}isenstein series.
\newblock Proc. Roy.~Soc. Edinburgh Ser.~A \textbf{76}(2), 107--117 (1976)

\bibitem{Resnikoff65}
Resnikoff, H.L.: A differential equation for the theta function.
\newblock Proc. Nat. Acad. Sci. USA \textbf{53}, 692--693 (1965)

\bibitem{Schoeneberg74}
Schoeneberg, B.: Elliptic Modular Functions, \emph{Die Grundlehren der
  mathematischen Wissenschaften}, vol. 203.
\newblock Springer-Verlag, New York/Berlin (1974)

\bibitem{Stiller88}
Stiller, P.F.: Classical automorphic forms and hypergeometric functions.
\newblock J.~Number Theory \textbf{28}(2), 219--232 (1988)

\bibitem{Takeuchi77}
Takeuchi, K.: Arithmetic triangle groups.
\newblock J.~Math. Soc. Japan \textbf{29}(1), 91--106 (1977)

\bibitem{Takeuchi77a}
Takeuchi, K.: Commensurability classes of arithmetic triangle groups.
\newblock J.~Fac. Sci. Univ. Tokyo Sect. IA Math. \textbf{24}(1), 201--212
  (1977)

\bibitem{Takhtajan92}
Takhtajan, L.A.: A simple example of modular forms as tau-functions for
  integrable equations.
\newblock Theoret. and Math. Phys. \textbf{93}(2), 1308--1317 (1992).
\newblock {R}ussian original in {T}eoret. {M}at. {F}iz. 93 (2) (1992), 330--341

\bibitem{Tricomi51}
Tricomi, F.: Funzioni Ellittiche, 2nd edn.
\newblock N.~Zanichelli, Bologna, Italy (1951)

\bibitem{Zudilin2003}
Zudilin, W.: The hypergeometric equation and {R}amanujan functions.
\newblock Ramanujan~J. \textbf{7}(4), 435--447 (2003)

\end{thebibliography}

\end{document}